\documentclass[12pt]{amsart}

\usepackage[utf8]{inputenc}

\usepackage[top=3cm, bottom=3cm, left=2.5cm, right=2.5cm]{geometry}
\usepackage[expansion=false]{microtype}

\usepackage{amsmath}
\usepackage{amssymb}
\usepackage{amsthm}
\usepackage{etoolbox}

\usepackage{float}
\usepackage{multirow}
\usepackage{subcaption}   % loads caption
\usepackage{enumitem}

\usepackage{tikz}
\usetikzlibrary{quotes, positioning, calc}
\usepackage{tikz-cd}

\usepackage{xcolor}

\usepackage{hyperref}
\hypersetup{
  colorlinks = true,
  linkcolor  = blue!55!black,
  citecolor  = green!45!black,
  urlcolor   = blue!55!black,
}

\newcommand{\Z}{\mathbb{Z}}

\newcommand{\R}{\mathbb{R}}
\newcommand{\C}{\mathbb{C}}

\newcommand{\into}{\hookrightarrow}

\DeclareRobustCommand\longtwoheadrightarrow{\relbar\joinrel\twoheadrightarrow}

\newcommand{\del}{\partial}
\newcommand{\delbar}{\bar\partial}

\newcommand{\om}{\omega}
\newcommand{\Om}{\Omega}

\newcommand{\oml}{\omega_\lambda}
\newcommand{\omlc}{\omega_{\lambda,c}}
\newcommand{\jj}{\mathcal{J}}

\newcommand{\jomlc}{\mathcal{J}_{\om_{\lambda,c}}}

\newcommand{\jsomc}{\mathcal{J}^{S^1}_{\om_{\lambda,c}}}

\newcommand{\Ii}{\mathcal{I}}

\newcommand{\iomlc}{\mathcal{I}_{\om_{\lambda,c}}}

\newcommand{\isomlc}{\mathcal{I}^{S^1}_{\om_{\lambda,c}}}

\newcommand{\Gg}{\mathcal{G}}

\newcommand{\Oo}{\mathcal{O}}
\newcommand{\bX}{\mathbb{X}}

\def\eps{\varepsilon}

\newcommand{\wt}{\widetilde}
\newcommand{\tom}{\widetilde\omega}
\newcommand{\wom}{\widetilde{\om}}
\newcommand{\wJ}{\widetilde{J}}
\newcommand{\wC}{\widetilde{\C}^2}
\newcommand{\wM}{\widetilde{M}}

\newcommand{\ttT}{\widetilde \T}

\newcommand{\CP}{\mathbb{C}P}

\newcommand{\T}{\mathbb{T}^2}
\newcommand{\SSS}{S^2 \times S^2}
\newcommand{\CCC}{\CP^2\# \overline{\CP^2}}
\newcommand{\SSSC}{(S^2 \times S^2)\# \overline{\CP^2}}
\newcommand{\fixedcomp}{\mathcal{F}}   % fixed point set component
\newcommand{\exdivisor}{Z}             % exceptional divisor of a blow-up

\newcommand{\wstd}{\omega_{\mathrm{std}}}

\DeclareMathOperator{\im}{im}
\DeclareMathOperator{\id}{id}
\DeclareMathOperator{\Inv}{Inv}
\DeclareMathOperator{\Wt}{Wt}
\DeclareMathOperator{\PD}{PD}
\DeclareMathOperator{\Emb}{Emb}
\DeclareMathOperator{\IEmb}{\Im\mathrm{Emb}}

\DeclareMathOperator{\SCont}{SCont}

\DeclareMathOperator{\Diff}{Diff}

\DeclareMathOperator{\Symp}{Symp}
\DeclareMathOperator{\Ext}{Ext}

\DeclareMathOperator{\Stab}{Stab}
\DeclareMathOperator{\Fix}{Fix}

\DeclareMathOperator{\Aut}{Aut}

\DeclareMathOperator{\Sp}{Sp}
\DeclareMathOperator{\U}{U}
\DeclareMathOperator{\SO}{SO}
\DeclareMathOperator{\PU}{PU}

\DeclareMathOperator{\pr}{pr}
\DeclareMathOperator{\near}{near}

\DeclareMathOperator{\Gauge1}{Gauge}
\DeclareMathOperator{\Isom}{Isom}

\newcommand{\EmbG}{{\Emb^G(M,\fixedcomp, c)}}
\newcommand{\IEmbG}{{\IEmb^G(M,\fixedcomp, c)}}
\newcommand{\EmbS}{{\Emb^{S^1}(M,\fixedcomp, c)}}

\newcommand{\EmbGiota}{{\Emb^G_\iota(M,\fixedcomp, c)}}

\theoremstyle{plain}
\newtheorem{thm}{Theorem}[section]
\newtheorem{prop}[thm]{Proposition}
\newtheorem{cor}[thm]{Corollary}
\newtheorem{lemma}[thm]{Lemma}

\theoremstyle{definition}
\newtheorem{defn}[thm]{Definition}
\newtheorem{exmp}[thm]{Example}

\newtheorem{Setting}[thm]{Setting}

\theoremstyle{remark}
\newtheorem{Remark}[thm]{Remark}
\newtheorem{Remarks}[thm]{Remarks}
\newtheorem*{remark}{Remark}

\newcommand{\eoesymbol}{$\between$}
\DeclareRobustCommand{\eoe}{%
  \ifmmode \mathqed
  \else
    \leavevmode\unskip\penalty9999 \hbox{}\nobreak\hfill
    \quad\hbox{\eoesymbol}%
  \fi
}

\AtEndEnvironment{Remark}{\eoe}
\AtEndEnvironment{Remarks}{\eoe}
\AtEndEnvironment{remark}{\eoe}
\AtEndEnvironment{exmp}{\eoe}

\begin{document}
\title[Homotopy type of equivariant embeddings]{Homotopy type of equivariant embeddings}
\date{November 12, 2026}

\author{Pranav Chakravarthy}
\address{Department of Mathematics\\Universit\'e Libre de Bruxelles \\ Brussels, Belgium}
\email{pranav.vijay.chakravarthy@ulb.be}

\author{River Chiang}
\address{Department of Mathematics, National Cheng Kung University, Tainan 701, Taiwan}
\email{riverch@mail.ncku.edu.tw}

\author{Liat Kessler}
\address{Department of Mathematics, Physics, and Computer Science, 
University of Haifa, Tivon 36006, Israel}
\email{lkessler@math.haifa.ac.il}

\author[M. Pinsonnault]{Martin Pinsonnault}
\address{Department of Mathematics \\  University of Western Ontario \\ Canada}
\email{mpinson@uwo.ca}



\begin{abstract}
We study spaces of symplectic embeddings of a closed four-ball into a closed
symplectic four-manifold that are equivariant under a compact group action and
centred on a component of the fixed point set. Adapting McDuff's argument, we
reduce the path-connectedness of these spaces, and the uniqueness of equivariant
blow-ups up to isotopy, to an equivariant deformation-to-isotopy statement on
the blow-up. For compact abelian groups, we identify the stabilizer of a ball
with the group of equivariant symplectomorphisms of the blow-up that preserve
the exceptional divisor. For Hamiltonian circle actions on $\mathbb{C}P^2$,
$S^2\times S^2$ and $\mathbb{C}P^2\#\overline{\mathbb{C}P^2}$, we obtain the
invariant curves needed for inflation from toric extensions of the blown-up
action, and we prove that the embedding spaces are path-connected. When the ball
is centred at an isolated fixed point of $S^2\times S^2$ or
$\mathbb{C}P^2\#\overline{\mathbb{C}P^2}$, the space of unparametrized balls is
contractible or homotopy equivalent to $\Omega S^3\times S^1$, according to the
numbers of toric actions that extend the circle action and its blow-up. On
$\mathbb{C}P^2$, it is homotopy equivalent to the fixed component. Along the
way, we compute the centralizers of Hamiltonian circle actions on the one-point
blow-up of $S^2\times S^2$, whose homotopy type depends only on the number of
toric extensions.
\end{abstract}

\maketitle

\tableofcontents

%%%%%%%%%%%%%%%%%%%%%%%%%%%%%%%%%%%%%%%%%%%%%%%%%%%%%%%%%%%%%%%%%%%%%%%%%%%%%%%%
%%%%%%%%%%%%%%%%%%%%%%%%%%%%%%%%%%%%%%%%%%%%%%%%%%%%%%%%%%%%%%%%%%%%%%%%%%%%%%%%
%%%%%%%%%%%%%%%%%%%%%%%%%%%%%%%%%%%%%%%%%%%%%%%%%%%%%%%%%%%%%%%%%%%%%%%%%%%%%%%%

\section{Introduction}\label{sec:intro}

The study of spaces of symplectic embeddings of balls is a central theme of
symplectic topology. Deciding whether such a space is non-empty is the
symplectic packing problem~\cite{gromovcurves,MP,Traynor,Biran}, which can be seen as an extension of the nonsqueezing theorem, and which generalizes the notion of capacity.  Deciding whether this space is path-connected is the famous symplectic camel problem~\cite{McD-Traynor}.

Although very difficult in general, these problems are more accessible on four-dimensional uniruled manifolds for which
Taubes--Seiberg--Witten theory~\cite{Taubes-SWGr} together with positivity of intersections gives a precise control of $J$-holomorphic curves. For the projective plane, for the rational ruled surfaces and for their blow-ups, the packing and the connectedness problems are well understood. The packing of $\CP^2$ by equal balls was settled in~\cite{MP,Biran}, and in general the packing problem for these manifolds reduces to the description of the symplectic cone of their blow-ups, which is given in terms of exceptional classes in~\cite{liliu}. Path-connectedness of the space of embeddings of one ball was established by McDuff for $\CP^2$ and for rational ruled surfaces~\cite{McD:isotopy}, then by Biran for several equal balls in $\CP^2$~\cite{BiranC}, and for any number of balls in these manifolds and their blow-ups in~\cite[Corollary~1.5]{deftoiso}. An important corollary is that symplectic blow-ups of these manifolds are unique up to isotopy. Other homotopical properties of the embedding spaces, including in some cases their full homotopy type, were obtained in~\cite{LP,P-com,Anj-Lal-Pin,AKP,AKP26}. For instance, the space of unparametrized embeddings of a ball into $\CP^2$ is homotopy equivalent to $\CP^2$ for every admissible capacity, whereas on $\SSS$ endowed with a split form $\mu\sigma\oplus\sigma$, $\mu>1$ not an integer, there is a critical capacity below which the space of unparametrized balls is homotopy equivalent to $\SSS$ itself, and above which it does not have the homotopy type of a finite CW-complex~\cite{P-com,Anj-Lal-Pin}.

This body of work shows that the existence of packings and the path-connectedness of embedding spaces can be decided by means of $J$-holomorphic curves that exist for \emph{generic} tamed or compatible almost complex structures, and whose existence is guaranteed by non-vanishing Gromov--Taubes invariants. Finer homotopy invariants of embedding spaces, on the other hand, also depend on curves that only exist for \emph{non-generic} almost complex structures, which typically form strata of positive codimension in the space of almost complex structures. On rational $4$-manifolds, the resulting stratification depends only on the cohomology class of the symplectic form.

In this paper we study spaces of symplectic embeddings of balls in the presence of group actions. Since invariant almost complex structures form a subspace of infinite codimension in the space of all almost complex structures, they cannot be assumed to be generic in the usual Gromov--Taubes sense. One of the main features of the equivariant setting is therefore the need to work with non-generic structures whose $J$-holomorphic curves are not determined by Gromov--Taubes invariants alone. This is already the case for the equivariant versions of the packing and connectedness problems. The space of invariant almost complex structures inherits a partition into strata, which depends on the action and not only on the cohomology class of the symplectic form.
Another feature, absent in the non-equivariant setting, is the existence of local invariants of the action, such as its isotropy representation, on which the homotopy type of the group of equivariant reparametrizations of the ball depends. Together, these two features make the equivariant problem genuinely different from the non-equivariant one, and substantial parts of the non-equivariant framework have to be reformulated.

\subsection{The equivariant framework}
Let $(M,\omega)$ be a closed symplectic $4$-manifold on which a compact Lie group $G$ acts symplectically, and let $\fixedcomp$ be a connected component of the fixed point set $M^G$. Let $B^4(c) \subset \R^4$ be the closed standard ball of capacity $c$, that is of radius $\sqrt{c/\pi}$, endowed with the restriction of the standard symplectic form $\wstd$ and with a linear $G$-action. We denote by
\[
\Emb^{G}(M,\fixedcomp,c)
\]
the space of $G$-equivariant symplectic embeddings $(B^4(c),0)\into (M,\fixedcomp)$, equipped with the $C^\infty$ topology. Here a
symplectic embedding of the closed ball $B^4(c)$ is the restriction of a symplectic embedding of an open neighbourhood of $B^4(c)$.
The space of \emph{unparametrized} equivariant embeddings centred on $\fixedcomp$ is the quotient
\[
\IEmb^G(M,\fixedcomp,c):= \Emb^G(M,\fixedcomp,c)/ \Symp^G(B^{4}(c)),
\]
where the group $\Symp^G(B^{4}(c))$ of equivariant symplectomorphisms of the closed ball acts by reparametrization. Equivariance forces the linear action on the ball to be equivalent to the isotropy representation at the centre.

The first question is to determine the capacities for which $\EmbG$ is non-empty. Their supremum is the \emph{$G$-equivariant capacity} of $\fixedcomp$ (Definition~\ref{defn:capacity}). For toric actions, it can be read off the moment polygon, and for Hamiltonian circle actions from their Karshon graphs. On rational ruled surfaces, it is easy to see that the $S^1$-equivariant capacities can also be read off the moment polygons of their toric extensions (Lemma~\ref{lemma:capacity basics} and Corollary~\ref{cor:capacity is max}).

The second question is whether $\EmbG$ is path-connected. One of our motivations is the uniqueness of equivariant symplectic blow-ups. Recall that a symplectic blow-up of $(M,\omega)$ is obtained by removing the interior of a symplectically embedded closed ball $B^4(c)$ and collapsing its boundary along the Hopf fibration $S^3\to\CP^1$, see~\cite{GS3}, \cite{Lerman-cut}, \cite[\S 7.1]{McD-Salamon-Intro}, or Section~\ref{sec:blowup}. When the embedding is $G$-equivariant, the blow-up inherits a $G$-action. Unlike the complex blow-up, the symplectic blow-up is not canonical. It depends on the embedded ball and on auxiliary choices, which themselves depend on local properties of the action. We say that equivariant symplectic blow-ups of $(M,\omega)$ are unique up to equivariant isotopy if, for any two embeddings $\iota_0,\iota_1\in\EmbG$ centered at the same point, the corresponding blow-up forms $\widetilde{\omega}_0$ and $\widetilde{\omega}_1$ on the blown-up manifold are joined by a path of cohomologous invariant symplectic forms. For Hamiltonian toric and circle actions on $4$-manifolds, the classifications of Delzant~\cite{delzant} and Karshon~\cite{karshon} show that equivariant blow-ups at a given fixed component and of a given size are unique up to equivariant symplectomorphism (Remark~\ref{rem:ball independence}). Since isotopic forms on a closed manifold are diffeomorphic (by Moser's argument), uniqueness up to isotopy is a strictly stronger classification than just classification up to equivariant symplectomorphisms. The non-equivariant uniqueness up to isotopy does not imply the equivariant uniqueness as averaging the forms of a non-equivariant isotopy over $G$ does not produce an equivariant isotopy, since an average of symplectic forms need not be non-degenerate.

Our first result reduces both the connectedness of the embedding space $\EmbG$ and the uniqueness of equivariant blow-ups up to isotopy to an equivariant ``deformation-to-isotopy'' statement for families of blow-up forms
(Propositions~\ref{prop:deformation from embeddings} and~\ref{prop:isotopy implies connectedness}). Its proof adapts the argument of
McDuff~\cite{deftoiso} to the equivariant setting and is valid for every compact subgroup $G\subset\U(2)$. The key point is that the centralizer of $G$ in $\Sp(4)$ is connected for all such subgroups (Lemma~\ref{Lemma:Sp^g(4)isconnected}). 

Deformation-to-isotopy statements are usually proved by inflation, which requires embedded $J$-holomorphic curves in specific homology classes. In the equivariant setting, these curves must be invariant and pseudo-holomorphic for an invariant almost complex structure. Since Gromov--Taubes invariants no longer guarantee that such curves exist, they must instead be found by an analysis specific to the action. For toric and Hamiltonian circle actions on $4$-manifolds, the classifications of Delzant and Karshon again provide the necessary input (Sections~\ref{sec:toric} and~\ref{sec:deformation to isotopy}).

Once the space $\EmbG$ is known to be path-connected, its homotopy type, and that of $\IEmbG$, can be approached as in the non-equivariant case~\cite{LP,P-com}. Let $\Symp^G(M, \fixedcomp,\om)$ be the group of $G$-equivariant symplectomorphisms of $(M,\om)$ that send $\fixedcomp$ to itself and let $\Symp^G(M,\om)$ be the symplectic centralizer of $G$. By an equivariant version of the symplectic isotopy of balls~\cite[Theorem 3.3.1]{McD-Salamon-Intro}, the group $\Symp^G(M, \fixedcomp,\om)$ acts transitively on $\EmbG$ and there is a ladder of fibrations
\begin{equation}\label{eq:MainFibration}
\begin{tikzcd}
\Fix^G(\iota) \arrow[hookrightarrow]{r} \arrow[hookrightarrow]{d}&\Symp^G(M, \fixedcomp,\om) \arrow{r} \arrow[equal]{d}& \EmbG\arrow[twoheadrightarrow]{d}\\
\Stab^G(\iota) \arrow[hookrightarrow]{r}&\Symp^G(M, \fixedcomp,\om) \arrow[twoheadrightarrow]{r}& \IEmbG
\end{tikzcd}
\end{equation}
where $\iota:(B^4(c),0) \into (M,\fixedcomp)$ is a given equivariant symplectic embedding, $\Fix^G(\iota)$ is the subgroup of equivariant symplectomorphisms that restrict to the identity on the image $\iota(B^4(c))$, and $\Stab^G(\iota)$ is the subgroup of equivariant symplectomorphisms that take $\iota(B^4(c))$ to itself. When $G$ is abelian, the stabilizer $\Stab^G(\iota)$ is homotopy equivalent to the group of equivariant symplectomorphisms of the blow-up that preserve the exceptional divisor (Proposition~\ref{prop:homotopy equivalence Stab(B) and Stab(E)}), so that the computation of the homotopy type of $\IEmbG$ essentially reduces to that of two centralizers, one before and one after the blow-up.

This identification, as well as the analysis of centralizers, relies on an equivariant version of Gromov's theorem stating that the group of compactly supported symplectomorphisms of the open $4$-ball is contractible~\cite{gromovcurves}. At the time of writing, this equivariant version is known only for compact abelian groups (Proposition~\ref{prop:equivariant Gromov}, see~\cite[Theorem~3.29]{ChPin-Memoirs}). For such groups, it implies that the group $\Symp^G(B^4(c))$ of equivariant reparametrizations deformation retracts onto the centralizer $\U^G(2)$ of $G$ in $\U(2)$ (Proposition~\ref{prop:homotopy type equivariant symplectomorphism of the ball}). This is another reason to restrict our analysis to Hamiltonian circle and toric actions.

%%%%%%%%%%%%%%%%%%%%%%%%%%%%%%%%%%%%%%%%%%%%%%%%%%%%%%%%%%%%%%%%%%%%%%%%%%%%%%%%
\subsection{Main results}
\label{subsec:intro-results}
%%%%%%%%%%%%%%%%%%%%%%%%%%%%%%%%%%%%%%%%%%%%%%%%%%%%%%%%%%%%%%%%%%%%%%%%%%%%%%%%

We carry out this programme for toric four-manifolds, and for Hamiltonian circle actions on the projective plane and on the rational ruled surfaces $\SSS$ and $\CCC$. Symplectic forms on these manifolds are classified, up to rescaling and diffeomorphism, in~\cite{Lalonde-McDuff}.

\subsubsection{Toric actions.}
As a first illustration of the method, Section~\ref{sec:toric} gives a new proof, in dimension four, of Pelayo's theorem~\cite{pelayo} that the space of $\T$-equivariant symplectic embeddings of a ball centred at a fixed point of a toric four-manifold is homotopy equivalent to the torus $\T$ (Proposition~\ref{prop:Pelayo}). There, the invariant curves needed for the deformation-to-isotopy lemma are the spheres over the edges of the moment polygon.
\\

\subsubsection{The projective plane.} Equivariant embeddings with respect to Hamiltonian circle actions on $\CP^2$ are easy to unederstand.

\begin{thm}\label{thm:main-a}
Let $S^1$ act effectively and symplectically on
$(\CP^2,\lambda \omega_{\operatorname{FS}})$, where $\omega_{\operatorname{FS}}$
gives a line area one, let $\fixedcomp$ be a connected component of the fixed
point set, and let $0<c<\lambda$. Then $\Emb^{S^1}(\CP^2,\fixedcomp,c)$ is
non-empty and path-connected, and
\[
\IEmb^{S^1}(\CP^2,\fixedcomp,c)\simeq\fixedcomp .
\]
In other words, the space of unparametrized invariant balls is contractible
when $\fixedcomp$ is an isolated fixed point, and homotopy equivalent to $S^2$
when $\fixedcomp$ is a pointwise fixed line.
\end{thm}

Theorem~\ref{thm:main-a} is the equivariant counterpart of the fact that the
space of unparametrized balls in $\CP^2$ is homotopy equivalent to $\CP^2$ for
every capacity~\cite{P-com}. In both cases the space of balls retracts onto the
space of their possible centres.

\subsubsection{Rational ruled surfaces.}
On $\SSS$ and $\CCC$ we normalize the symplectic form as
in~\S\ref{subsec:ruled}: the fibre has area $1$, and a section of
self-intersection $0$, respectively $+1$, has area $\lambda$.

\begin{thm}\label{thm:main-b}
Let $(M,\omega)$ be $\SSS$ or $\CCC$ with a symplectic form normalized as above,
let $S^1$ act effectively and symplectically on $(M,\omega)$, let $\fixedcomp$
be a connected component of its fixed point set, and let
$0<c<c^{S^1}(M,\fixedcomp)$.
\begin{enumerate}
\item The space $\EmbS$ is path-connected.
\item Assume that $\fixedcomp=\{p\}$ is an isolated fixed point. Let $N$ and
$\widehat N$ be the numbers of inequivalent toric actions extending,
respectively, the circle action on $M$ and its equivariant blow-up at a ball of
capacity $c$ centred at $p$. Then $N,\widehat N\leq 2$, and
\[
\IEmb^{S^1}(M,p,c)\simeq
\begin{cases}
\Omega S^3\times S^1 & \text{if } N=2 \text{ and } \widehat N=1,\\
\{*\} & \text{otherwise}.
\end{cases}
\]
\end{enumerate}
\end{thm}

The number of extensions $N$ and $\widehat N$ can be explicitly determined in terms of the weights of the action, the parameters $\lambda$, the capacity $c$, and  the position of $p$, see Table~\ref{table:NumberOfStrata} and~\S\ref{section:reading}. The blown-up action has at most as many toric extensions as the original one, that is $\widehat N\leq N$, except in the \emph{decomposable case} in which $M=\CCC$ and $c>\lambda-1$. In this regime, the capacity of the ball exceeds the area $\lambda-1$ of the section of self-intersection $-1$, so that the class of the exceptional divisor becomes decomposable. The stabilizer of a ball then need not be the full centralizer of the blow-up, but our framework still shows that the space of unparametrized balls is contractible (see~\S\ref{subsubsec:decomposable case}).

It is interesting to compare Theorem~\ref{thm:main-b} with the non-equivariant results of~\cite{P-com,Anj-Lal-Pin}. In both settings, the homotopy type changes when the capacity crosses a critical value, beyond which the space of balls no longer has the homotopy type of a finite CW-complex (the homology of the factor $\Omega S^3$ is non-trivial in every even degree). Non-equivariantly, the threshold is the capacity at which the stratification of the space of almost complex structures on the blow-up loses a stratum~\cite{P-com}. Equivariantly, it is the capacity at which a toric extension is lost, which is also equivalent to the loss of a stratum in the space of invariant almost complex structures. Contrary to the non-equivariant case, the equivariant threshold depends on the fixed point at which the ball is centred,
see the example at the end of~\S\ref{section:reading}.

\subsubsection{Uniqueness of equivariant blow-ups.}
As in the non-equivariant case, we obtain the first assertions of
Theorems~\ref{thm:main-a} and~\ref{thm:main-b} by producing equivariant
isotopies between blow-up forms, see Section~\ref{sec:steps} and
Section~\ref{sec:deformation to isotopy}. The proof therefore also gives the
following result.

\begin{thm}\label{thm:isotopy}
Let $(M,\omega)$ be either $(\CP^2,\lambda \omega_{\operatorname{FS}})$ or a
ruled symplectic four-manifold over $S^2$, endowed with an effective
symplectic circle action. Then $S^1$-equivariant symplectic blow-ups of $(M,\omega)$ of a given
size at a given fixed point component are unique up to equivariant symplectic
isotopy.
\end{thm}

\subsubsection{Centralizers.}
The computation of the stabilizers of balls rests on the following description
of centralizers, which extends to the one-point blow-up $\bX_2$ the results
of~\cite{ChPin-Memoirs} for $\SSS$ and $\CCC$. We write
$\Symp^{S^1}_h(\bX_2,\omlc)$ for the subgroup of equivariant symplectomorphisms
acting trivially on homology, the forms $\omlc$ being those
of~\S\ref{subsec:blowup}.

\begin{thm}\label{thm:main-centralizers}
\begin{enumerate}
\item Let $S^1$ act effectively and symplectically on $(\CP^2,\omega)$, where
$\omega$ is a Fubini--Study form. Then the centralizer
$\Symp^{S^1}(\CP^2,\omega)$ is homotopy equivalent to $\U(2)$ if the action
fixes a line pointwise, and to the torus $\T$ otherwise.
\item Let $S^1$ act effectively and symplectically on $(\bX_2,\omlc)$. The action extends to at
most two inequivalent toric actions. If it extends to a single toric action
$\ttT$, then $\Symp^{S^1}_h(\bX_2,\omlc)$ is homotopy equivalent to $\ttT$.
If it extends to two toric actions $\ttT_i$ and
$\ttT_j$, then $\Symp^{S^1}_h(\bX_2,\omlc)$ is homotopy equivalent to the
homotopy pushout, in the category of topological groups, of the inclusions
$\ttT_i\hookleftarrow S^1\hookrightarrow\ttT_j$, which has the homotopy type of
$\Omega S^3\times\mathbb{T}^3$.
\end{enumerate}
\end{thm}

Non-equivariantly, the homotopy type of the symplectomorphism group of $\bX_2$ is governed by a stratification of the space of almost complex structures whose number of strata grows with $\lambda$~\cite{P-com}. Equivariantly, at most two strata are met, which is why only two homotopy types occur. The circle actions on $\CP^2$, $\SSS$ and $\CCC$ that fix a line or a section pointwise have centralizers homotopy equivalent to $\U(2)$ or to $S^1\times\SO(3)$, by part~(1) of Theorem~\ref{thm:main-centralizers} and~\cite{ChPin-Memoirs}. On $\bX_2$, however, a curve of negative self-intersection fixed pointwise by the circle meets an invariant exceptional curve, and the centralizer of an action with a single toric extension is always a torus. By contrast with toric actions, whose centralizer is always homotopy equivalent to the torus itself~\cite{Pin-MaxTori}, the centralizer of a circle action thus depends on the symplectic form, but only through the number of its toric extensions.

\subsubsection{Questions left open.}
On the manifolds considered in this paper, the only case left open is that of
circle actions with non-isolated fixed points. When $\fixedcomp$ is a fixed
surface in $\SSS$ or $\CCC$, Theorem~\ref{thm:main-b} shows that $\EmbS$ is
path-connected, but we do not determine its homotopy type.

Beyond these manifolds, the case of circle actions on general symplectic
$4$-manifolds, rational or irrational, is open. Our analysis of invariant almost
complex structures relies on the fact that every Hamiltonian circle action on
the manifolds considered here extends to a toric action~\cite{finTori}.
This is no longer true on blow-ups of $\CP^2$ at three or more points, nor on irrational ruled surfaces, and it is not clear how much of the analysis carries over to maximal circle actions.

More generally, the local results (Propositions~\ref{prop:homotopy type equivariant symplectomorphism of the ball} and~\ref{prop:homotopy equivalence Stab(B) and Stab(E)}) hold for compact abelian groups, possibly finite, and those of Section~\ref{sec:steps} for every compact subgroup of $\U(2)$, whereas the global analysis of invariant almost complex structures is specific to toric and circle actions. Natural next steps are the extension of this analysis to finite cyclic groups, whose centralizers on rational ruled surfaces were computed in~\cite{ChPin-finitecyclic}, see also~\cite{LiatRiver,LiatRiver_noExtension}, and to several invariant balls. One may also hope to relate the topology of these embedding spaces to other invariants of the action, such as its equivariant cohomology.

%%%%%%%%%%%%%%%%%%%%%%%%%%%%%%%%%%%%%%%%%%%%%%%%%%%%%%%%%%%%%%%%%%%%%%%%%%%%%%%%
\subsection{Organization of the paper}
\label{subsec:intro-organization}
%%%%%%%%%%%%%%%%%%%%%%%%%%%%%%%%%%%%%%%%%%%%%%%%%%%%%%%%%%%%%%%%%%%%%%%%%%%%%%%%

Section~\ref{section:equivariant embeddings} is devoted to the local models. It shows that the
reparametrization group $\Symp^G(B^4(c))$ retracts onto $\U^G(2)$
(Proposition~\ref{prop:homotopy type equivariant symplectomorphism of the ball})
and identifies the stabilizer of an invariant ball with that of the exceptional
divisor of the blow-up
(Proposition~\ref{prop:homotopy equivalence Stab(B) and Stab(E)}).

Section~\ref{sec:steps} reduces path-connectedness to a deformation-to-isotopy
statement. Embeddings are first normalized near their centre
(Proposition~\ref{Lemma:NormalisationofEmbeddings}), two normalized embeddings
then give rise to a deformation of invariant blow-up forms
(Proposition~\ref{prop:deformation from embeddings}), and an isotopy homotopic
to that deformation yields a path of
embeddings (Proposition~\ref{prop:isotopy implies connectedness}).
Subsection~\ref{section:role of 3.12 and 3.13} states precisely what remains to
be proved.

Section~\ref{sec:toric} carries out the programme for toric actions and
recovers Pelayo's theorem in dimension four.

Section~\ref{section:circle actions} collects what is needed about Hamiltonian
circle actions on $\SSS$, $\CCC$ and their one-point blow-up $\bX_2$: toric
extensions and the different presentations of a circle action, the equivariant
capacity of a fixed component and its relation with toric extensions
(Proposition~\ref{prop:capacity rule}), and the identification of stabilizers of
balls with centralizers (Theorem~\ref{thm:stabilizer}).

Section~\ref{sec:Centralizers on bX2} computes the centralizers of Hamiltonian circle actions on $\bX_2$:
\S\ref{section:Strataofinvariantacs} relates the strata met by the invariant
almost complex structures to toric extensions, \S\ref{section:codimension}
computes their codimension and compares the two notions of genericity, and
\S\ref{section:homotopy type SSSC} computes the homotopy types.

Section~\ref{section:Centralizers on CP2} computes the homotopy type of equivariant symplectomorphisms for Hamiltonian circle actions on $\CP^2$.

Section~\ref{sec:deformation to isotopy} proves the deformation-to-isotopy
statement for circle actions (Proposition~\ref{prop:special}), which completes the proof of
path-connectedness (Theorem~\ref{thm:main0}). 

Section~\ref{section:Homotopy types of embedding spaces} combines these results
to determine the homotopy types of the embedding spaces, for $\CP^2$
in~\S\ref{subsec:CP2 embeddings} and for the rational ruled surfaces
in~\S\ref{subsec:ruled embeddings}, and explains in~\S\ref{section:reading} how
to read them off the data.

Finally, the first appendix computes the equivariant gauge groups used in the local models,
and the second one lists the decorated graphs of all the Hamiltonian circle actions on $\bX_2$.

\begin{Remark}\label{rem:CW}
All the spaces occurring in this paper are
metrizable manifolds modelled on a locally convex spaces, hence absolute neighbourhood retracts by~\cite[Theorem 5]{Palais-Homotopy}, and therefore has the homotopy type of a CW complex.
In particular a weak homotopy equivalence between two of them is a homotopy equivalence, and we use the two notions interchangeably.
\end{Remark}

%%%%%%%%%%%%%%%%%%%%%%%%%%%%%%%%%%%%%%%%%%%%%%%%%%%%%%%%%%%%%%%%%%%%%%%%%%%%%%%%
\section{Equivariant embeddings in dimension $4$}\label{section:equivariant embeddings}
%%%%%%%%%%%%%%%%%%%%%%%%%%%%%%%%%%%%%%%%%%%%%%%%%%%%%%%%%%%%%%%%%%%%%%%%%%%%%%%%

%%%%%%%%%%%%%%%%%%%%%%%%%%%%%%%%%%%%%%%%%%%%%%%%%%%%%%%%%%%%%%%%%%%%%%%%%%%%%%%%
\subsection{Stabilizers of balls and stabilizers of exceptional curves}
%%%%%%%%%%%%%%%%%%%%%%%%%%%%%%%%%%%%%%%%%%%%%%%%%%%%%%%%%%%%%%%%%%%%%%%%%%%%%%%%

Let $(M,\om)$ be a compact symplectic $4$-manifold on which a compact group $G$ is acting effectively through symplectomorphisms. Suppose $G$ is acting unitarily on the standard ball $B^4(c)$. Let $(\tilde{M},\tom)$ denote the equivariant symplectic blow-up of $M$ along a fixed equivariant symplectic embedding $\iota:B^4(c)\into M$. Let $\exdivisor$ denote the exceptional divisor inside $\tilde{M}$ for this blow-up. Let $\Symp^G(\tilde{M},\exdivisor,\tom)$ denote the group of equivariant symplectomorphisms of the blow-up that take the  embedded symplectic $2$-sphere $\exdivisor$ to itself. In this section we show that the group $\Stab^G(\iota)$ of equivariant symplectomorphisms preserving the image of $\iota$ is homotopy equivalent to $\Symp^G(\tilde{M},\exdivisor,\tom)$, see Proposition~\ref{prop:homotopy equivalence Stab(B) and Stab(E)}.

The following standard fact will be used several times to upgrade a homotopy equivalence to a deformation retraction.

\begin{lemma}\label{lemma:cofibration to deformation retraction}
Let $M$ be a (Hausdorff, paracompact) Fréchet manifold and let $N$ be a closed submanifold. Suppose the inclusion $N\into M$ is a homotopy equivalence. Then $M$ deformation retracts onto $N$.
\end{lemma}
\begin{proof}
By Palais~\cite[Theorem 5]{Palais-Homotopy} both $M$ and $N$ are ANRs. By~Aguilar \& al.~\cite[Theorem~4.2.13]{A-G-P}, the inclusion is a cofibration with closed image. By Hatcher~\cite[Corollary~0.20]{Hatcher-Book}, $N$ is a deformation retract of $M$.
\end{proof}

Here and below, a superscript $G$ on a group of transformations denotes the subgroup of $G$-equivariant elements.
For $G\subset\U(2)$, we write $\overline G$ for its image in $\PU(2)=\SO(3)$, which acts on $\CP^1=S^2$. A subscript $0$ denotes the identity component of a topological group. By Remark~\ref{rem:CW}, a weak homotopy equivalence between the groups considered in this section is a homotopy equivalence, and we use this without further mention.

\begin{lemma}\label{lemma:linear centralizers}
Let $K$ be a closed subgroup of $\SO(3)$, acting on the standard two-sphere $(S^2,\sigma)$. The group $\Symp^K(S^2,\sigma)$ of equivariant symplectomorphisms deformation retracts onto its subgroup $\SO^K(3)=\PU^K(2)$ of equivariant unitary maps.
\end{lemma}
\begin{proof}
M.-T. Wang~\cite{Wang-2001, Wang-2008} constructs a deformation retraction of $\Symp(S^2,\sigma)$ onto its subgroup $\SO(3)$ of orthogonal transformations by running the mean curvature flow, in $(S^2,\sigma)\times(S^2,-\sigma)$ with the product round metric, on the Lagrangian graph of a symplectomorphism: the flow exists for all time, remains a graph, depends continuously on its initial datum, and converges to a totally geodesic graph, that is, to the graph of an isometry.

The construction is equivariant. The diagonal action of $K\subset\SO(3)$ preserves the product round metric, and the graph of an equivariant $f$ is invariant. The mean curvature flow is determined by the ambient metric alone, hence commutes with isometries. By uniqueness of the flow, the evolution of an invariant submanifold stays invariant, and the maps of which it is the graph are equivariant. Their limit is an equivariant isometry, that is, an element of $\SO^K(3)$. Finally the graph of an isometry is totally geodesic, hence stationary for the flow, so the retraction fixes $\SO^K(3)$ pointwise.
\end{proof}

\begin{Remark}
The centralizer $\SO^K(3)=\PU^K(2)$ may not be connected: if $K\cong \Z_2$ is generated by a rotation by a half-turn, its centralizer $\PU^K(2)$ is $\mathrm{O}(2)$.
\end{Remark}

\begin{lemma}\label{Lemma:Sp^g(4)isconnected}
Let $G\subset\U(2)\subset\Sp(4)$ be a closed subgroup, and let $\Sp^G(4)$ and $\U^G(2)$ denote the centralizers of $G$ in $\Sp(4)$ and in $\U(2)$. Then $\Sp^G(4)$ deformation retracts onto $\U^G(2)$, and
\[
\U^G(2)=
\begin{cases}
\U(2) & \text{if~} G\subset \mathcal{Z}(\U(2)),\\
\U(L_1)\times\U(L_2) & \text{if~$G$ is abelian and~} G\not\subset \mathcal{Z}(\U(2)),\\
\mathcal{Z}(\U(2)) & \text{if~$G$ is not abelian,}
\end{cases}
\]
where, in the second case, $\C^2=L_1\oplus L_2$ is the unique decomposition into $G$-invariant complex lines on which $G$ acts with distinct weights and $\mathcal{Z}(\U(2))$ stands for the center of $\U(2)$ . In particular, $\Sp^G(4)$ is connected.
\end{lemma}
\begin{proof}
The polar decomposition defines a deformation retraction of $\Sp(4)$ onto $\U(2)$ that commutes with conjugation by elements of $\U(2)$. Since $G\subset\U(2)$, it restricts to a deformation retraction of $\Sp^G(4)$ onto $\U^G(2)$.

The group $\U^G(2)$ consists of the unitary $G$-equivariant automorphisms of the representation of $G$ on $\C^2$. If this representation is irreducible, Schur's lemma implies that every $G$-equivariant endomorphism of $\C^2$ is a complex multiple of the identity, so that $\U^G(2)$ is the group $\mathcal{Z}(\U(2))$ of unitary scalar matrices. If the representation is reducible, let $L_1$ be a $G$-invariant complex line, and let $L_2$ be its Hermitian orthogonal complement (which is also invariant). Then $G$ acts on $L_j$ through a character $\chi_j\colon G\to\U(1)$, so that $G$ is a subgroup of the torus $\U(L_1)\times\U(L_2)$ and is, therefore, abelian. If $\chi_1=\chi_2$, then $G$ acts by scalars and $\U^G(2)=\U(2)$. Otherwise, $L_1$ and $L_2$ are preserved by every $G$-equivariant endomorphism, so that $\U^G(2)=\U(L_1)\times\U(L_2)$.

Every abelian subgroup of $\U(2)$ acts reducibly on $\C^2$. Hence $G$ acts irreducibly if and only if it is not abelian. Finally, as $\U(2)$, the torus $\U(L_1)\times\U(L_2)$ and the circle $\mathcal{Z}(\U(2))$ are connected, so is $\Sp^G(4)$.
\end{proof}

Let $S\subset M$ be either a compact invariant symplectic submanifold contained in the interior of $M$, or a union of connected components of $\del M$. Inside $\Symp^G(M)$ we consider the three subgroups
\[
\Symp^G_{\id}(M,S)\ \supset\ \Symp^G_{\id,N}(M,S)\ \supset\ \Symp^G_{\id,\near}(M,S),
\]
consisting, from left to right, of the equivariant symplectomorphisms $\phi$ that restrict to the identity on $S$, of those that moreover satisfy $d\phi|_{T_SM}=\id$, the derivative being the identity on the whole restriction of $TM$ to $S$, and of those that are the identity on a neighbourhood of $S$. Statements of the following kind are standard, but we could not find a symplectic version that covers the equivariant setting and the boundary components variant at once. We therefore give a complete proof.

\begin{lemma}\label{lemma:HomotopyEquivalenceSympN-SympNear}
Let $G$ be a compact Lie group acting symplectically on a compact symplectic manifold $(M,\om)$, and let $S$ be as above. Then the inclusion
\[
\Symp^G_{\id,\near}(M,S)\into\Symp^G_{\id,N}(M,S)
\]
is a weak homotopy equivalence. 
\end{lemma}

\begin{proof}
\emph{Euler-like vector fields.} If $S$ lies in the interior of $M$, let $\mathcal N\to S$ be its normal bundle, and let $E$ be the Euler vector field of $\mathcal N$ which generates the fibrewise dilations. Recall from~\cite{BLM-Splitting} that a vector field $X$ on $M$ is \emph{Euler-like} along $S$ if it vanishes on $S$ and its linear approximation along $S$ is $E$. Such a vector field determines a tubular neighbourhood embedding of a neighbourhood of the zero section of $\mathcal N$ into $M$ that intertwines $E$
and $X$, and this embedding is equivariant when $X$ is invariant~\cite{BLM-Splitting}. We fix a $G$-invariant Euler-like vector field, for instance the image of $E$ under the normal exponential map of a $G$-invariant metric, and we use the corresponding identification of a neighbourhood of $S$ with $\mathcal N(r_0)$, the set of normal vectors of length less than $r_0$ for that metric, with $r_0>0$ small. If $S\subset\del M$, we use instead an equivariant collar $\mathcal N(r_0)=S\times[0,r_0)$, with $\mathcal N=S\times[0,\infty)$, in which the Euler-like vector field is $E=u\,\partial_u$, where $u$ is the collar coordinate.

In both cases, the flow of $E$ at time $\log s$ is the fibrewise multiplication $m_s$ by $s\in(0,1]$, and $m_0$ is the projection onto $S$. Thus $\partial_sm_s=s^{-1}\,E\circ m_s$. For a $k$-form $\beta$ on $\mathcal N(r)$, set
\[
P\beta=\int_0^1\frac1s\,m_s^*(\iota_E\beta)\,ds \in\Omega^{k-1}(\mathcal N(r)).
\]
$P$ is the standard homotopy operator of the retraction $m_s$, see~\cite[\S4.1]{Meinrenken-EulerLike}. If $\beta$ is closed and vanishes at every point of $S$, then $dP\beta=\beta$, and $P\beta$ vanishes along $S$ together with its first derivatives. The operator $P$ is linear, continuous and equivariant, it preserves smooth dependence on parameters, and $P\beta$ vanishes on $\mathcal N(\rho)$ whenever $\beta$ does.

\emph{Rescaling.} Let $\mathcal K$ be a compact subset of $\Symp^G_{\id,N}(M,S)$. By definition, its elements are tangent to the identity along $S$. Hence there is $r\in(0,r_0)$ such that the conjugate of $g$ by the flow of $E$,
\[
g^t=m_t^{-1}\circ g\circ m_t,
\]
is defined on $\mathcal N(r)$, with values in $\mathcal N(r_0)$, for all $g\in\mathcal K$ and $t\in(0,1]$. In a local trivialization of $\mathcal N$, write $g(x,v)=\big(x+a(x,v),\,v+b(x,v)\big)$, where $a$ and $b$ are smooth maps that vanish at $v=0$ together with their first derivatives. Then
\[
g^t(x,v)=\big(x+a(x,tv),\,v+t^{-1}b(x,tv)\big),
\]
so that $g^t$ extends smoothly to $t\in[0,1]$. Its value $g^0$ is the differential of $g$ along $S$, that is, $g^0=\id$. The maps $g^t$ are equivariant, fix $S$, are tangent to the identity along $S$, and depend continuously on $g$. Moreover $g^1=g$, and $g^t$ is the identity on $\mathcal N(\rho)$ whenever $g$ is.

\emph{Symplectic correction.} The closed forms $\tau_t=(g^t)^*\om-\om$ on $\mathcal N(r)$ vanish at every point of $S$. Hence $\beta_t=P\tau_t$ satisfies $d\beta_t=\tau_t$ and vanishes along $S$ with its first derivatives. After shrinking $r$, the forms $\om+s\tau_t$ are non-degenerate on $\mathcal N(r)$ for all $g\in\mathcal K$ and $s,t\in[0,1]$. Let $\mu_{t,s}$ be the flow of the vector field $W_{t,s}$ defined by $\iota_{W_{t,s}}(\om+s\tau_t)=-\beta_t$. Since $W_{t,s}$ vanishes along $S$ with its first derivatives, there is $r_1<r$ such that $\mu_{t,s}$ is defined on $\mathcal N(r_1)$ for $s\in[0,1]$, fixes $S$, and is tangent to the identity along $S$. Moser's argument gives $\mu_{t,1}^*(g^t)^*\om=\om$. Thus $\psi_t=g^t\circ\mu_{t,1}$ is an equivariant symplectic embedding of $\mathcal N(r_1)$ into $M$. It fixes $S$, is tangent to the identity along $S$, and depends smoothly on $t$ and continuously on $g$. Since $\tau_0=\tau_1=0$, we have $\psi_0=\id$ and $\psi_1=g$. If $g$ is the identity on $\mathcal N(\rho)$, then $\tau_t$ and $\beta_t$ vanish there, and $\psi_t$ is the identity there.

\emph{Extension.} The time-dependent vector field $Y_t=(\partial_t\psi_t)\circ \psi_t^{-1}$ is symplectic and vanishes on $S$. Hence $\psi_t^*(\iota_{Y_t}\om)=dk_t$ with $k_t=P\big(\psi_t^*\iota_{Y_t}\om\big)$, and $Y_t$ is the Hamiltonian vector field of $H_t=k_t\circ \psi_t^{-1}$, which vanishes along $S$ with its first derivatives. Given an invariant cut-off function $\chi$, let $\Phi_t$ be the Hamiltonian flow of $\chi H_t$, normalized by $\Phi_1=\id$. By construction, $\Phi_t\circ g=\psi_t$ on some smaller neighborhood $\mathcal N(\rho_2)$ and is the identity outside $\mathcal N(\rho_1)$. 

\emph{Conclusion.} Set $h_t=\Phi_t\circ g$. By construction, $h_t$ agrees with $\psi_t$ on $\mathcal N(\rho_2)$, so that $h_t\in\Symp^G_{\id,N}(M,S)$. Moreover $h_1=g$, $h_0$ is the identity on $\mathcal N(\rho_2)$, and if $g$ is the identity near $S$, then so is $h_t$ for every $t$. Since $h_t$ depends continuously on $g$, it defines a relative homotopy of any family $f:(D^k,\del D^k)\to\big(\Symp^G_{\id,N}(M,S),\Symp^G_{\id,\near}(M,S)\big)$ into $\Symp^G_{\id,\near}(M,S)$.
\end{proof}

We now study the group $\Symp^G(B^4)$ of equivariant symplectomorphisms of the
closed unit ball through its restriction to the boundary sphere $S^3=\del B^4$.
We first fix some notation for a collar of $S^3$ in $B^4$.

Let $\alpha$ denote the standard contact form on $S^3=\del B^4$, that is, the restriction to $S^3$ of $\iota_X\om_0$ where $X$ is the radial Liouville vector field of $(\R^4,\om_0)$, and let $R$ be its Reeb vector field. The flow $\theta_s$ of $R$ is the Hopf circle action $S^1(1,1)$. Let $T>0$ be its period. Set
\[
\mathcal{D}^G(S^3,d\alpha)=\{\psi\in\Diff^G(S^3)~|~\psi^*d\alpha=d\alpha\},
\]
the group of equivariant diffeomorphisms of $S^3$ preserving $d\alpha$, and let $\SCont^G(S^3,\alpha)\subset\mathcal{D}^G(S^3,d\alpha)$ be the subgroup of strict contactomorphisms, that is, of those $\psi$ with $\psi^*\alpha=\alpha$. Let 
\[
\Psi:S^3\times(-\delta,0]\to B^4,\qquad \Psi^*\om_0=d\big(e^\tau\alpha\big),\qquad \Psi_*\partial_\tau=X,
\]
be the symplectization collar obtained by flowing $S^3$ along $X$, where $\tau$ denotes the coordinate on $(-\delta,0]$. On this collar, the Hamiltonian vector field of a function $H$ is
\begin{equation}\label{eq:collar hamiltonian}
X_H=e^{-\tau}\,R(H)\,\partial_\tau\;-\;e^{-\tau}\,(\partial_\tau H)\;R\;+\;e^{-\tau}\big(d\alpha|_\xi\big)^{-1}d_\xi H ,
\end{equation}
where $\xi=\ker\alpha$ and $d_\xi$ denotes the restriction of the differential to $\xi$.

\begin{lemma}\label{lemma:ExtensionOfInfinitesimalAutomorphisms}
Let $Y$ be a $G$-invariant vector field on $S^3$ with $\mathcal{L}_Yd\alpha=0$. There is a $G$-invariant function $H\in C^\infty(B^4)$, depending linearly and continuously on $Y$, whose Hamiltonian vector field $X_{H}$ is tangent to $S^3$ and restricts to $Y$ along $S^3$.
\end{lemma}
\begin{proof}
The $1$-form $\iota_Yd\alpha$ is closed
so that $\iota_Yd\alpha=dF$ for a unique $G$-invariant function $F$ normalized by $\int_{S^3}F\,\alpha\wedge d\alpha=0$. Evaluating on $R$ gives $R(F)=d\alpha(Y,R)=0$.
Set $a=\alpha(Y)$ and decompose $Y=aR+Y_\xi$.  Then $\iota_{Y_\xi}d\alpha=dF$, that is, $Y_\xi=(d\alpha|_\xi)^{-1}d_\xi F$. Choose a smooth cut-off $\chi:\R\to [0,1]$ with $\chi\equiv1$ near $\tau=0$ and $\chi\equiv0$ near $\tau=-\delta$, and set
\[
H=\chi(\tau)\big(F-\tau\,a\big)\ \text{ on the collar},\qquad H=0\ \text{ elsewhere.}
\]
By construction, $H$ is $G$-invariant. On $S^3$ we have $H=F$, $\partial_\tau H=-a$, and $R(H)|_{S^3}=R(F)=0$. By~\eqref{eq:collar hamiltonian}, $X_{H}$ is tangent to $S^3$ and
\[
X_{H}|_{S^3}=a\,R+\big(d\alpha|_\xi\big)^{-1}d_\xi F=Y .
\]
Finally, the assignment $Y\mapsto(F,a)\mapsto H$ is clearly linear and continuous.
\end{proof}

\begin{lemma}\label{lemma:RestrictionIsAFibration}
The image of the restriction homomorphism $\rho:\Symp^G(B^4)\to \mathcal{D}^G(S^3,d\alpha)$ is a subgroup containing the identity component. Moreover, $\rho$ is a fibration over $\im(\rho)$ whose fibre is the group $\Symp^G_{\id}(B^4,S^3)$.
\end{lemma}
\begin{proof}
Let $\phi_t$ be a smooth isotopy in $\mathcal{D}^G(S^3,d\alpha)$ and let $\tilde{\phi}_0\in\Symp^G(B^4)$ be a lift of $\phi_0$. The path $\gamma_t=\phi_t\circ\phi_0^{-1}$ starts at the identity, and its generating field $Y_t=(\partial_t\gamma_t)\circ\gamma_t^{-1}$ is $G$-invariant and satisfies $\mathcal{L}_{Y_t}d\alpha=0$. Let $H_t$ be the extension given by Lemma~\ref{lemma:ExtensionOfInfinitesimalAutomorphisms}. The field $X_{H_t}$ is tangent to $\del B^4$, so its time-dependent flow $\Theta_t$ is defined on the closed ball $B^4$, with $\rho(\Theta_t)=\gamma_t$. As $\rho$ is a homomorphism, $\Theta_t\circ\tilde{\phi}_0$ is a lift of $\phi_t$ extending $\tilde{\phi}_0$.

Finally, since $\mathcal{D}^G(S^3,d\alpha)$ is locally contractible, this construction provides local sections over $\im(\rho)$, which implies that $\rho$ is a fibration over its image.
\end{proof}

\begin{lemma}\label{lemma:StrictContactRetraction}
The group $\mathcal{D}^G(S^3,d\alpha)$ deformation retracts onto its subgroup $\SCont^G(S^3,\alpha)$.
\end{lemma}

\begin{proof}
Let $\psi\in\mathcal{D}^G(S^3,d\alpha)$. The form $\psi^*\alpha-\alpha$ being closed, there exists a unique function $h_\psi$, normalized by $\int_{S^3}h_\psi\,\alpha\wedge d\alpha=0$, such that $\psi^*\alpha=\alpha+dh_\psi$. This potential  $h_\psi$ depends continuously on $\psi$. The Jacobian of $\psi$ is then
\[
\mathrm{Jac}(\psi) = 1+R(h_\psi).
\]
Since the Reeb flow preserves the volume, $\int_{S^3}R(h_\psi)\,\alpha\wedge d\alpha=0$, so that $1+R(h_\psi)$ is everywhere positive. In particular,  $1+sR(h_\psi)>0$ for every $s\in[0,1]$.

For $f\in C^\infty(S^3)$ consider the twist map  $\kappa_f(p)=\theta_{f(p)}(p)$. From 
\[
d(\kappa_f)_p(v)=df(v)\,R+d(\theta_{f(p)})_p(v)\quad \text{and}\quad \theta_s^*\alpha=\alpha ,
\]
we get
\begin{equation}\label{eq:kappa alpha}
\kappa_f^*\alpha=\alpha+df,
\end{equation}
and, decomposing $df=R(f)\,\alpha+\beta$ with $\beta(R)=0$ so that $\beta\wedge d\alpha=0$,
\begin{equation}\label{eq:kappa volume}
\kappa_f^*(\alpha\wedge d\alpha)=\big(1+R(f)\big)\,\alpha\wedge d\alpha .
\end{equation}
We claim that
$r_s(\psi)=\psi\circ\kappa_{sh_\psi}^{-1}$, $s\in[0,1]$,
is the desired retraction.

For $f=sh_\psi$, applying~\eqref{eq:kappa volume} to $\kappa_{f}$ shows that 
$\mathrm{Jac}(\kappa_{f}) = 1+sR(h_\psi)>0$. Thus the map $\kappa_f$ is a local diffeomorphism. It is also bijective. Indeed, $\kappa_f$ maps each Reeb orbit into itself. Now, for $p\in S^3$ and $c\in\R$, set $g_p(c)=f(\theta_c(p))+c$, so that $\kappa_f(\theta_c(p))=\theta_{g_p(c)}(p)$. Since $g_p'(c)=1+R(f)(\theta_c(p))>0$ and $g_p(c+T)=g_p(c)+T$, the map $g_p$ induces a bijection of $\R/T\Z$, and $\kappa_f$ maps the orbit of $p$ bijectively onto itself. Hence $\kappa_f$ is a diffeomorphism.

By~\eqref{eq:kappa alpha} we have $\kappa_f^*d\alpha=d\alpha$, so $r_s(\psi)$ lies in $\mathcal{D}(S^3,d\alpha)$. Clearly, $r_0(\psi)=\psi$. At time $s=1$ equation~\eqref{eq:kappa alpha} gives $(\kappa_{h_\psi}^{-1})^*(\alpha+dh_\psi)=\alpha$, so that $r_1(\psi)^*\alpha=\alpha$. Moreover, for a strict contactomorphism $\psi$, we have $dh_\psi=0$, hence $h_\psi=0$ and $r_s(\psi)=\psi$.

Finally, $\kappa_{sh_\psi}$ is equivariant. Indeed, every $g\in G\subset\U(2)$ preserves $\alpha$, $R$, the flow $\theta_s$, and the volume. From $\psi=g^{-1}\circ\psi\circ g$ we get $d(h_\psi\circ g)=dh_\psi$, so $h_\psi\circ g=h_\psi$ by the normalisation, and $g$ commutes with $\kappa_{sh_\psi}$.
\end{proof}

\begin{Remark}
The argument only uses that $\alpha$ is a contact form on a closed manifold $V$
with $H^1(V;\R)=0$ whose Reeb flow is periodic, and that $G$ acts by
diffeomorphisms of $V$ preserving $\alpha$. It applies, for instance, to the
lens spaces with the contact form induced by $\alpha$.
\end{Remark}

We now recall a fundamental result on equivariant symplectomorphisms of the standard $4$-ball whose proof is given in~\cite{ChPin-Memoirs}.

\begin{prop}{\rm (Equivariant Gromov contractibility theorem \cite[Theorem 3.29]{ChPin-Memoirs})}\label{prop:equivariant Gromov}
Let $V\subset\R^4$ be a compact star-shaped domain with $0$ in its interior, and let $\om$ be a symplectic form on $V$ that agrees with $\om_0$ near $\del V$. Let $G$ be a compact abelian group acting on $V$ symplectically, preserving $\del V$ and acting linearly near $\del V$. Then the group $\Symp^G_c(V,\om)$ of equivariant symplectomorphisms with compact supports in the interior of $V$ is non-empty and contractible. In particular, for a compact abelian subgroup $G$ of $\U(2)$ acting (effectively) on the standard symplectic ball $B^4(c)$, the group $\Symp^G_{\id,\near}(B^4(c), S^3)$ is contractible.\qed
\end{prop}

\begin{lemma}\label{lemma:KernelOfRestriction}
Suppose $G\subset\U(2)$ is abelian. The group $\Symp^G_{\id}(B^4,S^3)$ is homeomorphic to $\Symp^G_{\id,N}(B^4,S^3)\times C^\infty(S^3)^G$. In particular it is contractible.
\end{lemma}

\begin{proof}
Let $\phi\in\Symp^G_{\id}(B^4,S^3)$ and pick any point $p\in S^3$. By definition, $\phi(p)=p$ and $d\phi|_{T_pS^3}=\id$. We also have the identities
\[
\om_0(\phi_* X, \phi_* Y)=\om_0(X,Y)=\alpha(Y)\quad\text{and}\quad \om_0|_{TS^3}=d\alpha.
\]
Write $d\phi(X)=bX+v$ with $v\in T_pS^3$. Taking $Y\in\xi$, we get $\om_0(v,Y) = \om_0(bX+v,Y) = \alpha(Y) =0$, showing that $v$ is a multiple of $R$. Taking $Y=R$, we obtain $b=\om_0(bX+v,R)=\alpha(R)=1$. 
Hence
\[
d\phi(X)=X+\nu_\phi\,R\ \text{ along } S^3,\qquad \nu_\phi:=\alpha\big(d\phi(X)\big)\big|_{S^3},
\]
where $\nu_\phi$ is $G$-invariant because $G$ commutes with $\phi$ and preserves $X$, $R$, and $\alpha$. We easily check that $\nu_{\phi_1\circ\phi_2}=\nu_{\phi_1}+\nu_{\phi_2}$, so that the assignment
\begin{align*}
\mathcal{V}: \Symp^G_{\id}(B^4,S^3)&\to \big(C^\infty(S^3)^G,+\big) \\
\phi&\mapsto \nu_\phi
\end{align*}
is a continuous homomorphism with kernel $\Symp^G_{\id,N}(B^4,S^3)$.

It admits a continuous section. Choose a smooth cut-off $\chi:\R\to [0,1]$ with $\chi\equiv1$ near $\tau=0$ and $\chi\equiv0$ near $\tau=-\delta$. Given $\nu\in C^\infty(S^3)^G$, set $H_\nu=-\tfrac12\tau^2\,\nu\,\chi(\tau)$ on the collar and $H_\nu=0$ elsewhere. Then $dH_\nu$ vanishes along $S^3$, so that the flow of $X_{H_\nu}$ fixes $S^3$ pointwise and preserves $B^4$. Let $\Theta(\nu)$ be its time-one map. 
Near $\tau=0$, the formula~\eqref{eq:collar hamiltonian} gives
\begin{equation}\label{eq:X_H near S3}
X_{H_\nu}=\tau e^{-\tau}\nu\,R+O(\tau^2).
\end{equation}
As $X_{H_\nu}=0$ along $S^3$,
its linearization $A_p:T_pB^4\to T_pB^4$ vanishes on $T_pS^3$, and differentiating \eqref{eq:X_H near S3} in $\tau$ gives
$A_p(\partial_\tau)=\nu(p)\,R(p)$. Hence $A_p^2=0$.
As the differential of the time-one flow at a zero of the field is $\exp(A_p)=\id+A_p$, we get $\mathcal{V}(\Theta(\nu))=\nu$.

Therefore $(\phi,\nu)\mapsto\phi\circ\Theta(\nu)$ is a homeomorphism from $\Symp^G_{\id,N}(B^4,S^3)\times C^\infty(S^3)^G$ onto $\Symp^G_{\id}(B^4,S^3)$, with inverse $\phi\mapsto\big(\phi\circ\Theta(\mathcal{V}(\phi))^{-1},\mathcal{V}(\phi)\big)$.
The group $\Symp^G_{\id,N}(B^4,S^3)$ is homotopy equivalent to $\Symp^G_{\id,\near}(B^4,S^3)$ by Lemma~\ref{lemma:HomotopyEquivalenceSympN-SympNear}, hence contractible by Proposition~\ref{prop:equivariant Gromov}. The contractibility of  $\Symp^G_{\id}(B^4,S^3)$ follows.
\end{proof}

\begin{lemma}\label{lemma:SCont equivalent to U}
For any closed subgroup $G\subset\U(2)$ the inclusion
\[
\U^G(2)\into\SCont^G(S^3,\alpha)
\]
is a homotopy equivalence. In particular, $\SCont^G(S^3,\alpha)$ is connected.
\end{lemma}
\begin{proof}
The group $\SCont^G(S^3,\alpha)$ acts on $S^3$ preserving the Reeb flow which, in this case, is the circle action $S^1(1,1)$ that defines the Hopf fibration. Consequently, $\SCont^G(S^3,\alpha)$ acts by symplectomorphisms on the quotient $\CP^1\simeq S^2$. This action commutes with the induced action of $G$ on $\CP^1$.
We obtain a commutative diagram
\begin{equation}
\begin{tikzcd}
S^1(1,1) \arrow[hookrightarrow]{r} \arrow[equal]{d}&\SCont^G(S^3, \alpha) \arrow{r}{p} & \im(p)\subset\Symp^{\overline G}(S^2)\\
S^1(1,1) \arrow[hookrightarrow]{r}&\U^G(2) \arrow{r}{q} \arrow[hookrightarrow]{u}& \PU_0^{\overline G}(2)\arrow[hookrightarrow]{u}.
\end{tikzcd}
\end{equation}
Note that the kernel of $p$ consists of the strict contactomorphisms that preserve each fibre of the Hopf bundle $\beta:S^3\to S^2$. Since they also preserve the Reeb flow, these maps have the form $x\mapsto\lambda(\beta(x))\,x$, with $\lambda:S^2\to S^1$ smooth, and such a map preserves $\alpha$ only if $\lambda$ is constant. The kernel of $p$ is therefore $S^1(1,1)$.

As $\U^G(2)$ is connected by Lemma~\ref{Lemma:Sp^g(4)isconnected}, the map $q$ is a principal bundle over the identity component of $\PU^{\overline G}(2)$. 

We now show that $p$ is a fibration over the identity component $\Symp^{\overline G}_0(S^2)$. First, an equivariant symplectic isotopy of $S^2$ is generated by $\overline G$-invariant Hamiltonians $H_t$, and the contact Hamiltonians $H_t\circ\beta$ are invariant under $G$ and under the Reeb flow, so that their flows consist of elements of $\SCont^G(S^3,\alpha)$ lifting the given isotopy. As in Lemma~\ref{lemma:RestrictionIsAFibration}, $p$ is therefore a fibration over its image.
Second, let $\phi\in\SCont^G(S^3,\alpha)$ cover an element
$[u]\in\PU^{\overline G}(2)$. Since $u$ preserves $\alpha$, the map
$\phi\circ u^{-1}$ is a strict contactomorphism covering the identity. The
computation of the kernel of $p$ above shows that  $\phi\circ u^{-1}=\lambda$ for a constant $\lambda\in S^1$. Hence $\phi=\lambda u$ is a unitary map commuting with $G$, so that $u\in\U^G(2)$. Consequently, $[u]$ belongs to the identity component $\PU_0^{\overline G}(2)$.
Third, by Lemma~\ref{lemma:linear centralizers}, every element of $\Symp^{\overline G}(S^2)$ is joined by a path to an element $[u]$ of $\PU^{\overline G}(2)$. If it lies in the image of $p$, so does $[u]$.
Hence $p$ is a fibration over $\Symp^{\overline G}_0(S^2)$.

Finally, since $\PU^{\overline G}_0(2)\into\Symp^{\overline G}_0(S^2)$ is a homotopy equivalence by Lemma~\ref{lemma:linear centralizers}, the middle inclusion $\U^G(2)\into\SCont^G(S^3,\alpha)$ is a homotopy equivalence.
\end{proof}

We are now in a position to prove one of the main results of this section.

\begin{prop}\label{prop:homotopy type equivariant symplectomorphism of the ball}
Suppose $G\subset\U(2)$ is abelian. 
Then the group $\Symp^G(B^4(c))$ of $G$-equivariant symplectomorphisms of the standard ball $B^4(c)$ deformation retracts onto the linear centralizer subgroup $\U^G(2)$.
\end{prop}
\begin{proof}
Clearly it is enough to consider the unit ball $B^4$. Putting Lemma~\ref{lemma:StrictContactRetraction} and Lemma~\ref{lemma:SCont equivalent to U} together, we have homotopy equivalences
\[
\U^G(2)\into\SCont^G(S^3,\alpha)\into\mathcal{D}^G(S^3,d\alpha).
\]
In particular, $\mathcal{D}^G(S^3,d\alpha)$ is connected. Combining this with Lemma~\ref{lemma:RestrictionIsAFibration}, the restriction map 
\[
\rho:\Symp^G(B^4)\longrightarrow\mathcal{D}^G(S^3,d\alpha)
\]
is therefore a fibration whose fibre is the group $\Symp^G_{\id}(B^4, S^3)$. By Lemma~\ref{lemma:KernelOfRestriction}, $\Symp^G_{\id}(B^4,S^3)$ is contractible, hence $\rho$ is a homotopy equivalence as well. Combining the two steps, the inclusion $\U^G(2)\into\Symp^G(B^4)$ is a homotopy equivalence. As $\U^G(2)$ is a finite dimensional closed subgroup of $\Symp^G(B^4)$, Lemma~\ref{lemma:cofibration to deformation retraction} implies that $\Symp^G(B^4)$ deformation retracts onto $\U^G(2)$.
\end{proof}

%%%%%%%

Given an equivariant symplectic embedding $\iota:B^4(c)\into M$, let $\Stab^G_U(\iota)\subset \Stab^G(\iota)$ be the subgroup of equivariant symplectomorphisms that restrict to unitary maps on the image $\iota(B^4(c))\simeq B^4(c)$.

\begin{cor}\label{cor:stabilizer of a ball}
For $G\subset \U(2)$ abelian, there is a deformation retraction of $\Stab^G(\iota)$ onto its subgroup $\Stab^G_U(\iota)$.
\end{cor}
\begin{proof}
Consider the fibration
\[\Stab^G_U(\iota)\into\Stab^G(\iota)\to \Symp^G(B^4(c)) / \U^G(2)\]
obtained by projecting to the quotient the restriction to the image $\iota(B^4)$ of a symplectomorphism $\phi\in\Stab^G(\iota)$. By Proposition~\ref{prop:homotopy type equivariant symplectomorphism of the ball}, the quotient $\Symp^G(B^4(c)) / \U^G(2)$ is contractible, so that the fibre is homotopy equivalent to the total space. Since $\Stab^G_U(\iota)$ is a closed submanifold of $\Stab^G(\iota)$, Lemma~\ref{lemma:cofibration to deformation retraction} turns it into a deformation retraction.
\end{proof}

%%%%%

\begin{lemma}\label{lemma:equivalence Aut(N(E))}
Let $G\subset\U(2)$ be a closed subgroup, acting on the tautological bundle
$\Oo(-1)\to\CP^1$ through the linear action of $\U(2)$ on $\C^2$. Let
$\Aut^G_\om(\Oo(-1))$ be the group of symplectic bundle automorphisms of
$\Oo(-1)$ that commute with the action of $G$ and which cover symplectomorphisms
$\CP^1\to\CP^1$. There is a deformation retraction of $\Aut^G_\om(\Oo(-1))$ onto
its subgroup $\U^G(2)$.
\end{lemma}
\begin{proof}
The proof follows the same argument as the proof of Lemma~\ref{lemma:SCont equivalent to U}. 
Restriction to the zero section and projection to $\PU(2)$ give the commutative diagram
\begin{equation}
\begin{tikzcd}
\Gg_\om^G(\Oo(-1)) \arrow[hookrightarrow]{r} &\Aut^G_\om(\Oo(-1)) \arrow{r}{p} & \Symp^{\overline{G}}(S^2)\\
S^1(1,1) \arrow[hookrightarrow]{u}\arrow[hookrightarrow]{r}&\U^G(2) \arrow{r}{q} \arrow[hookrightarrow]{u}& \PU^{\overline{G}}_0(2)\arrow[hookrightarrow]{u}
\end{tikzcd}
\end{equation}
whose vertical maps are inclusions.

We start by showing the map $p$ is a fibrations over the identity component $\Symp^{\overline{G}}_0(S^2)$. Fix a $G$-invariant unitary connection on $\Oo(-1)$.
Let $h_t$ be a smooth path in $\Symp^{\overline G}(S^2)$ with $h_0=p(\phi_0)$. The vector fields generating $h_t\circ h_0^{-1}$ are ${\overline G}$-invariant, and so are their horizontal lifts, whose flow $\psi_t$ consists of equivariant unitary bundle automorphisms covering $h_t\circ h_0^{-1}$. Hence $\psi_t\circ\phi_0$ lifts $h_t$. Since $\Symp^{\overline G}(S^2)$ is locally contractible, this construction defines local sections of $p$ over its image.
Hence $p$
is a fibration over a union of path components of $\Symp^{\overline G}(S^2)$.

Choose $\phi\in\Aut^G_\om(\Oo(-1))$ such that $p(\phi)=[u]\in\PU^{\overline G}(2)$, where $u\in\U(2)$. We claim that $u\in\U^G(2)$. Each fibre $L_x$ of $\Oo(-1)$ is a complex line in $\C^2$, with the Hermitian product and the area form induced from $\C^2$. It follows that any symplectic linear map $A:L_x\to L_y$ admits a unique polar decomposition $A=U\,(A^*A)^{1/2}$, where the adjoint is taken with respect to these Hermitian products.
Applied fibrewise to $\phi$, it gives a deformation of $\phi$ to a unitary bundle automorphism $\phi^U$ that covers $[u]$,
and that commutes with $G$.
Then $u^{-1}\circ\phi^U$ acts on each fibre $L_x$ as an element of $\U(1)$, that is, 
$\phi^U=u\circ\tau$ with $\tau:S^2\to\U(1)$. Since $[u]$ commutes with $\overline G$, we have $gug^{-1}=c(g)\,u$ for a homomorphism $c:G\to\U(1)$, and the equivariance of $\phi^U$ implies $\tau(\bar g^{-1}x)\,c(g)=\tau(x)$ for all $x\in S^2$ and $g\in G$. But every $g\in G$ has an eigenvector, that is, a fixed point of $\bar g$ in $S^2$, at which this identity gives $c(g)=1$. Hence $u\in\U^G(2)$. 

Next, consider an arbitrary $\phi_1\in\Aut^G_\om(\Oo(-1))$. By Lemma~\ref{lemma:linear centralizers}, $p(\phi_1)$ is joined by a path in $\Symp^{\overline G}(S^2)$ to an element $[u]$ of $\PU^{\overline G}(2)$. Lifting this path, we obtain an element $\phi\in\Aut^G_\om(\Oo(-1))$ that covers $[u]$, and by the previous step, $p(\phi)\in q(\U^G(2))=\PU^{\overline G}_0(2)\subset\Symp^G_0(S^2)$. This shows that $p$ is a fibration over the identity component $\Symp^{\overline G}_0(S^2)$.

We now look at the fibres. The gauge group $\Gg_\om^G(\Oo(-1))$ equivariantly retracts onto the unitary gauge group $\Gg^G_U(\Oo(-1))$ by fibrewise polar decomposition. The latter is the equivariant gauge group of the unit circle bundle of $\Oo(-1)$. By Lemmas~\ref{lemma:gauge group as invariant maps} and~\ref{lemma:invariant circle valued maps}, it fits into homotopy equivalences
\[
\Gg^G_U(\Oo(-1)) \cong \Inv(S^2, S^1)\simeq S^1(1,1)\subset \U(2),
\]
where $\Inv(S^2, S^1)$ denotes the space of smooth maps from $S^2$ to $S^1$ that are invariant under the action of $G$ on $S^2$. Consequently, the leftmost vertical inclusion is also a homotopy equivalence.

Finally, since the inclusion $\PU^{\overline G}_0(2)\into\Symp^{\overline G}_0(S^2)$ is a homotopy equivalence by Lemma~\ref{lemma:linear centralizers} 
it follows that $\U^G(2)\into\Aut^G_\om(\Oo(-1))$ is a homotopy equivalence, and
Lemma~\ref{lemma:cofibration to deformation retraction} gives the deformation
retraction.
\end{proof}

%%%%%

\begin{lemma}\label{Lemma:Stab(E) vs StabU(E)}
Let $G\subset\U(2)$ be a compact subgroup. Let $\exdivisor$ be an invariant exceptional curve in $(M^4,\om)$. Let $\Symp^G(M,\exdivisor)$ be the group of equivariant symplectomorphisms of $M$ that send $\exdivisor$ to itself, and let $\Symp_U^G(M,\exdivisor)$ be the subgroup of those symplectomorphisms that act $\U(2)$-linearly near $\exdivisor$ with respect to a fixed symplectic embedding of a standard neighbourhood $L(\exdivisor)$. Then $\Symp^G(M,\exdivisor)$ retracts onto $\Symp^G_U(M,\exdivisor)$.
\end{lemma}
\begin{proof}
Restricting the derivative of $\phi\in \Symp^G(M,\exdivisor)$ to $T_\exdivisor M$ defines two fibrations
\begin{equation}
\begin{tikzcd}
\Symp^G_{\id,N}(M,\exdivisor) \arrow[hookrightarrow]{r} &\Symp^G(M,\exdivisor) \arrow[twoheadrightarrow]{r} & \Aut^G_\om(N(\exdivisor))\\
\Symp^G_{\id,\near}(M,\exdivisor) \arrow[hookrightarrow]{u}\arrow[hookrightarrow]{r}&\Symp^G_U(M,\exdivisor) \arrow[twoheadrightarrow]{r} \arrow[hookrightarrow]{u}& \U^G(2)\arrow[hookrightarrow]{u}.
\end{tikzcd}
\end{equation}
The leftmost inclusion is an equivalence by Lemma~\ref{lemma:HomotopyEquivalenceSympN-SympNear} while the rightmost inclusion is an equivalence by Lemma~\ref{lemma:equivalence Aut(N(E))}. The result follows from Lemma~\ref{lemma:cofibration to deformation retraction}.
\end{proof}

\begin{prop}\label{prop:homotopy equivalence Stab(B) and Stab(E)} Let $G$ be a compact abelian group acting unitarily on $B^4(c)$ and symplectically on $(M,\om)$. Let $\iota:B^4(c)\into M$ be an equivariant embedding. The equivariant stabilizer group $\Stab^G(\iota)$ of the ball $\iota(B^4(c))$ is homotopy equivalent to the equivariant stabilizer group $\Symp^G(\tilde{M}, \exdivisor,\tom)$ of the invariant exceptional curve in the equivariant symplectic blow-up.
\end{prop}
\begin{proof}
There is a chain of homotopy equivalences
\[
\Stab^G(\iota)\simeq \Stab^G_U(\iota) \cong \Symp^G_U(\tilde{M}, \exdivisor,\tom) \simeq \Symp^G(\tilde{M}, \exdivisor,\tom)
\]
The first homotopy equivalence follows from Corollary~\ref{cor:stabilizer of a ball}, and the third is Lemma~\ref{Lemma:Stab(E) vs StabU(E)}.

The second map is the homeomorphism induced by the blow-up construction described in more detail in the next section. Indeed, let $\phi\in\Stab^G_U(\iota)$. Outside $\iota(B^4(c))$ the blow-up $\pi$ is a diffeomorphism identifying $\tom$ with $\om$, by part~(1) of Proposition~\ref{prop:symplectic blowup on M}, so $\phi$ lifts there. On $\pi^{-1}(\iota(B^4(c)))$ the map $\phi$ is unitary, hence commutes with the blow-up of the local model, and its lift preserves $\tom$ by part~(2) of Proposition~\ref{prop:symplectic blowup on M} together with the $\U(2)$-invariance of the model form, which is part~(5) of Lemma~\ref{lemma:local model blow-up}. The two lifts agree on the overlap, so they define $\tilde\phi\in\Symp^G_U(\tilde{M},\exdivisor,\tom)$. Conversely, an equivariant symplectomorphism of $\tilde M$ preserving $\exdivisor$ and unitary near it descends to an element of $\Stab^G_U(\iota)$, and the two assignments are inverse to one another and continuous.
\end{proof} 

%%%%%%%%%%%%%%%%%%%%%%%%%%%%%%%%%%%%%%%%%%%%%%%%%%%%%%%%%%%%%%%%%%%%%%%%%%%%%%%%
\section{Deformations, isotopies, and connectedness of embedding spaces}
\label{sec:steps}
%%%%%%%%%%%%%%%%%%%%%%%%%%%%%%%%%%%%%%%%%%%%%%%%%%%%%%%%%%%%%%%%%%%%%%%%%%%%%%%%

In this section we show that, for a $4$-dimensional manifold $M$ endowed with a symplectic and homologically trivial action by a compact abelian group $G$ (possibly finite), the connectedness of the space $\EmbG$ of equivariant embeddings can be reduced to an equivariant version of McDuff's deformation to isotopy lemma~\cite{deftoiso}. Throughout this section we work in the following setting.

\begin{Setting}\label{set:1}
Let $(M,\omega)$ be a compact connected symplectic $4$-manifold. Let $G$ be a compact abelian group, possibly disconnected, acting as a subgroup of $\U(2) \subset \Sp(4,\R)$ on the ball $B^4(c)$.  Assume that $G$ acts symplectically on $(M,\om)$ in a homologically trivial way, and let $\fixedcomp \subset M$ be a connected component of the fixed point set $M^G$. 
\end{Setting}

\subsection{Reducing to embeddings with fixed germs at the origin}
In order to show connectedness of $\EmbG$, we first show that we can always equivariantly isotope an arbitrary embedding in $\EmbG$ so that it coincides with a given one near the origin.

\begin{lemma}\label{Lemma:Isotopywithinconnectedcomponent}
Let  $\iota \in \EmbG$ be an embedding and let $x \in \fixedcomp$ be a point  such that $\iota(0) \neq x$. Then there is an ambient equivariant Hamiltonian isotopy $\psi_t: M \rightarrow M$ with $t \in [0,1]$ such that \begin{enumerate}
        \item $\psi_0 = \id$,
        \item $\psi_1(\iota(0)) = x$,
        \item $\psi_t(\fixedcomp) = \fixedcomp$ for all $t$.
    \end{enumerate}
\end{lemma} 
\begin{proof}
The component $\fixedcomp$ is a closed symplectic submanifold of $M$, with $T_p\fixedcomp=(T_pM)^G$ for $p\in\fixedcomp$. Choose a Hamiltonian $H_t$ on $\fixedcomp$ whose time-one flow maps $\iota(0)$ to $x$. Extend $H_t$ smoothly to $M$ and average the extension over $G$. Since $G$ fixes $\fixedcomp$ pointwise, this gives invariant functions $\hat H_t$ on $M$ with $\hat H_t|_{\fixedcomp}=H_t$. The Hamiltonian vector field of $\hat H_t$ is $G$-invariant, hence tangent to $\fixedcomp$ at every point of $\fixedcomp$, and restricts to the Hamiltonian vector field of $H_t$. Its flow $\psi_t$ is therefore an equivariant Hamiltonian isotopy of $M$ that preserves $\fixedcomp$ and satisfies $\psi_1(\iota(0))=x$.
\end{proof}

\begin{prop}\label{Lemma:NormalisationofEmbeddings}
Let $\iota_0,\iota_1:(B^4(c),\wstd)\rightarrow (M,\om)$ be two equivariant symplectic embeddings in $\EmbG$. Then there exists an ambient equivariant Hamiltonian isotopy $\psi_t:(M,\om) \rightarrow (M,\om)$ such that $\psi_1 \circ \iota_1$ agrees with $\iota_0$ in a small neighbourhood of the origin. 
\end{prop}
\begin{proof}%
By Lemma~\ref{Lemma:Isotopywithinconnectedcomponent} we may assume that $\iota_0(0)=\iota_1(0)=p$. As $\iota_0(B^4(c))$ is a neighbourhood of $p$, the map $g=\iota_0^{-1}\circ\iota_1$ is defined on a small ball $B^4(\delta_0)$. It is an equivariant symplectic embedding into $B^4(c)$ with $g(0)=0$, and its derivative $L=dg(0)$ lies in $\Sp^G(4)$. We follow now the rescaling and extension steps of the proof of Lemma~\ref{lemma:HomotopyEquivalenceSympN-SympNear} with $S=\{0\}\subset B^4(c)$, so that $E$ is the radial vector field $E(x)=x$ and $m_t(x)=tx$.

\emph{Rescaling.} For $t\in[0,1]$ and $x\in B^4(\delta_0)$, set $g^t(x)=\int_0^1dg(stx)\,x\,ds$. This map is smooth in $(t,x)$, with $g^t=m_t^{-1}\circ g\circ m_t$ for $t>0$ and $g^0=L$. Each $g^t$ is equivariant, and is symplectic since $m_t^*\wstd=t^2\wstd$. 

\emph{Extension.} Let $j_t$, $t\in[0,1]$, be a smooth path of equivariant symplectic embeddings of a ball $B^4(\delta)$ into $B^4(c)$ such that $j_0=\id$ and $j_t(0)=0$. We claim that there are $\delta'\in(0,\delta)$ and an ambient equivariant Hamiltonian isotopy $\phi_t$ such that $\phi_t\circ\iota_0=\iota_0\circ j_t$ on $B^4(\delta')$. Set $\kappa_t=\iota_0\circ j_t$, and define $Y_t$ on $\kappa_t(B^4(\delta))$ by $\partial_t\kappa_t=Y_t\circ\kappa_t$. This vector field is $G$-invariant, because $\kappa_t$ is equivariant, and vanishes at $p$. As $\kappa_t^*\om=\wstd$ does not depend on $t$, the invariant $1$-form $\beta_t=\kappa_t^*\big(\om(Y_t,\cdot)\big)$ satisfies $d\beta_t=\kappa_t^*(\mathcal L_{Y_t}\om)=\partial_t(\kappa_t^*\om)=0$, and it vanishes at the origin. By the properties of $P$, the function $k_t=P\beta_t$ is then $G$-invariant, smooth in $t$, and satisfies $dk_t=\beta_t$. Let $\chi$ be a radial function supported in the interior of $B^4(\delta)$ and equal to $1$ on $B^4(\delta')$. The function $F_t$ equal to $(\chi k_t)\circ\kappa_t^{-1}$ on $\kappa_t(B^4(\delta))$ and to $0$ elsewhere is smooth and $G$-invariant. Its flow $\phi_t$ is the equivariant Hamiltonian isotopy of the claim.

\emph{Linear correction.} The derivative $L$ need not be the identity. However, since $\Sp^G(4)$ is connected by Lemma~\ref{Lemma:Sp^g(4)isconnected}, there is a smooth path $A_t$ in $\Sp^G(4)$ from the identity to $L$. After reparametrizations, concatenating $A_t$ with $g^t$, gives a smooth path $j_t$ of equivariant symplectic embeddings fixing the origin, with $j_0=\id$ and $j_1=g$. By continuity, these maps send $B^4(\delta)$ into $B^4(c)$ for $\delta$ small enough. The extension step above then gives $\phi_t$ with $\phi_1\circ\iota_0=\iota_0\circ g=\iota_1$ on $B^4(\delta')$, and $\phi_t^{-1}$ is the required isotopy.
\end{proof}

\begin{remark}
In Proposition~\ref{Lemma:NormalisationofEmbeddings} the neighbourhood of the origin on which two embeddings agree actually depends on the embeddings.
\end{remark}

\begin{defn}
Let $\iota:(B^4(c),\wstd)\rightarrow (M,\om)$ be an equivariant symplectic embedding centered on $\fixedcomp$. We define the space of normalized symplectic embeddings locally modelled on $\iota$ as
\[
\EmbGiota:=\left\{\phi\in \EmbG~|~\phi\text{~coincides with~}\iota\text{~near the origin}\right\}
\]
\end{defn}

Proposition~\ref{Lemma:NormalisationofEmbeddings} immediately implies the following statement.

\begin{cor}\label{cor:nsubspace}
The embedding space $\EmbG$ is path connected whenever any two elements of its
subspace $\EmbGiota$ can be joined by a path in $\EmbG$.\qed
\end{cor}

\subsection{Embeddings, symplectic blow-ups, and restrictions}
\label{sec:blowup}

In this subsection, we recall the symplectic blow-up construction in a formulation suitable for our study of embedding spaces. It is used in \S\ref{subsec:embtodef} to associate with two embeddings having the same germ at the origin a deformation of symplectic forms on the blow-up $\wM$. We follow \cite[Section 7.1]{McD-Salamon-Intro} with slightly different conventions.

\subsubsection{Conventions}
We measure balls by their capacities, and we normalize the Fubini--Study form $\om_{FS}$ so that $\CP^1$ has area one, as in Section~\ref{section:Centralizers on CP2}. In \cite[Section 7.1]{McD-Salamon-Intro}, balls are measured by their radii and $\CP^1$ has area $\pi$. The construction there depends on two parameters. The first one, $\lambda$, determines the area $\pi\lambda^2$ of the exceptional divisor, and the second one, $\varepsilon$, determines the ball of radius $\lambda+\varepsilon$ in which the construction takes place. We use instead the area $c$ of the exceptional divisor and the capacity $b>c$ of a ball containing the support of the construction. Given $0<c<b$, we set $\lambda=\sqrt{c/\pi}$, and we fix $\varepsilon=\varepsilon(c,b)>0$, depending smoothly on $(c,b)$, such that $\pi(\lambda+\varepsilon)^2<b$. 
Beside this translation, Lemma~\ref{lemma:local model blow-up} and Proposition~\ref{prop:symplectic blowup on M} below are statements of \cite[Section 7.1]{McD-Salamon-Intro}, and only the $G$-invariance requires an additional argument.

\subsubsection{Local model} Let $\pi:\wC\to\C^2$ be the complex blow-up of $\C^2$ at the origin, and let $\pr:\wC\to\CP^1$ be the tautological bundle. Let $\exdivisor=\pi^{-1}(0)$ be the exceptional divisor and, given $a>0$, let $L(a)=\pi^{-1}(B^4(a))$. For $c>0$, let $\wom_c=\pi^*\om_0+c\,\pr^*\om_{FS}$, so that $\exdivisor$ has area $c$. If $G$ acts unitarily on $\C^2$, the blow-up $\wC$ inherits a $G$-action for which the exceptional divisor is invariant.

\begin{lemma}[{\cite[Definition 7.1.12 and Proposition 7.1.13]{McD-Salamon-Intro}}]\label{lemma:local model blow-up}
Given $0<c<b$, there is a symplectic form $\wom^b_c$ on $\wC$ such that
\begin{enumerate}
\item $\wom^b_c$ coincides with $\wom_c$ on a neighbourhood of $\exdivisor$,
\item $\wom^b_c$ coincides with $\pi^*\om_0$ outside $L(b')$ for some $b'<b$,
\item $\wom^b_c$ is compatible with the standard complex structure $\wJ$ on $\wC$,
\item $\wom^b_c$ depends smoothly on $c$ and $b$,
\item $\wom^b_c$ is invariant under the standard $\U(2)$ action on $\wC$.
\end{enumerate}
\end{lemma}
\begin{proof}
Set $\wom^b_c=\widetilde\omega_{\lambda,\varepsilon}$, with $\lambda$ and $\varepsilon$ as above. This form is K\"ahler, and it satisfies (1), (2) and (4) by \cite[Proposition 7.1.13]{McD-Salamon-Intro}, with $b'=\pi(\lambda+\varepsilon)^2$. Part (5) holds because the form is obtained from $\om_0$ and $\om_{FS}$, which are $\U(2)$-invariant, through the map $z\mapsto f_{\lambda,\varepsilon}(|z|)\,z/|z|$, which is $\U(2)$-equivariant.
\end{proof}

\subsubsection{Symplectic blow-up along a single embedding} Let $g: B^4(b) \into (M,\omega)$ be an equivariant symplectic embedding centered at $p=g(0) \in \fixedcomp$. Let $J'$ be the push-forward on $g(B^4(b))$ of the standard complex structure, and let $J''$ be any (compatible) extension of $J'$ to $M$. By averaging the metric associated to the pair $(\om, J'')$, one obtains an invariant almost complex structure $J$ that extends the standard structure on $g(B^4(b))$. By construction, the embedding $g$ is then a Kähler parametrization of $(M,\om,J)$ near $p$. Let $\pi:(\widetilde{M},\widetilde{J})\to (M,J)$ be the  complex blow-up of $M$ at $p$. Since the action of $G$ on $g(B^4(b))$ is unitary, it lifts to an action on $\widetilde{M}$ that leaves the exceptional divisor $\exdivisor$ invariant and that is holomorphic near $\exdivisor$. The next proposition describes how we can construct symplectic forms on $\wM$ that interpolate between the symplectic form $\pi^*\om$ on $\wM\setminus \pi^{-1}(g(B^4(b)))$ and the model forms $\wom^b_c$ near $\exdivisor$.

\begin{prop}[{\cite[Theorem 7.1.21]{McD-Salamon-Intro}}]\label{prop:symplectic blowup on M}
For every $0<c<b$ there exists a unique symplectic form $\wom_{g,c}$ on $\widetilde{M}$ such that
\begin{enumerate}
\item the restriction of $\wom_{g,c}$ to $\wM\setminus \pi^{-1}(g(B^4(b)))$ coincides with $\pi^{*}\om$,
\item the restriction of $\wom_{g,c}$ to $\pi^{-1}(g(B^4(b)))$ coincides with the push-forward, under the blow-up of $g$, of the model form $\wom^b_c$ of Lemma~\ref{lemma:local model blow-up},
\item $\wom_{g,c}$ is compatible with $\wJ$,
\item $\wom_{g,c}$ represents the cohomology class
\[
\pi^*[\om]-c\,\PD([Z]),
\]
\item the form $\wom_{g,c}$ depends smoothly on $c$,
\item the form $\wom_{g,c}$ is invariant under the action of $G$ on~$\wM$.
\end{enumerate}
\end{prop}
\begin{proof}
Parts (1) -- (5) are \cite[Theorem 7.1.21]{McD-Salamon-Intro} for the normalized embedding $\psi=\psi_0=g$, with $J$ the push-forward of the standard complex structure by $g$. Part (6) follows from uniqueness: for $h\in G$, the form $\tilde h^*\wom_{g,c}$ satisfies (1) and (2), since $g$ is equivariant, $\om$ is invariant, and the model form is $\U(2)$-invariant by Lemma~\ref{lemma:local model blow-up}.
\end{proof}

\subsection{From embeddings to deformations}\label{subsec:embtodef}

Throughout this subsection we fix $\iota\in\EmbG$ and denote by $\wM$ the complex blow-up of $M$ at $p=\iota(0)$ with respect to the push-forward $\iota_*J_0$ of the standard complex structure.
By definition, every $g\in\EmbG$ extends to a symplectic embedding of $B^4(c+\eta)$ for some $\eta>0$, and the extension can be chosen equivariant.

Given two equivariant symplectic embeddings $g_0$ and $g_1$ in $\EmbG$, there is no canonical identification between the corresponding blow-ups, even as smooth manifolds.
When the two embeddings coincide with $\iota$ on a smaller ball $B^4(\delta')$, however, their blow-ups can be compared on the single manifold $\wM$ through a rescaling procedure taken from the proof of \cite[Corollary~1.5]{deftoiso}.
Instead of blowing up along $g$, we first deform $\om$ by a diffeomorphism $\phi_{g,a}$ supported in $g(B^4(c+\eta))$, chosen so that the small ball $\iota(B^4(\delta'))$ becomes a ball of capacity $a\delta'/\delta$ for the new form $\phi_{g,a}^*\om$. We then blow up along this small ball, with an exceptional divisor of area $a$. As the small ball and its complex structure do not depend on $g$, all these forms live on the same manifold $\wM$. For $a=\delta$ there is no deformation, and the form does not depend on $g$. For $a=c$, the map $\phi_{g,c}$ sends the small ball onto $g(B^4(c\delta'/\delta))$, and the form is the pull-back, by the lift of $\phi_{g,c}$, of the blow-up of $(M,\om)$ along $g$. Letting $a$ vary from $c$ to $\delta$ thus joins the blow-up along $g$ to a form that is the same for all $g$.

\begin{lemma}{\rm (Rescaled blow-up)}\label{lemma:rescaled blow-up}
Let $g\in\EmbGiota$ be extended to $B^4(c+\eta)$, and let $0<\delta<\delta'$ be such that $g$ coincides with $\iota$ on $B^4(\delta')$ and $c\,\delta'<\delta\,(c+\eta)$. For $\delta\leq a\leq c$ there are equivariant diffeomorphisms $\phi_{g,a}$ of $M$ and invariant symplectic forms $\wt\sigma_{g,a}$ on $\wM$, depending smoothly on $a$, such that
\begin{enumerate}
\item $\phi_{g,a}^*\om=\frac{a}{\delta}\,\om$ on $\iota(B^4(\delta'))$, and $\phi_{g,\delta}=\id$,
\item $\wt\sigma_{g,a}$ is the form of Proposition~\ref{prop:symplectic blowup on M} for the symplectic manifold $(M,\phi_{g,a}^*\om)$, the embedding $z\mapsto\iota(\sqrt{\delta/a}\,z)$ of $B^4(a\delta'/\delta)$, and the area $a$,
\item $\wt\sigma_{g,a}$ represents the class $\pi^*[\om]-a\,\PD([Z])$,
\item on a neighbourhood of $\exdivisor$ that depends neither on $g$ nor on $a$, the form $\wt\sigma_{g,a}$ is the push-forward of $\frac{a}{\delta}\,\wom_\delta$ under the blow-up of $\iota$,
\item $\wt\sigma_{g,\delta}$ does not depend on $g$.
\end{enumerate}
\end{lemma}
\begin{proof}
Choose an increasing diffeomorphism $k$ of $[0,\infty)$ equal to $u\mapsto\frac{c}{\delta}u$ on $[0,\delta']$ and to the identity on $[\bar c,\infty)$, where $c\,\delta'/\delta<\bar c<c+\eta$. For $\delta\leq a\leq c$, the convex combination $k_a=(1-\tau)\,\id+\tau\,k$ with $\tau=(a-\delta)/(c-\delta)$ is again an increasing diffeomorphism. It equals $u\mapsto\frac{a}{\delta}u$ on $[0,\delta']$, and $k_\delta=\id$. The radial map
\[
\mu_a(z)=\bigl(k_a(\pi|z|^2)/\pi|z|^2\bigr)^{1/2}\,z
\]
is therefore a $\U(2)$-equivariant diffeomorphism of $\C^2$, equal to $\sqrt{a/\delta}\,z$ on $B^4(\delta')$ and to the identity outside $B^4(\bar c)$. Let $\phi_{g,a}$ be $g\circ \mu_a\circ g^{-1}$ on $g(B^4(c+\eta))$ and the identity elsewhere. It is an equivariant diffeomorphism depending smoothly on $a$, and (1) holds because $g$ is symplectic and $\mu_a$ multiplies $\wstd$ by $a/\delta$ on $B^4(\delta')$.

By (1), the formula $h_a(z)=\iota(\sqrt{\delta/a}\,z)$ defines an equivariant symplectic embedding of $B^4(a\delta'/\delta)$ into $(M,\phi_{g,a}^*\om)$ with image $\iota(B^4(\delta'))$. Since dilations are holomorphic, the push-forward of the standard complex structure under $h_a$ is $\iota_*J_0$, so the blow-up of Proposition~\ref{prop:symplectic blowup on M} takes place on $\wM$ for all $g$ and $a$. As $a<a\delta'/\delta$, that proposition applies with $b=a\delta'/\delta$ and the area $a$, and defines $\wt\sigma_{g,a}$ as in (2). It is invariant by part~(6) there, and depends smoothly on $a$ by parts~(1) and~(2), which describe it explicitly. Point (3) follows from part~(4), since $\phi_{g,a}$ is isotopic to the identity. Near $\exdivisor$, the form is the push-forward of $\wom_a$ under the blow-up of $h_a$, by part~(1) of Lemma~\ref{lemma:local model blow-up}. The dilation $z\mapsto\sqrt{a/\delta}\,z$ lifts to $\wC$, commutes with $\pr$, and multiplies $\om_0$ by $a/\delta$, so it pulls $\wom_a$ back to $\frac{a}{\delta}\pi^*\om_0+a\,\pr^*\om_{FS}=\frac{a}{\delta}\,\wom_\delta$. This gives point (4), the neighbourhood being chosen uniformly for $a$ in the compact interval $[\delta,c]$, and independently of $g$ because on $\pi^{-1}(\iota(B^4(\delta')))$ the construction involves only $\iota$. Finally, for $a=\delta$ we have $\phi_{g,\delta}=\id$ and $h_\delta=\iota$, so $\wt\sigma_{g,\delta}$ is the blow-up of $(M,\om)$ along $\iota$, which proves (5).
\end{proof}

\begin{prop}\label{prop:deformation from embeddings}
Let $g_0,g_1\in\EmbGiota$. There are $0<\delta<c$, a neighbourhood $V$ of $\exdivisor$ in $\wM$, and a smooth path $\wt\sigma_t$, $0\leq t\leq 1$, of invariant symplectic forms on $\wM$ such that
\begin{enumerate}
\item $\wt\sigma_0=\wt\sigma_{g_0,c}$ and $\wt\sigma_1=\wt\sigma_{g_1,c}$,
\item $\wt\sigma_t$ represents $\pi^*[\om]-c_t\,\PD([Z])$, for a smooth function $\delta\leq c_t\leq c$ with $c_0=c_1=c$,
\item on $V$, the form $\wt\sigma_t$ is the push-forward of $\frac{c_t}{\delta}\,\wom_\delta$ under the blow-up of $\iota$,
\item for every $\delta\leq a\leq c$, the truncation $\wt\sigma^a_t$ of $\wt\sigma_t$, that is, its restriction to the interval $\{t\mid c_t\leq a\}$, is, up to reparametrization, the deformation given by the proposition for the restricted pair $(g_0|_{B^4(a)},g_1|_{B^4(a)})$ when the choices made for $(g_0,g_1)$ are restricted to it, with the same $\delta$ and $V$. These deformations are thus determined by $\wt\sigma_t$, and $\wt\sigma^\delta_t$ is constant.
\end{enumerate}
\end{prop}
\begin{proof}
Extend $g_0$ and $g_1$ to $B^4(c+\eta)$ for a common $\eta>0$, and choose $0<\delta<\delta'$ such that both coincide with $\iota$ on $B^4(\delta')$ and $c\,\delta'<\delta\,(c+\eta)$. This is possible since both coincide with $\iota$ near the origin. Let $c_t$ be a smooth function on $[0,1]$ that decreases from $c$ to $\delta$ on $[0,\frac12]$, increases back to $c$ on $[\frac12,1]$, and has all its derivatives equal to zero at $t=\frac12$. Set $\wt\sigma_t=\wt\sigma_{g_0,c_t}$ for $t\leq\frac12$ and $\wt\sigma_t=\wt\sigma_{g_1,c_t}$ for $t\geq\frac12$. The two definitions agree at $t=\frac12$ by part~(5) of Lemma~\ref{lemma:rescaled blow-up}, and the path is smooth there because $c_t$ is flat. Properties (2) and (3) follow from parts (3) and (4) of the lemma. For (4), the choices made for $g_i$ restrict to $g_i|_{B^4(a)}$: take $g_i$ itself as the extension, the same $\delta$ and $\delta'$, and $k_a$ in place of $k$ in the proof of Lemma~\ref{lemma:rescaled blow-up}. The hypotheses of that lemma hold, since $a\,\delta'\leq c\,\delta'<\delta\,(c+\eta)$. The convex combinations of the identity and $k_a$ are the $k_{a'}$, $\delta\leq a'\leq a$, so the forms obtained are the $\wt\sigma_{g_i,a'}$, and the path given by the proposition for the restricted pair is, up to reparametrization, the truncation $\wt\sigma^a_t$.
\end{proof}

\subsection{From isotopies to connectedness of embedding spaces}\label{sub:fromisotopy}

\begin{prop}\label{prop:isotopy implies connectedness}
Let $g_0,g_1\in\EmbGiota$, let $\wt{\sigma}_t$ be a deformation given by Proposition~\ref{prop:deformation from embeddings}, and let $\wt{\sigma}^a_t$, $\delta\leq a\leq c$, be its truncations, as in part~(4) of that proposition. Assume that the truncations are homotopic relative to their endpoints, continuously in $a$, to isotopies $\wt{\tau}^a_t$ of invariant forms for which the exceptional divisor $\exdivisor$ remains symplectic at all times, the homotopy being constant for $a=\delta$. Then $g_0$ and $g_1$ lie in the same path component of $\EmbG$. Consequently, if this holds for every such pair, then $\EmbG$ is path connected, since by Proposition~\ref{Lemma:NormalisationofEmbeddings} every element of $\EmbG$ is joined to an element of $\EmbGiota$.
\end{prop}
\begin{proof}
We follow the proof of \cite[Corollary~1.5]{deftoiso}, and write $\wt{\tau}_t:=\wt{\tau}^c_t$, an isotopy from $\wt{\sigma}_0$ to $\wt{\sigma}_1$. For each $t$, the equivariant symplectic neighbourhood theorem, applied to $\wt{\tau}_t$ and $\wt{\tau}_0$ near $\exdivisor$, is a Moser argument at fixed $t$ that depends only on $\wt{\tau}_t$ near $\exdivisor$ and gives the identity when $\wt{\tau}_t=\wt{\tau}_0$ there. It yields diffeomorphisms $\wt{h}_t$, depending continuously on $t$, supported in a neighbourhood $N\subset V$ of $\exdivisor$, with $\wt{h}_t^*\wt{\tau}_t=\wt{\tau}_0$ on a smaller neighbourhood $N'\subset N$. Since $\wt{\tau}_0=\wt{\sigma}_0$ and $\wt{\tau}_1=\wt{\sigma}_1$ coincide on $V$, $\wt{h}_0=\wt{h}_1=\id$.
Applying an equivariant Moser argument to the cohomologous family $\wt{h}_t^*\wt{\tau}_t$, we find an isotopy $\wt{\psi}_t$ on $\wM$ such that $\wt{\psi}_t^*\wt{h}_t^*\wt{\tau}_t = \wt{\sigma}_0$, with $\wt{\psi}_0=\id$ and $\wt{\psi}_t=\id$ on $N'$, shrinking $N'$ if necessary: the derivative of the family vanishes on $N'$, so, as $H^1(N')=0$, its primitives can be chosen to vanish there, and invariant after averaging over $G$. 

The diffeomorphism $\wt{h}_1\circ \wt{\psi}_1$ being the identity on $N'$, it descends, through the description of the blow-up by surgery \cite[Theorem 7.1.21(ii) and (7.1.14)]{McD-Salamon-Intro}, in which the neighbourhood of $\exdivisor$ can be taken as thin as needed, to a diffeomorphism $h\circ\psi$ of $M$ which is the identity on a neighbourhood of $\iota(B^4(\delta))$ and which, by construction, pulls back $\phi_{g_1,c}^*\om$ to $\phi_{g_0,c}^*\om$. It follows that $H=\phi_{g_1,c}\circ h\circ\psi \circ\phi_{g_0,c}^{-1}$ is a symplectomorphism that takes the embedding $g_0$ to $g_1$.

Finally, the construction above depends continuously on the isotopy: the diffeomorphisms $\wt{h}_t$ are obtained at fixed $t$ by a canonical construction, and the Moser primitives can be chosen continuously. Applied to the isotopies $\wt{\tau}^a$ and to the restrictions $g_i|_{B^4(a)}$, whose deformations are the truncations of $\wt{\sigma}_t$ by part~(4) of Proposition~\ref{prop:deformation from embeddings}, it yields equivariant symplectomorphisms $H^a$ depending continuously on $a$, with $H^c=H$, and $H^\delta=\id$ since $\wt{\tau}^\delta=\wt{\sigma}^\delta$ is constant and $\phi_{g_i,\delta}=\id$. The path $H^a\circ g_0$, $\delta\leq a\leq c$, joins $g_0$ to $g_1$ in $\EmbG$.
\end{proof}

%%%%%%%%%%%%%%%%%%%%%%%%%%%%%%%%%%%%%%%%%%%%%%%%%%%%%%%%%%%%%%%%%%%%%%%%%%%%%%%%
\subsubsection{How the two propositions are used}
\label{section:role of 3.12 and 3.13}
%%%%%%%%%%%%%%%%%%%%%%%%%%%%%%%%%%%%%%%%%%%%%%%%%%%%%%%%%%%%%%%%%%%%%%%%%%%%%%%%

Together, Propositions~\ref{prop:deformation from embeddings} and~\ref{prop:isotopy implies connectedness} turn a deformation-to-isotopy statement into connectedness of a space of equivariant embeddings. 
What is left to be constructed is therefore a $G$-equivariant deformation-to-isotopy lemma of the following form: Given a deformation $\wt\sigma_t$ as in Proposition~\ref{prop:deformation from embeddings}, it must produce isotopies $\wt\tau^a_t$ of invariant forms such that
\begin{enumerate}
\item $\wt\tau^a$ is homotopic to the truncation $\wt\sigma^a$ relative to its endpoints;
\item the exceptional divisor stays symplectic along $\wt\tau^a$;
\item the homotopy depends continuously on $a$ and is constant for $a=\delta$.
\end{enumerate}
Section~\ref{sec:toric} proves one such lemma for $4$-dimensional toric manifolds, while Section~\ref{sec:deformation to isotopy} proves one for $S^1$ actions on ruled surfaces and on $\CP^2$.

\begin{Remark}
Since the path joining $g_0$ to $g_1$ is not obtained by blowing down the isotopy $\wt{\tau}_t$, the blow-up construction does not, a priori, induce an identification of path components between embedding spaces and spaces of blow-up forms.
\end{Remark}

%%%%%%%%%%%%%%%%%%%%%%%%%%%%%%%%%%%%%%%%%%%%%%%%%%%%%%%%%%%%%%%%%%%%%%%%%%%%%%%%
\section{Hamiltonian toric actions on $4$-manifolds}
\label{sec:toric}
%%%%%%%%%%%%%%%%%%%%%%%%%%%%%%%%%%%%%%%%%%%%%%%%%%%%%%%%%%%%%%%%%%%%%%%

We illustrate the framework developed in the previous sections by giving an alternative proof of Pelayo's result on ball embeddings equivariant under a toric action. Although our proof only applies to dimension $4$ (the original proof relies solely on moment map techniques and holds in any dimension), it serves as a preparation for the case of circle actions treated in the following sections.

Let $(M,\om)$ be a symplectic $4$-manifold equipped with a toric action $\rho:\T\times M\to M$ with moment map $\phi:M\to\mathfrak{t}^*$ and moment polygon $\Delta=\phi(M)$. The moment map can be normalized by sending the center of mass of $\Delta$ to the origin in $\mathfrak{t}^*$, which makes $\phi$ completely canonical. Let $p\in M$ be a fixed point and endow the ball $B^4(c)$ with a toric action with the same weights as the ones at $p$. Let $c_p$ be the symplectic area of the smallest invariant sphere containing $p$ (i.e. the rational length of the shortest edge meeting at $\phi(p)$). There is an equivariant symplectic embedding of $B^4(c)$ centered at $p$ if and only if $0<c<c_p$. In other words, $c_p$ is the equivariant capacity $c^{\T}(M,p)$ of Definition~\ref{defn:capacity} and Lemma~\ref{lemma:capacity basics}.

\begin{prop}[{\cite[Theorem 1.1]{pelayo}}]\label{prop:Pelayo}
For any capacity $0<c<c_p$, the space of equivariant embeddings $\Emb^{\T}(M,p,c)$ is homotopy equivalent to $\T$.
\end{prop}

We first recall two results that we will use repeatedly in our analysis of equivariant embeddings. The first one gives a sufficient condition for an invariant surface to be $J$-holomorphic with respect to invariant almost complex structures.

\begin{lemma}[{\cite[Lemma~2.5]{LiatRiver_noExtension}}]\label{LR}
Let $(M,J)$ be an almost complex manifold on which a compact Lie group $G$ acts leaving $J$ invariant. Let $S$ be an embedded sphere in $M$. Assume that $S$ is a connected component of the fixed point set of a non-trivial subgroup $H\subseteq G$. Then $S$ is $J$-holomorphic.\qed
\end{lemma}

The second result is the equivariant version of the $J$-tame inflation lemma. It is our main tool for turning symplectic deformations into isotopies.

\begin{lemma}[Equivariant inflation, {\cite[Corollary~1.4]{CPP}}]\label{lem:inf}
Let $(M,\om)$ be a closed symplectic $4$-manifold with a symplectic action of a compact Lie group $G$, and let $S\subset M$ be a $G$-invariant embedded symplectic surface which is $J$-holomorphic for some invariant $\om$-tame almost complex structure $J$. Set $s_S=\infty$ if $S\cdot S\geq 0$ and $s_S=\om(S)/(-S\cdot S)$ if $S\cdot S<0$. Then for every $0\leq s<s_S$ there is an invariant symplectic form taming $J$ in the class $[\om]+s\,\PD[S]$.\qed
\end{lemma}

The construction of~\cite{CPP} adapts to smooth families of data. When $J$ is compatible with $\om$ along $S$, it is the Thom-form construction of~\cite[\S3 and Remark~3.1]{CPP} (See also \cite{Buse-inflation} and \cite{McDuff-ALS}).

The proof of Proposition~\ref{prop:Pelayo} occupies the rest of this section. Let $\wM$ be the equivariant blow-up of $M$ at $p$, with the lifted toric action $\ttT$, and let $S_1,\ldots,S_k$ be the invariant spheres over the edges of its moment polygon (the exceptional divisor $\exdivisor$ is one of them), with primitive inward normals $n_1,\ldots,n_k$. Each $S_i$ is a component of the fixed point set of a circle subgroup of $\ttT$, so by Lemma~\ref{LR} it is $J$-holomorphic for every invariant almost complex structure $J$, hence symplectic for every invariant symplectic form.

\begin{lemma}\label{lem:toric correction}
To every $\ttT$-invariant symplectic form $\sigma$ on $\wM$ and every class $\kappa$ of such a form one can associate a deformation $C_t(\sigma,\kappa)$, $0\leq t\leq 1$, of invariant symplectic forms from $\sigma$ to a form in the class $\kappa$, such that
\begin{enumerate}
\item every $S_i$ is symplectic for every $C_t(\sigma,\kappa)$,
\item $C_t(\sigma,\kappa)$ depends continuously on $(\sigma,\kappa)$,
\item $C_t(\sigma,[\sigma])=\sigma$ for all $t$.
\end{enumerate}
\end{lemma}
\begin{proof}
Since $H^1(\wM)=0$, every invariant symplectic form $\sigma$ has a moment map, which we normalize by sending the center of mass of its image to the origin. The resulting polygon has the normals $n_i$, its edges being the images of the $S_i$, and positive support constants $h_i(\sigma)$, on which the class $[\sigma]$ depends linearly. Any other class $\kappa$ of invariant forms determines constants $h_i(\kappa)$ in the same way. Set
\[
\mu=\min\Bigl(1,\ \min_i\frac{h_i(\kappa)}{h_i(\sigma)}\Bigr),\qquad s_i=h_i(\kappa)-\mu\,h_i(\sigma)\geq 0 .
\]
The path rescales $\sigma$ to $\mu\sigma$, then inflates along the $S_i$ so as to increase the support constants by the amounts $s_i$, using Lemma~\ref{lem:inf} with, before each step, the invariant almost complex structure compatible with the current form obtained by polar decomposition with respect to a fixed invariant metric. The support constants move from $\mu h(\sigma)$ to $h(\kappa)$ along a segment of the open convex cone of support constants of polygons with normals $n_i$. Performing the inflations in small alternating steps keeps them near this segment, hence within the bounds of Lemma~\ref{lem:inf}, the bound for $S_i$ expressing that the $i$-th edge keeps a positive length. The final class is $\kappa$, and (1) holds by the remark preceding the lemma. All choices depend continuously on the data, which gives (2). If $[\sigma]=\kappa$, then $\mu=1$ and every $s_i=0$, and an inflation by zero along a surface on which $J$ is compatible with the form returns that form, which gives (3).
\end{proof}

\begin{cor}\label{cor:toric connected}
Let $g_0$ and $g_1$ be two embeddings in $\Emb^{\T}_\iota(M,p,c)$, and let $\wt\sigma_t$ be a deformation given by Proposition~\ref{prop:deformation from embeddings}. Its truncations satisfy the hypothesis of Proposition~\ref{prop:isotopy implies connectedness}. Consequently, $\Emb^{\T}(M,p,c)$ is path connected for every $0<c<c_p$.
\end{cor}
\begin{proof}
The endpoints of each truncation $\wt\sigma^a_t$ lie in the class $\kappa_a=\pi^*[\om]-a\,\PD([\exdivisor])$, by part~(3) of Lemma~\ref{lemma:rescaled blow-up}. By Lemma~\ref{lem:toric correction}, the map $(t,u)\mapsto C_u(\wt\sigma^a_t,\kappa_a)$ is a homotopy from $\wt\sigma^a_t$ to a path in the class $\kappa_a$ along which $\exdivisor$ stays symplectic. It is relative to the endpoints by part~(3), continuous in $a$ by part~(2), and constant for $a=\delta$, since $\wt\sigma^\delta_t$ is constant. After the truncations are reparametrized to be constant near their endpoints, smoothing these paths in $t$ yields isotopies, so Proposition~\ref{prop:isotopy implies connectedness} applies.
\end{proof}

\begin{cor}\label{cor:toric fibrations}
There are evaluation fibrations
\begin{equation}\label{eq:ToricFibration}
\begin{tikzcd}
\Fix^{\T}(\iota) \arrow[hookrightarrow]{r} \arrow[hookrightarrow]{d}&\Symp^{\T}(M,\om) \ar[r,twoheadrightarrow] \arrow[equal]{d}& \Emb^{\T}(M,p,c)\arrow[twoheadrightarrow]{d}\\
\Stab^{\T}(\iota) \arrow[hookrightarrow]{r}&\Symp^{\T}(M,\om) \arrow[twoheadrightarrow]{r}& \IEmb^{\T}(M,p,c)
\end{tikzcd}
\end{equation}
where $\Symp^{\T}(M,\om)$ is the group of equivariant symplectomorphisms, that is, the centralizer of the toric action in $\Symp(M,\om)$.
\end{cor}
\begin{proof}
By~\cite[Proposition 3.21]{Pin-MaxTori}, the group $\Symp^{\T}(M,\om)$ retracts onto the torus $\T$. In particular, it is connected and therefore sends each isolated fixed point to itself, so that it coincides with $\Symp^{\T}(M,\{p\},\om)$. The fibrations are those of diagram~\eqref{eq:MainFibration}, which applies since $\Emb^{\T}(M,p,c)$ is connected by Corollary~\ref{cor:toric connected}.
\end{proof}

\begin{lemma}\label{lemma:toric IEmb contractible}
The space of unparametrized equivariant embeddings $\IEmb^{\T}(M,p,c)$ is contractible.
\end{lemma}
\begin{proof}
By Proposition~\ref{prop:homotopy equivalence Stab(B) and Stab(E)}, the stabilizer $\Stab^{\T}(\iota)$ of the ball $\iota(B^4(c))$ is homotopy equivalent to the group $\Symp^{\ttT}(\wM,\exdivisor,\tom)$ of equivariant symplectomorphisms of the equivariant blow-up at $p$ that preserve the exceptional curve $\exdivisor$. By~\cite[Proposition 3.21]{Pin-MaxTori}, the group $\Symp^{\ttT}(\wM,\tom)$ retracts onto the torus $\ttT$. In particular it is connected, hence sends each invariant sphere to itself, so that
\[
\Symp^{\ttT}(\wM,\tom) = \Symp^{\ttT}(\wM,\exdivisor,\tom).
\]
The equivariant blow-down identifies $\ttT$ with $\T$ itself and yields natural homotopy equivalences
\[
\Stab^{\T}(\iota)\simeq \Symp^{\ttT}(\wM,\tom)\simeq \ttT\simeq \T.
\]
Therefore, there is a ladder of fibrations
\begin{equation}\label{eq:Final diagram for toric actions}
\begin{tikzcd}
\Stab^{\T}(\iota) \ar[r, hookrightarrow]&\Symp^{\T}(M,\om) \arrow{r} & \IEmb^{\T}(M,p,c)\\
\T \ar[u,hookrightarrow,swap,"\simeq"]\ar[r] &\T \ar[u,hookrightarrow,swap,"\simeq"]\ar[r]& *\ar[u,hookrightarrow]
\end{tikzcd}
\end{equation}
in which the first two vertical arrows are equivalences. In the top row, the inclusion of the fibre $\Stab^{\T}(\iota)$ into $\Symp^{\T}(M,\om)$ is therefore a homotopy equivalence. As $\IEmb^{\T}(M,p,c)$ is connected by Corollary~\ref{cor:toric connected}, the long exact sequence of homotopy groups shows that all its homotopy groups vanish. Having the homotopy type of a CW complex, it is contractible.
\end{proof}

\begin{proof}[Proof of Proposition~\ref{prop:Pelayo}]
Since $\IEmb^{\T}(M,p,c) := \Emb^{\T}(M,p,c) / \Symp^{\T}(B^4(c))$ is contractible, the space $\Emb^{\T}(M,p,c)$ is homotopy equivalent to the equivariant reparametrization group $\Symp^{\T}(B^4(c))$. By Proposition~\ref{prop:homotopy type equivariant symplectomorphism of the ball}, the latter is homotopy equivalent to the torus $\T$.
\end{proof}

%%%%%%%%%%%%%%%%%%%%%%%%%%%%%%%%%%%%%%%%%%%%%%%%%%%%%%%%%%%%%%%%%%%%%%%%%%%%%%%%
\section{Hamiltonian circle actions on $\SSS$, $\CCC$ and their blow-up}
\label{section:circle actions}
%%%%%%%%%%%%%%%%%%%%%%%%%%%%%%%%%%%%%%%%%%%%%%%%%%%%%%%%%%%%%%%%%%%%%%%%%%%%%%%%

This section collects what the study of equivariant ball embeddings in the
rational ruled surfaces requires. In~\S\ref{subsec:ruled} we fix the
normalisations and recall, from~\cite[\S2.2]{ChPin-Memoirs}, the list of toric
and circle actions on $\SSS$ and $\CCC$. In~\S\ref{subsec:blowup} we describe
the one-point blow-up $\bX_2$, its two blow-downs and its toric actions.
In~\S\ref{subsec:capacities} we introduce the equivariant capacity of a fixed
component and show that it governs which toric actions survive an equivariant
blow-up. In~\S\ref{subsec:stabilizers} we identify the stabilizer of an
invariant ball with the centralizer of the blown-up action.

%%%%%%%%%%%%%%%%%%%%%%%%%%%%%%%%%%%%%%%%%%%%%%%%%%%%%%%%%%%%%%%%%%%%%%%%%%%%%%%%
\subsection{Toric and circle actions on rational ruled surfaces}
\label{subsec:ruled}
%%%%%%%%%%%%%%%%%%%%%%%%%%%%%%%%%%%%%%%%%%%%%%%%%%%%%%%%%%%%%%%%%%%%%%%%%%%%%%%%

\begin{thm}[{\cite[Theorem 1]{finTori}}]\label{thm:extoric}
Any symplectic $S^1$ action on a rational surface $M$ of Euler number
$\chi(M)\leq 5$ extends to a Hamiltonian $\T$-action.\qed
\end{thm}

Given an integer $r\geq0$, the $r^{\text{th}}$ Hirzebruch surface $W_r$ is the
complex submanifold of $\CP^1\times\CP^2$ defined by $x_1^ry_2-x_2^ry_1=0$. The
first projection makes it a $\CP^1$-bundle over $\CP^1$, diffeomorphic to $\SSS$
when $r$ is even and to $\CCC$ when $r$ is odd, and every complex structure on
$\SSS$ or $\CCC$ is biholomorphic to some $W_r$, see~\cite{Qin}. We write
$\epsilon_r=0$ when $r$ is even and $\epsilon_r=1$ when $r$ is odd.

Throughout, symplectic forms on $\SSS$ and $\CCC$ are normalized so that the
fibre class $F$ has area $1$ and so that
\[
\oml(B)=\lambda\geq1 \quad\text{on } \SSS,
\qquad
\oml(B_{+1})=\lambda>1 \quad\text{on } \CCC,
\]
where $B$ is the class of a section of self-intersection $0$ and $B_{\pm1}$ are
the classes of the sections of self-intersection $\pm1$, so that
$F=B_{+1}-B_{-1}$. By the classification of Lalonde--McDuff, every symplectic
form on $\SSS$ or $\CCC$ is, after rescaling, diffeomorphic to exactly one such
$\oml$.

The torus $\T$ acts on $\CP^1\times\CP^2$ by
\[
(u,v)\cdot\bigl([x_1,x_2],[y_1,y_2,y_3]\bigr)
=\bigl([ux_1,x_2],[u^ry_1,y_2,vy_3]\bigr),
\]
preserving $W_r$ and its K\"ahler structure. We write $\T_r$ for the induced
toric action on $(\SSS,\oml)$ or $(\CCC,\oml)$. Its moment polygon is, up to
$\mathrm{AGL}(2,\Z)$,
\begin{equation}\label{eq:polygon}
\Delta_r=\Bigl\{(x,y)\ :\ 0\leq y\leq1,\ x\geq0,\ x+ry\leq\Lambda\Bigr\},
\qquad \Lambda=\lambda+\tfrac{r-\epsilon_r}{2},
\end{equation}
with vertices $P=(0,0)$, $Q=(0,1)$, $R=(\Lambda-r,1)$ and $S=(\Lambda,0)$, see
Figures~\ref{fig:even hirz} and~\ref{fig:Odd hirz}. The two edges $PQ$ and $RS$
are fibres of area $1$, the edge $PS$ is the section of self-intersection $+r$
of area $\Lambda$, and the edge $QR$ is the section of self-intersection $-r$ of
area
\begin{equation}\label{eq:negative section}
\om_\lambda(QR)=\lambda-\tfrac{r+\epsilon_r}{2}.
\end{equation}
In particular $\T_r$ acts on $(M,\oml)$ if and only if $2\lambda>r+\epsilon_r$.

\begin{figure}[H]
\centering       
\begin{tikzpicture}
\node[left] at (0,2) {$Q=(0,1)$};
\node[left] at (0,0) {$P=(0,0)$};
\node[right] at (4,2) {$R= (\lambda - k ,1)$};
\node[right] at (6,0) {$S=(\lambda + k ,0)$};
\node[above] at (2,2) {$B-kF$};
\node[right] at (5.15,1) {$F$};
\node[left] at (0,1) {$F$};
\node[below] at (3,0) {$B+ kF$};
\draw (0,2) -- (4,2) ;
\draw (0,0) -- (0,2) ;
\draw (0,0) -- (6,0) ;
\draw (4,2) -- (6,0) ;
\end{tikzpicture}
\caption{An even Hirzebruch polygon. The labels along the edges refer to the homology classes of the  $\T_r$ invariant spheres in $\SSS$, and the vertices $P$,$Q$,$R$,$S$ are the fixed points of the torus action.}
\label{fig:even hirz}
\end{figure}

\noindent Similarly, when $r=2k+1$ is odd, the moment map image is given in Figure~\ref{fig:Odd Hirzebruch}. We denote by $F$ the class of a fibre, while we denote by $B_{k}$ the class of a section of self-intersection $k$. Note that $F=B_{+1}-B_{-1}$. With our conventions, $\oml(F)=1$, $\oml(B_{+1})=\lambda$, and $\oml(B_{-1})=\lambda-1$. 
\begin{figure}[H]
\centering   
\begin{tikzpicture}
\node[left] at (0,2) {$Q=(0,1)$};
\node[left] at (0,0) {$P=(0,0)$};
\node[right] at (4,2) {$R= (\lambda - (k+1) ,1)$};
\node[right] at (6,0) {$S=(\lambda + k ,0)$};
\node[above] at (2,2) {$B_{+1}-(k+1)F$};
\node[right] at (5.15,1) {$F$};
\node[left] at (0,1) {$F$};
\node[below] at (3,0) {$B_{+1}+kF$};
\draw (0,2) -- (4,2) ;
\draw (0,0) -- (0,2) ;
\draw (0,0) -- (6,0) ;
\draw (4,2) -- (6,0) ;
\end{tikzpicture}
\caption{An odd Hirzebruch polygon. The labels along the edges refer to the homology classes of the  $\T_r$ invariant spheres in $\CCC$, and the vertices $P$,$Q$,$R$,$S$ are the fixed points of the torus action.}\label{fig:Odd Hirzebruch}
\label{fig:Odd hirz}
\end{figure}

\begin{lemma}[{\cite[Lemma 2.15]{ChPin-Memoirs}}]\label{lemma:list of tori}
Let $\ell$ be the largest integer strictly smaller than $\lambda$. Up to
equivariant symplectomorphism, the toric actions on $(\SSS,\oml)$ are the
$\T_{2k}$ with $0\leq k\leq\ell$, and the toric actions on $(\CCC,\oml)$ are the
$\T_{2k+1}$ with $0\leq k\leq\ell-1$.\qed
\end{lemma}

Given $(a,b)\in\Z^2$ with $\gcd(a,b)=1$, we write $S^1(a,b;r)\subset\T_r$ for the
subcircle $t\mapsto(t^a,t^b)$. Replacing $(a,b)$ by $(-a,-b)$ gives the same
subgroup, so we normalise $a\geq0$, and $b=1$ when $a=0$. The circle
$S^1(a,b;r)$ fixes pointwise the sphere over a facet of $\Delta_r$ exactly when
$(a,b)$ is parallel to the inward normal of that facet, that is exactly when
\[
a=0,\qquad\text{or}\qquad a=1 \text{ and } b\in\{0,r\}.
\]
Outside these cases the fixed point set consists of the four vertices. By
Theorem~\ref{thm:extoric} every Hamiltonian circle action on $(\SSS,\oml)$ or
$(\CCC,\oml)$ is of this form. The possible identifications between them are
described in~\cite[\S2.3]{ChPin-Memoirs}.

\begin{prop}[{\cite[Proposition 2.21]{ChPin-Memoirs}}]\label{prop:two tori}
The action $S^1(a,b;r)$ extends to $\T_r$ and to no other toric action, unless
$a=1$ and $b\notin\{0,r\}$, in which case it extends also to $\T_{n}$ with
$n=|2b-r|$, provided $2\lambda>n+\epsilon_n$, and to no other.\qed
\end{prop}

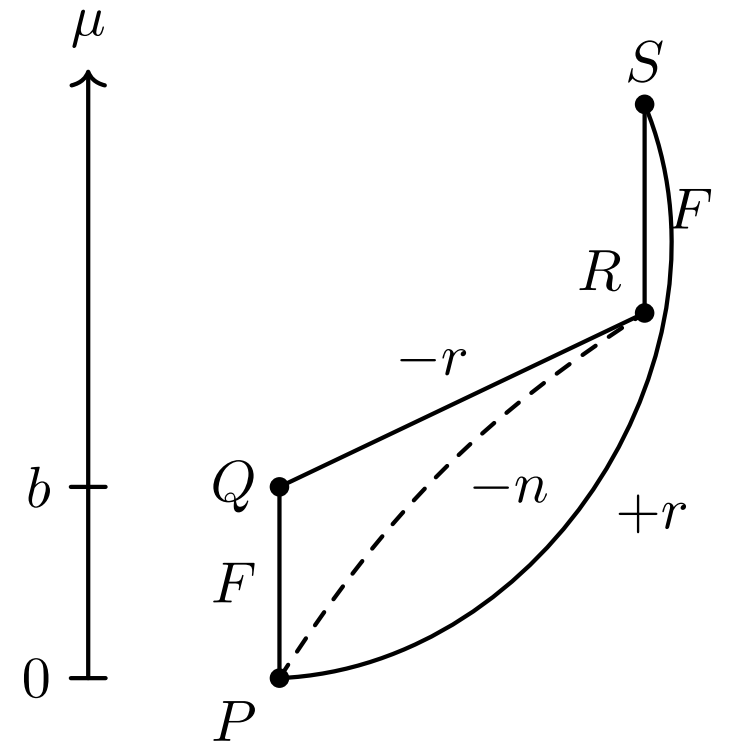
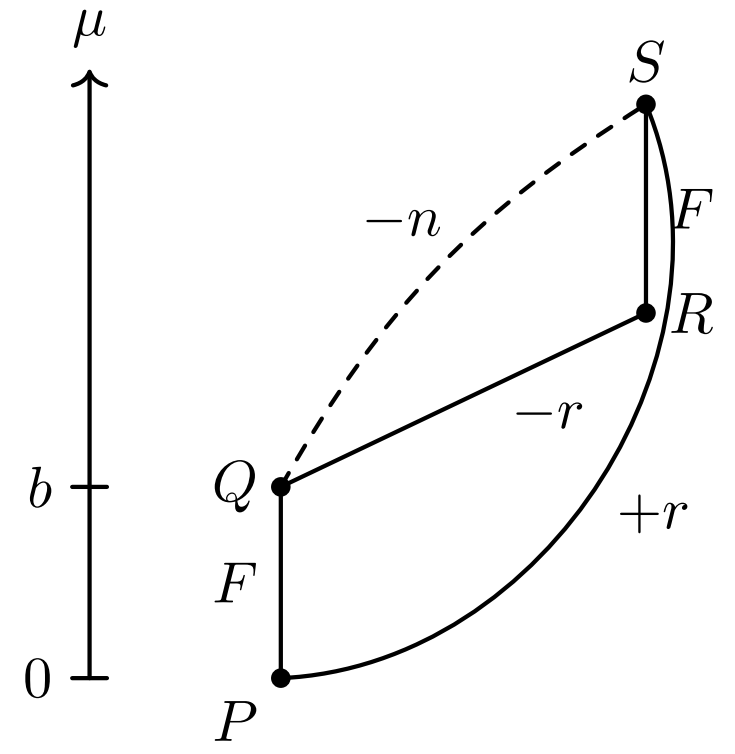
\begin{figure}[H]
\centering
\begin{subfigure}[b]{0.46\textwidth}
\centering
\begin{tikzpicture}[scale=1.15, thick, line cap=round, line join=round,
  every node/.style={font=\small}]
  \coordinate (P) at (0,0);
  \coordinate (Q) at (0,1.1);
  \coordinate (R) at (2.1,2.1);
  \coordinate (S) at (2.1,3.3);
  \draw (P) -- (Q) node[midway, left] {$F$};
  \draw (Q) -- (R) node[pos=0.55, above left=-3pt] {$-r$};
  \draw (R) -- (S) node[midway, right] {$F$};
  \draw (P) to[bend right=55] node[pos=0.5, below right=-2pt] {$+r$} (S);
  \draw[dashed] (P) to[bend left=12] node[pos=0.55, below right=-3pt] {$-n$} (R);
  \foreach \p/\l/\pos in {P/P/below left, Q/Q/left, R/R/above left, S/S/above}
    \fill (\p) circle (1.6pt) node[\pos] {$\l$};
  \draw[->] (-1.1,0) -- (-1.1,3.5) node[above] {$\mu$};
  \draw (-1.2,0) -- (-1.0,0) node[left=6pt] {$0$};
  \draw (-1.2,1.1) -- (-1.0,1.1) node[left=6pt] {$b$};
\end{tikzpicture}
\caption{$2b<r$: the sphere joins $P$ and $R$}
\end{subfigure}
\hfill
\begin{subfigure}[b]{0.46\textwidth}
\centering
\begin{tikzpicture}[scale=1.15, thick, line cap=round, line join=round,
  every node/.style={font=\small}]
  \coordinate (P) at (0,0);
  \coordinate (Q) at (0,1.1);
  \coordinate (R) at (2.1,2.1);
  \coordinate (S) at (2.1,3.3);
  \draw (P) -- (Q) node[midway, left] {$F$};
  \draw (Q) -- (R) node[pos=0.6, below right=-3pt] {$-r$};
  \draw (R) -- (S) node[midway, right] {$F$};
  \draw (P) to[bend right=55] node[pos=0.5, below right=-2pt] {$+r$} (S);
  \draw[dashed] (Q) to[bend left=14] node[pos=0.55, above left=-3pt] {$-n$} (S);
  \foreach \p/\l/\pos in {P/P/below left, Q/Q/left, R/R/right, S/S/above}
    \fill (\p) circle (1.6pt) node[\pos] {$\l$};
  \draw[->] (-1.1,0) -- (-1.1,3.5) node[above] {$\mu$};
  \draw (-1.2,0) -- (-1.0,0) node[left=6pt] {$0$};
  \draw (-1.2,1.1) -- (-1.0,1.1) node[left=6pt] {$b$};
\end{tikzpicture}
\caption{$2b>r$: the sphere joins $Q$ and $S$}
\end{subfigure}
\caption{The two regimes of Lemma~\ref{lemma:second torus}, for $a=1$ and
$0<b<r$, drawn as graphs: the vertices are the fixed points, placed at their
moment values, and the edges are the invariant spheres joining them. The solid
edges are the spheres over the edges of $\Delta_r$, labelled by their
self-intersections, or by $F$ for the two fibres; the dashed edge is the sphere
of self-intersection $-n$, $n=|2b-r|$, which is the negative section of the
second torus $\T_n$. Its two ends are the fixed points corresponding to the
vertices $Q'$ and $R'$ of $\Delta_n$. For $b<0$ or $b>r$ the moment values
reorder --- the minimum becomes $Q$, or the maximum becomes $R$ --- but the
pairing is unchanged.}
\label{fig:second torus}
\end{figure}

\begin{lemma}[The second torus]\label{lemma:second torus}
Let $a=1$ and $b\notin\{0,r\}$, put $n=|2b-r|$ and assume $2\lambda>n+\epsilon_n$,
so that $\T_n$ acts on $(M,\oml)$. Then $S^1(1,b;r)$ and $S^1(1,b';n)$ are
conjugate subgroups, with
\[
b'=-b \quad\text{if } 2b<r,\qquad b'=b \quad\text{if } 2b>r,
\]
and the induced correspondence between the vertices of $\Delta_r$ and those of
$\Delta_n$, written $P',Q',R',S'$, is, see Figure~\ref{fig:second torus},
\[
\begin{array}{ll}
2b<r:& P\mapsto Q',\quad Q\mapsto P',\quad R\mapsto R',\quad S\mapsto S',\\[2pt]
2b>r:& P\mapsto P',\quad Q\mapsto Q',\quad R\mapsto S',\quad S\mapsto R'.
\end{array}
\]
\end{lemma}
\begin{proof}
A Hamiltonian circle action on a closed symplectic $4$-manifold with isolated
fixed points is determined, up to equivariant symplectomorphism, by the moment
values and the isotropy weights at its fixed points, see~\cite[\S4]{karshon}. For
$u=(1,b)$ on $\Delta_r$ these are
\[
\begin{array}{c|cccc}
 & P & Q & R & S\\\hline
\mu & 0 & b & \Lambda-r+b & \Lambda\\
\text{weights} & (b,1) & (-b,1) & (-1,r-b) & (-1,b-r)
\end{array}
\]
and the same formulas hold on $\Delta_n$ for $u'=(1,b')$, with $n$ and
$\Lambda'=\lambda+\frac{n-\epsilon_n}{2}$ in place of $r$ and $\Lambda$. Since
$n\equiv r \pmod 2$, one has $\epsilon_n=\epsilon_r$ and
\[
\Lambda'=\Lambda-b \quad\text{if } 2b<r,\qquad
\Lambda'=\Lambda+b-r \quad\text{if } 2b>r .
\]
Substituting the stated value of $b'$ gives, when $2b<r$, so that $n=r-2b$,
$b'=-b$ and $\Lambda'=\Lambda-b$,
\[
\begin{array}{c|cccc}
 & P' & Q' & R' & S'\\\hline
\mu'+b & b & 0 & \Lambda-r+b & \Lambda\\
\text{weights} & (-b,1) & (b,1) & (-1,r-b) & (-1,b-r)
\end{array}
\]
and, when $2b>r$, so that $n=2b-r$, $b'=b$ and $\Lambda'=\Lambda+b-r$,
\[
\begin{array}{c|cccc}
 & P' & Q' & R' & S'\\\hline
\mu' & 0 & b & \Lambda & \Lambda-r+b\\
\text{weights} & (b,1) & (-b,1) & (-1,b-r) & (-1,r-b).
\end{array}
\]
In the first case the moment map has been translated by $b$, which is harmless.
Comparing each of these tables with the one above, entry by entry, gives the two
correspondences of the statement.
\end{proof}

We record the weights, in the convention of~\cite[Table 1]{ChPin-Memoirs}: at a
vertex the two weight vectors are the primitive directions of the two edges
leaving it, and the weights of $S^1(a,b;r)$ are their pairings with $(a,b)$.
Thus the weight vectors are
\[
P:\ \{(1,0),(0,1)\},\quad
Q:\ \{(1,0),(0,-1)\},\quad
R:\ \{(r,-1),(-1,0)\},\quad
S:\ \{(-1,0),(-r,1)\}.
\]
An equivariant blow-up at a vertex with weight vectors $w_1,w_2$ replaces it by
two vertices, on the two edges, with weight vectors $\{w_1,w_2-w_1\}$ and
$\{w_1-w_2,w_2\}$, see~\cite[\S6.1]{karshon}. The weights of $S^1(a,b;r)$ at the
two new fixed points are therefore as follows, where $p'$ denotes the new vertex
on the first edge listed above and $p''$ the one on the second.

\begin{table}[H]
\begin{center}
\begin{tabular}{|c||c|c|}
\hline
blow-up at & weights at $p'$ & weights at $p''$\\
\hline\hline
$P$ & $\{a,\,b-a\}$ & $\{b,\,a-b\}$ \rule{0pt}{10pt}\\
\hline
$Q$ & $\{a,\,-a-b\}$ & $\{-b,\,a+b\}$ \rule{0pt}{10pt}\\
\hline
$R$ & $\{ar-b,\,-a-ar+b\}$ & $\{-a,\,a+ar-b\}$ \rule{0pt}{10pt}\\
\hline
$S$ & $\{-a,\,a-ar+b\}$ & $\{-ar+b,\,-a+ar-b\}$ \rule{0pt}{10pt}\\
\hline
\end{tabular}
\end{center}
\caption{Weights of $S^1(a,b;r)$ at the two fixed points created by an
equivariant blow-up at $P$, $Q$, $R$ or $S$.}
\label{table:weights}
\end{table}

\begin{lemma}[The label of the exceptional divisor]\label{lemma:divisor label}
Let $p$ be a vertex of $\Delta_r$ with edge directions $v_1,v_2$, and let $Z$ be
the exceptional divisor of the equivariant blow-up of $S^1(a,b;r)$ at $p$. Then
$Z$ is the sphere over the chop edge, of primitive direction $v_2-v_1$, and the
circle acts on it with weight $\langle(a,b),v_2-v_1\rangle$, that is
\[
b-a \ \text{ at } P,\qquad
-(a+b) \ \text{ at } Q,\qquad
b-a(r+1) \ \text{ at } R,\qquad
b-a(r-1) \ \text{ at } S .
\]
We call the absolute value of that weight the \emph{label} of $Z$. If the label
is $0$, then $Z$ is pointwise fixed. If it is $k\geq2$, then $Z$ is a component
of the fixed point set of the subgroup $\Z_k\subset S^1$. In either case $Z$ is
invariant under every torus containing the circle, and is holomorphic for every
invariant almost complex structure, by Lemma~\ref{LR}. If the label is $1$, then
$Z$ is a free sphere and does not appear on the decorated graph of the blown-up
action, which records only the fixed components and the $\Z_k$-spheres with
$k\geq2$. Nothing then forces $Z$ to be invariant under a torus other than the
one produced by the equivariant blow-up.
\end{lemma}
\begin{proof}
The chop replaces $p$ by the two points $p+cv_1$ and $p+cv_2$, so the new edge
has primitive direction $v_2-v_1$, and the weight of the circle along an
invariant sphere is the pairing of $(a,b)$ with the primitive direction of the
corresponding edge. The four values follow from the edge directions listed
before Table~\ref{table:weights}. A sphere of label $0$ is a component of the
fixed point set, and one of label $k\geq2$ is a component of the fixed point set
of $\Z_k$. Both are therefore invariant under any torus containing $S^1$, and
Lemma~\ref{LR} applies. A sphere of label $1$ is free, and free spheres are not
part of the decorated graph: they are the spheres one adds to it when extending
the circle action to a torus.
\end{proof}

%%%%%%%%%%%%%%%%%%%%%%%%%%%%%%%%%%%%%%%%%%%%%%%%%%%%%%%%%%%%%%%%%%%%%%%%%%%%%%%%
\subsection{The one-point blow-up and its two blow-downs}
\label{subsec:blowup}
%%%%%%%%%%%%%%%%%%%%%%%%%%%%%%%%%%%%%%%%%%%%%%%%%%%%%%%%%%%%%%%%%%%%%%%%%%%%%%%%

Write $\bX_2=\CP^2\#2\overline{\CP^2}$ and let $L,E_1,E_2$ be the standard basis
of $H_2(\bX_2;\Z)$, with $L\cdot L=1$, $E_i\cdot E_i=-1$ and $L\cdot E_i=0$. The
exceptional classes of $\bX_2$ are $E_1$, $E_2$ and the characteristic class
\[
E_0:=L-E_1-E_2,
\]
with $E_1\cdot E_2=0$ and $E_0\cdot E_i=1$. Blowing down $E_0$ gives $\SSS$,
while blowing down $E_1$ or $E_2$ gives $\CCC$. It follows
from~\cite{liliu} that every symplectic form on $\bX_2$ is, up to
diffeomorphism and rescaling, obtained by blowing up a ball in one of these
three ways.

\begin{lemma}[The two blow-downs]\label{lemma:two blow-downs}
Let $0<c<1\leq\lambda$. The blow-up of $(\SSS,\oml)$ along a ball of capacity
$c$, with exceptional class $E:=E_0$, coincides with the blow-up of
$(\CCC,\om_{\lambda+1-c})$ along a ball of capacity $1-c$, with exceptional
class $F-E$; the fibre classes agree. Conversely, if $c\leq\lambda-1$, the
blow-up of $(\CCC,\oml)$ along a ball of capacity $c$ is
$(\bX_2,\om_{\lambda-c,\,1-c})$, with exceptional class $F-E$ and
$B_{-1}=B-E$.
\end{lemma}
\begin{proof}
In the basis $L,E_1,E_2$ one has $B=L-E_2$, $F=L-E_1$ and $E_0=L-E_1-E_2$, while
$B_{+1}=L$ and $B_{-1}=E_1$ for the blow-down of $E_2$. Comparing areas gives
both statements as the hypothesis $\lambda-c\geq1$ is
the normalisation $\om(B)\geq\om(F)$ for the resulting form.
\end{proof}

We write $(\bX_2,\omlc)$ for the blow-up of $(\SSS,\oml)$ at a ball of capacity
$c$, so that
\[
\omlc(B)=\lambda,\qquad \omlc(F)=1,\qquad \omlc(E)=c .
\]
For $i\geq0$ set
\begin{equation}\label{eq:Di}
D_0=B-E,\qquad D_{2k-1}=B-kF,\qquad D_{2k}=B-kF-E \quad(k\geq1),
\end{equation}
so that $D_i\cdot D_i=-(i+1)$. The indexing is that of the codimension of the
corresponding stratum in~\S\ref{section:Strataofinvariantacs} and is the
opposite of the one of~\cite[Lemma 2.5]{P-com}. Their areas,
\[
\omlc(D_0)=\lambda-c,\qquad
\omlc(D_{2k-1})=\lambda-k,\qquad
\omlc(D_{2k})=\lambda-k-c,
\]
alternate between expressions with and without $c$.

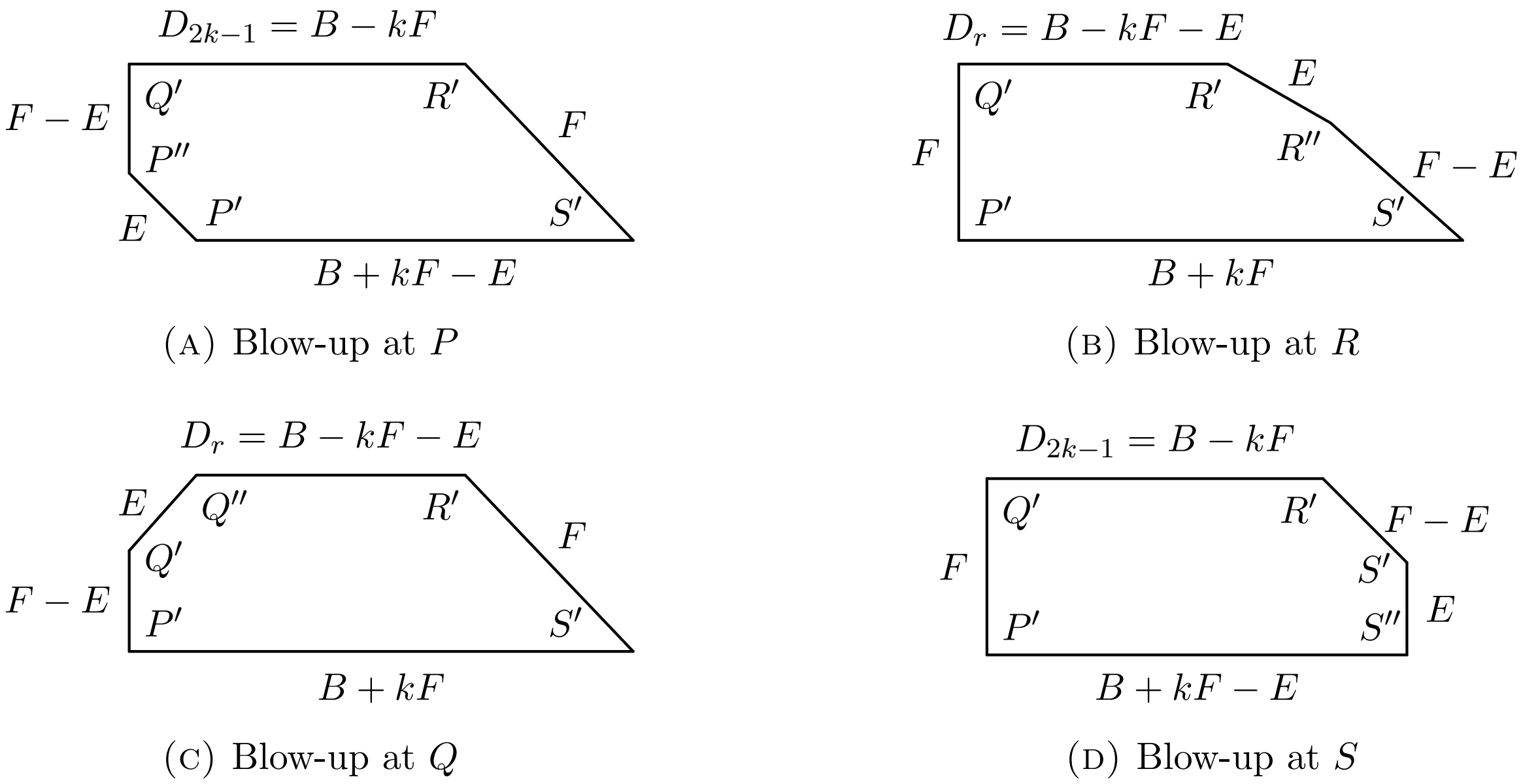
\begin{figure}[htbp]
\centering
\begin{subfigure}[b]{0.45\textwidth}
\centering
%
% Blow-up at P
%
\begin{tikzpicture}[
    scale=0.85,
    line cap=round,
    line join=round,
    >=stealth,
    thick,                             % thicker lines for visibility
    every node/.style={font=\small}    % smaller, uniform label font
]

  % === Coordinates ===
  \coordinate (Pp) at (0.8,0);
  \coordinate (Ppp) at (0,0.8);
  \coordinate (Qp) at (0,2.1);
  \coordinate (Rp) at (4,2.1);
  \coordinate (Sp) at (6,0);

  % === Edges ===
  \draw (Ppp) -- (Qp);
  \draw (Qp)  -- (Rp);
  \draw (Pp)  -- (Ppp);
  \draw (Pp)  -- (Sp);
  \draw (Rp)  -- (Sp);

  % === Node labels ===
  \node[above right=0pt and -2pt of Pp]{$P'$}; 
  \node[below right=0pt of Qp]{$Q'$};
  \node[below left=0pt and -3pt of Rp]{$R'$};
  \node[above left=0pt and 10pt of Sp]{$S'$};  
  \node[above right=-4pt and 0pt of Ppp]{$P''$};

  % === Edge labels (aligned midpoints) ===
  \node[above=1pt of $(Qp)!0.5!(Rp)$]{$D_{2k-1} = B - kF$};
  \node[above right=0pt and -2pt of $(Rp)!0.5!(Sp)$]{$F$};
  \node[left=1pt of $(Ppp)!0.5!(Qp)$]{$F-E$};
  \node[below=1pt of $(Pp)!0.5!(Sp)$]{$B + kF-E$};
  \node[below left=-2pt and 0pt of $(Pp)!0.5!(Ppp)$]{$E$};  
\end{tikzpicture}
\caption{Blow-up at $P$}
\end{subfigure}
\hfill
%  
% Blow-up at R  
%
\begin{subfigure}[b]{0.45\textwidth}
\centering
\begin{tikzpicture}[
    scale=0.85,
    line cap=round,
    line join=round,
    >=stealth,
    thick,                             % thicker lines for visibility
    every node/.style={font=\small}    % smaller, uniform label font
]

  % === Coordinates ===
  \coordinate (Pp) at (0,0);
  \coordinate (Qp) at (0,2.1);
  \coordinate (Rp) at (3.2,2.1);
  \coordinate (Rpp) at (4.425,1.4);  
  \coordinate (Sp) at (6,0);

  % === Edges ===
  \draw (Pp) -- (Qp);
  \draw (Qp)  -- (Rp);
  \draw (Rp)  -- (Rpp);
  \draw (Pp)  -- (Sp);
  \draw (Rpp)  -- (Sp);

    % === Node labels ===
  \node[above right=0pt and 0pt of Pp]{$P'$}; 
  \node[below right=0pt and 0pt of Qp]{$Q'$};
  \node[below left=0pt and -3pt of Rp]{$R'$};
  \node[above left=0pt and 12pt of Sp]{$S'$};  
  \node[below left=-2pt and -2pt of Rpp]{$R''$};

  % === Edge labels (aligned midpoints) ===
  \node[above=1pt of $(Qp)!0.5!(Rp)$]{$D_r = B - kF-E$};
  \node[above right=-4pt and 0pt of $(Rpp)!0.5!(Sp)$]{$F-E$};
  \node[left=1pt of $(Pp)!0.5!(Qp)$]{$F$};
  \node[below=1pt of $(Pp)!0.5!(Sp)$]{$B + kF$};
  \node[above right=-2pt and -2pt of $(Rp)!0.5!(Rpp)$]{$E$};

\end{tikzpicture}
\caption{Blow-up at $R$}
\end{subfigure}

    \vspace{1.5em} % Espace vertical entre les lignes
%
% Blow-up at Q
%
\begin{subfigure}[b]{0.45\textwidth}
\centering
\begin{tikzpicture}[
    scale=0.85,
    line cap=round,
    line join=round,
    >=stealth,
    thick,                             % thicker lines for visibility
    every node/.style={font=\small}    % smaller, uniform label font
]

  % === Coordinates ===
  \coordinate (Pp) at (0,0);
  \coordinate (Qp) at (0,1.2);
  \coordinate (Qpp) at (0.8,2.1);
  \coordinate (Rp) at (4,2.1);
  \coordinate (Sp) at (6,0);

  % === Edges ===
  \draw (Qpp) -- (Rp);
  \draw (Qp)  -- (Qpp);
  \draw (Pp)  -- (Qp);
  \draw (Pp)  -- (Sp);
  \draw (Rp)  -- (Sp);

  % === Node labels ===
  \node[above right=0pt and 0pt of Pp]{$P'$}; 
  \node[below right=-7pt and 0pt of Qp]{$Q'$};
  \node[below left=0pt and -3pt of Rp]{$R'$};
  \node[above left=0pt and 10pt of Sp]{$S'$};  
  \node[below right=0pt and -3pt of Qpp]{$Q''$};

  % === Edge labels (aligned midpoints) ===
  \node[above=2pt of $(Qpp)!0.5!(Rp)$] {$D_r = B - kF - E$};
  \node[above right=0pt and -2pt of $(Rp)!0.5!(Sp)$]{$F$};
  \node[left=1pt of $(Pp)!0.5!(Qp)$]   {$F-E$};
  \node[below=2pt of $(Pp)!0.5!(Sp)$]  {$B + kF$};
  \node[above left= -5pt and 0pt of $(Qp)!0.5!(Qpp)$]{$E$};  
\end{tikzpicture}
\caption{Blow-up at $Q$}
\end{subfigure}
\hfill
%
% Blow-up at S
%
\begin{subfigure}[b]{0.45\textwidth}
\centering
\begin{tikzpicture}[
    scale=0.85,
    line cap=round,
    line join=round,
    >=stealth,
    thick,                             % thicker lines for visibility
    every node/.style={font=\small}    % smaller, uniform label font
]

  % === Coordinates ===
  \coordinate (Pp) at (0,0);
  \coordinate (Qp) at (0,2.1);
  \coordinate (Rp) at (4,2.1);
  \coordinate (Sp) at (5,1.1);
  \coordinate (Spp) at (5,0);  

  % === Edges ===
  \draw (Pp) -- (Qp);
  \draw (Qp)  -- (Rp);
  \draw (Rp)  -- (Sp);
  \draw (Pp)  -- (Spp);
  \draw (Spp)  -- (Sp);

    % === Node labels ===
  \node[above right=0pt and 0pt of Pp]{$P'$}; 
  \node[below right=0pt and 0pt of Qp]{$Q'$};
  \node[below left=0pt and -3pt of Rp]{$R'$};
  \node[above left=-10pt and 0pt of Sp]{$S'$};  
  \node[above left=0pt and -3pt of Spp]{$S''$};

  % === Edge labels (aligned midpoints) ===
  \node[above=1pt of $(Qp)!0.5!(Rp)$] {$D_{2k-1} = B - kF$};
  \node[right=1pt of $(Rp)!0.5!(Sp)$]  {$F-E$};
  \node[left=1pt of $(Pp)!0.5!(Qp)$]   {$F$};
  \node[below=1pt of $(Pp)!0.5!(Spp)$]  {$B + kF-E$};
  \node[right=1pt of $(Sp)!0.5!(Spp)$]   {$E$}; 
\end{tikzpicture}
\caption{Blow-up at $S$}
\end{subfigure}

\caption{Delzant polygons of toric actions on $\bX_2$ with homology and fixed points labels. The homology labels refer to the identification of $\bX_2$ with the blow-up $\SSSC$.}
\label{fig:toric blow-ups}
\end{figure}

Every toric action on $\bX_2$ is an equivariant blow-up of a toric action on a
rational ruled surface at a vertex of its polygon, and the four possibilities
are drawn in Figure~\ref{fig:toric blow-ups}. We write $\ttT_i$ for the toric
action on $\bX_2$ whose invariant curve of most negative self-intersection
lies in the class $D_i$.

\begin{lemma}[Dictionary]\label{lemma:dictionary}
For $i\geq1$,
\[
\ttT_i\ \cong\ \T_{i+1}(P)\ \cong\ \T_{i+1}(S)\ \cong\ \T_i(Q)\ \cong\ \T_i(R),
\]
where $\T_m(p)$ denotes the equivariant blow-up of $\T_m$ at the vertex $p$. For
$i=0$ all four blow-ups of $\T_0$, and the blow-ups of $\T_1$ at $P$ and $S$,
give $\ttT_0$. The action $\ttT_i$ exists on $(\bX_2,\omlc)$ if and only if
$\omlc(D_i)>0$.
\end{lemma}
\begin{proof}
Blowing up at a vertex of the negative section of $\Delta_m$ lowers its
self-intersection by one and its area by $c$, producing a curve in the class
$D_m$. Blowing up at a vertex of the positive section leaves the negative
section untouched, so the most negative curve is in the class $D_{m-1}$. The
last assertion is the positivity of the area of that curve.
\end{proof}

\begin{defn}\label{defn:indecomposable}
Let a compact group $G$ act symplectically on a symplectic $4$-manifold
$(M,\om)$. An exceptional class $A\in H_2(M;\Z)$ is \emph{indecomposable} if,
for every $\om$-tame almost complex structure $J$, it is represented by an
embedded $J$-holomorphic sphere, and \emph{equivariantly indecomposable} if
this holds for every invariant $\om$-tame $J$.
\end{defn}

\begin{lemma}[{\cite[Lemmata 2.2 and 2.5]{P-com}}]\label{lemma:indecomposable}
On $(\bX_2,\omlc)$ the classes $E$ and $F-E$ are indecomposable, while $B-E$ is
indecomposable if and only if $\lambda=1$.\qed
\end{lemma}

\begin{cor}\label{cor:blow down}
Every Hamiltonian toric or circle action on $(\bX_2,\omlc)$ can be
equivariantly blown down along an invariant sphere in the class $E$ to an action
on $(\SSS,\oml)$,
and along an invariant sphere in the class $F-E$ to an action on $(\CCC,
\om_{\lambda+1-c})$.
\end{cor}
\begin{proof}
By Lemma~\ref{lemma:indecomposable}, the classes $E$ and $F-E$ are
represented, for every invariant compatible almost complex structure, by a
unique embedded $J$-holomorphic sphere, which is therefore invariant. Blowing
down along it gives the two assertions, by Lemma~\ref{lemma:two blow-downs}.
\end{proof}

By Corollary~\ref{cor:blow down} and Theorem~\ref{thm:extoric}, every
Hamiltonian circle action on $\bX_2$ is the equivariant blow-up of some
$\alpha:=S^1(a,b;r)$ at a component $\fixedcomp$ of its fixed point set. We denote it by
$S^1(a,b;r;\fixedcomp)$, and by $\widehat\alpha_{\fixedcomp,c}$ when we wish to
display the capacity. The same action may arise from several such data, see
Lemma~\ref{lemma:dictionary}, and when $\fixedcomp$ is a fixed surface the
blow-up does not depend on the chosen point of $\fixedcomp$, by Lemma~\ref{Lemma:Isotopywithinconnectedcomponent}.

\begin{Remark}\label{rem:ball independence}
The blow-up $\widehat\alpha_{\fixedcomp,c}$ does not depend on the chosen
invariant ball either. Indeed, a Hamiltonian circle action on a closed
symplectic $4$-manifold is determined, up to equivariant symplectomorphism, by
its graph, see~\cite[\S4]{karshon}, and the graph of an equivariant blow-up is
obtained from that of the action by the corner-chopping operation
of~\cite[\S6.1]{karshon}, which depends only on $\fixedcomp$ and on $c$. Two
invariant balls of capacity $c$ centred on $\fixedcomp$ therefore have
equivariantly symplectomorphic blow-ups.
\end{Remark}

%%%%%%%%%%%%%%%%%%%%%%%%%%%%%%%%%%%%%%%%%%%%%%%%%%%%%%%%%%%%%%%%%%%%%%%%%%%%%%%%
\subsection{Equivariant capacities}
\label{subsec:capacities}
%%%%%%%%%%%%%%%%%%%%%%%%%%%%%%%%%%%%%%%%%%%%%%%%%%%%%%%%%%%%%%%%%%%%%%%%%%%%%%%%

\begin{defn}\label{defn:capacity}
Let a compact group $G$ act symplectically on $(M,\om)$ and let $\fixedcomp$ be
a component of its fixed point set. The \emph{equivariant capacity of
$\fixedcomp$} is
\[
c^G(M,\fixedcomp)=\sup\bigl\{c\ :\ \text{there is a $G$-equivariant symplectic
embedding of } B^4(c) \text{ centred on }\fixedcomp\bigr\},
\]
the action on $B^4(c)$ being the isotropy representation.
\end{defn}

\begin{lemma}\label{lemma:capacity basics}
Let $H\subseteq G$ be compact groups acting on $(M,\om)$ and let $\fixedcomp$ be
a component of the fixed point set of $G$. Then
\[
c^{G}(M,\fixedcomp)\ \leq\ c^{H}(M,\fixedcomp)\ \leq\ c^{\{1\}}(M,\fixedcomp),
\]
the last quantity being the Gromov width of $M$ at a point of $\fixedcomp$.
Moreover, if $\T$ is a toric action and $p$ a fixed point of $\T$, then
$c^{\T}(M,p)$ is the smaller of the areas of the two invariant spheres through
$p$, that is the smaller of the lengths of the two edges of the moment polygon
meeting at $p$.
\end{lemma}
\begin{proof}
A $G$-equivariant ball is in particular $H$-equivariant, which gives the
inequalities. For a toric action, the equivariant balls centered at a vertex are
the standard triangles of side length equalling the capacity of the ball embedding $c$ in the moment polytope. This follows because, by the equivariant Darboux theorem such
a ball is standard near the vertex, so that its image under the moment map is
the simplex of size $c$ based there, and by convexity that simplex lies in the
polygon if and only if $c$ is smaller than the lengths of the two edges meeting
at the vertex (cf.~\cite{KKP07}). Finally, those two edges are the images of the
two invariant spheres through the vertex, so their lengths are the areas of
those spheres.
\end{proof}

\begin{prop}\label{prop:capacity rule}
Let $\alpha$ be a Hamiltonian circle action on a rational ruled surface
$(M,\oml)$, let $\fixedcomp$ be a component of its fixed point set and let $c$ be
the capacity of an equivariant ball centred on $\fixedcomp$. Write $\T(p)$ for
the equivariant blow-up of a toric action $\T$ at a fixed point $p$ of $\T$.
Then
\[
\Ext\bigl(\widehat\alpha_{\fixedcomp,c}\bigr)
\supseteq\bigl\{\,\T(p)\ :\ \T\in\Ext(\alpha),\ p\in\fixedcomp\ \text{a fixed
point of }\T,\ c<c^{\T}(M,p)\,\bigr\},
\]
and equality holds whenever the exceptional divisor is invariant under every
toric extension of $\widehat\alpha_{\fixedcomp,c}$. Furthermore, the following two conditions guarantee that the exceptional curve is invariant under every
toric extension of $\widehat\alpha_{\fixedcomp,c}$.
\begin{enumerate}
\item the exceptional divisor is not a free gradient sphere for the $S^1$-action ;
\item the class of the divisor is indecomposable in $(\bX_2,\tom)$, which by
Lemmas~\ref{lemma:two blow-downs} and~\ref{lemma:indecomposable} holds for
$M=\SSS$, and for $M=\CCC$ when $c\leq\lambda-1$.
\end{enumerate}
\end{prop}
\begin{proof} 
By Lemma~\ref{Lemma:Isotopywithinconnectedcomponent} and Remark~\ref{rem:ball independence}, the blow-up
$\widehat\alpha_{\fixedcomp,c}$ depends only on $\fixedcomp$ and on $c$: the
centre of the ball may be moved to any point of $\fixedcomp$, and the ball may
be replaced by any invariant ball of the same capacity. Let now
$\T\in\Ext(\alpha)$ and let $p\in\fixedcomp$ be a fixed point of $\T$ with
$c<c^{\T}(M,p)$. By Karshon's uniqueness theorem \cite{karshon} that the symplectomorphism type of the $S^1$-equivariant blow up only depends on the fixed point component and the capacity of the blow-up, hence we can now blow up along the $\T$-equivariant ball. The corner chop of size $c$ at $p$ is a $\T$-equivariant ball,
and the blow-up along it is toric with torus $\T(p)$ whose blow-up
is $\widehat\alpha_{\fixedcomp,c}$, so $\T(p)$ extends it. Conversely, assume the exceptional divisor is invariant
under every toric extension and let $\ttT$ extend
$\widehat\alpha_{\fixedcomp,c}$. Blowing the divisor down gives a toric action $\T$ on
$(M,\oml)$ extending $\alpha$,
together with a fixed point $p$ of $\T$ lying in $\fixedcomp$, at which the ball
is the corner chop of size $c$. Hence $\ttT=\T(p)$ with $\T\in\Ext(\alpha)$ and
$c<c^{\T}(M,p)$, by Lemma~\ref{lemma:capacity basics}. Finally, when (1) is satisfied, the exceptional curve is holomorphic for all invariant almost complex structures (and therby all toric extensions) by Lemma~\ref{LR}. Similarly when condition (2) is satisfied, 
the class is then represented, for every $\ttT$-invariant compatible almost
complex structure, by a unique embedded sphere, which is therefore
$\ttT$-invariant. 
\end{proof}

\begin{Remark}\label{rem:corner counterexample}
Consider $M=(\CCC,\om_\lambda)$ with
$\lambda=\tfrac32$, the action $\alpha=S^1(1,2;1)$, the ball at $P$ and
$c=\tfrac9{10}$, so that $c>\lambda-1$. Then $\Ext(\alpha)=\{\T_1\}$, since
$\T_3$ would require $\lambda>2$, and the label of the divisor at $P$ is
$|b-a|=1$. The blown-up
action has the five fixed points
\[
p'\ (c;\{1,1\}),\quad S\ (\lambda;\{-1,1\}),\quad p''\ (2c;\{2,-1\}),\quad
Q\ (2;\{1,-2\}),\quad R\ (\lambda+1;\{-1,-1\}),
\]
with a single $\Z_k$-sphere of $k\geq2$, namely $p''Q$. By
Proposition~\ref{prop:capacity rule}, one toric extension of $\widehat\alpha$ is
$\T_1(P)$, which is $\ttT_0$ in the normalisation of
Lemma~\ref{lemma:decomposable normalisation}. Consider also the toric action
whose Delzant pentagon has vertices
\[
(0,0),\quad (0.6,0.6),\quad (0.9,0.6),\quad (1.1,0.5),\quad (1.6,0).
\]
Its subcircle whose moment map is the first coordinate has the same moment
values up to the additive constant $c$, the same weights and the same unique
$\Z_2$-sphere, so by Karshon's uniqueness theorem~\cite[Theorem 4.1]{karshon}
it is the action $\widehat\alpha$. The self-intersections of the edges of this
pentagon are $(2,-1,-1,-3,0)$, so that it is the polygon of $\ttT_2$. Its edges
carry the classes $2L-E_1-E_2$, $E_1$, $L-E_1-E_2$, $2E_2-L$, $L-E_2$, of areas
$1.6$, $0.5$, $0.1$, $0.3$, $0.6$, in agreement with $\om$, and the class
$2E_2-L$ is $D_2$ in the normalisation of
Lemma~\ref{lemma:decomposable normalisation}. So $\widehat\alpha$ has two toric
extensions although $\alpha$ has one, and the inclusion of
Proposition~\ref{prop:capacity rule} is strict. In $\ttT_2$ the divisor is not an
edge. It is represented by the chain $S\,p''\,Q\,R$, with
$E_2=(2E_2-L)+(L-E_1-E_2)+E_1$ and areas $0.3+0.1+0.5=0.9=c$, see
Figure~\ref{fig:ball at P}. In particular the class of the divisor is not
equivariantly indecomposable.
\end{Remark}

\tikzset{
  fp/.style={circle,fill,inner sep=1.4pt},
  free/.style={thin},
  zk/.style={very thick,blue!70!black},
  lab/.style={font=\scriptsize,inner sep=1.5pt},
  fat/.style={draw,fill=black!15,ellipse,minimum width=1.6cm,minimum height=0.3cm,inner sep=0pt},
}

\begin{figure}[H]
\centering
\begin{subfigure}[b]{0.27\textwidth}\centering
\begin{tikzpicture}[yscale=1]
\node[fp,label={[lab]right:$p'\ (0.9)$}] (a) at (0,0) {};
\node[fp,label={[lab]right:$S\ (1.5)$}] (S) at (0.8,1.92) {};
\node[fp,label={[lab]left:$p''\ (1.8)$}] (b) at (-0.8,2.88) {};
\node[fp,label={[lab]left:$Q\ (2)$}] (Q) at (-0.8,3.52) {};
\node[fp,label={[lab]right:$R\ (2.5)$}] (R) at (0,5.12) {};
\draw[zk] (b) -- node[lab,right]{$2$} (Q);
\end{tikzpicture}
\caption{Karshon graph}
\end{subfigure}\hfill
\begin{subfigure}[b]{0.33\textwidth}\centering
\begin{tikzpicture}
\node[fp,label={[lab]below:$p'$}] (a) at (0,0) {};
\node[fp,label={[lab]right:$S$}] (S) at (0.8,1.92) {};
\node[fp,label={[lab]left:$p''$}] (b) at (-0.8,2.88) {};
\node[fp,label={[lab]left:$Q$}] (Q) at (-0.8,3.52) {};
\node[fp,label={[lab]above:$R$}] (R) at (0,5.12) {};
\draw[free] (a) -- node[lab,left]{$E_2\ (0.9)$} (b);
\draw[zk] (b) -- node[lab,right]{$2$} (Q);
\draw[free] (Q) -- node[lab,left]{$E_1\ (0.5)$} (R);
\draw[free] (a) -- node[lab,right]{$L-E_2\ (0.6)$} (S);
\draw[free] (S) -- node[lab,right]{$L-E_1\ (1)$} (R);
\end{tikzpicture}
\caption{Extension (i): $\T_1(P)=\ttT_0$}
\end{subfigure}\hfill
\begin{subfigure}[b]{0.39\textwidth}\centering
\begin{tikzpicture}
\node[fp,label={[lab]below:$p'$}] (a) at (0,0) {};
\node[fp,label={[lab]left:$S$}] (S) at (1.2,1.92) {};
\node[fp,label={[lab]left:$p''$}] (b) at (1.2,2.88) {};
\node[fp,label={[lab]left:$Q$}] (Q) at (1.2,3.52) {};
\node[fp,label={[lab]above:$R$}] (R) at (0,5.12) {};
\draw[free] (a) -- node[lab,right]{$L-E_2\ (0.6)$} (S);
\draw[free] (S) -- node[lab,right]{$2E_2-L\ (0.3)$} (b);
\draw[zk] (b) -- node[lab,right]{$2$} (Q);
\draw[free] (Q) -- node[lab,right]{$E_1\ (0.5)$} (R);
\draw[free] (a) -- node[lab,left]{$2L-E_1-E_2\ (1.6)$} (R);
\end{tikzpicture}
\caption{Extension (ii): $\ttT_2$}
\end{subfigure}
\caption{$S^1(1,2;1)$ blown up at $P$, for $\lambda=\tfrac32$ and $c=\tfrac9{10}$.
Vertices are the fixed points, placed at their moment values (shown in (a)).
Thick edges are $\Z_k$-spheres, labelled by $k$; they form Karshon's graph and
appear in every extension. Thin edges are the free invariant spheres of each
toric extension, labelled by their class in the basis $L,E_1,E_2$ and their
area. In (i) the divisor $E_2$ is the edge $p'p''$; in (ii) it is no longer an
edge and is represented by the chain $S\,p''\,Q\,R$.}
\label{fig:ball at P}
\end{figure}
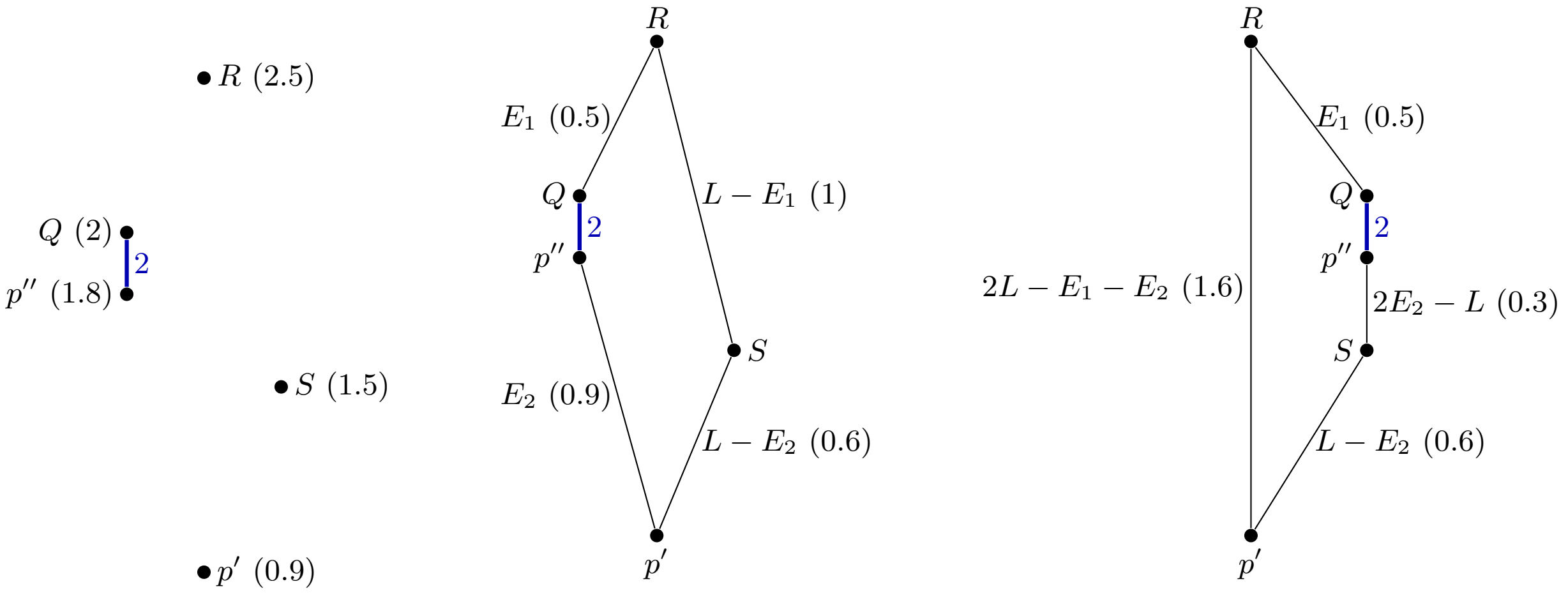

When $\fixedcomp$ is an isolated fixed point, $p=\fixedcomp$ is the only choice
and the statement reads
$\Ext(\widehat\alpha)=\{\T\in\Ext(\alpha):c<c^{\T}(M,\fixedcomp)\}$. When
$\fixedcomp$ is a fixed surface, its two ends are fixed points of every
$\T\in\Ext(\alpha)$, and a single torus may then contribute two distinct
extensions.

The equivariant capacity of $\fixedcomp$ is therefore at least the largest of the
capacities $c^{\T}(M,p)$, taken over the pairs $(\T,p)$ above. That it is equal
to it follows from the correspondence between toric extensions and strata, and
is proved in Corollary~\ref{cor:capacity is max}. 
\\

The following is an immediate consequence of invariance of capacities under equivariant symplectomorphisms. 

\begin{lemma}[Invariance]\label{lemma:invariance}
Let $\varphi:(M,\oml)\to(M',\om_{\lambda'})$ be a symplectomorphism between
rational ruled surfaces conjugating a Hamiltonian circle action $\alpha$ to a
circle action $\alpha'$ and mapping a fixed component $\fixedcomp$ onto a fixed
component $\fixedcomp'$. Then $\lambda=\lambda'$, conjugation by $\varphi$ is a
bijection from $\Ext(\alpha)$ onto $\Ext(\alpha')$, it preserves the capacities,
\[
c^{\T}(M,p)=c^{\varphi\T\varphi^{-1}}(M',\varphi(p))
\qquad\text{and}\qquad
c^{S^1}(M,\fixedcomp)=c^{S^1}(M',\fixedcomp'),
\]
and for every admissible $c$ the two blow-ups
$\widehat\alpha_{\fixedcomp,c}$ and $\widehat\alpha'_{\fixedcomp',c}$ are
equivariantly symplectomorphic, so that $\varphi$ also identifies
$\Ext(\widehat\alpha_{\fixedcomp,c})$ with
$\Ext(\widehat\alpha'_{\fixedcomp',c})$.
\end{lemma}

\begin{Remark}\label{rem:presentation independent}
The pair $(\alpha,\fixedcomp)$ may be written in several ways: by
Proposition~\ref{prop:two tori} and Lemma~\ref{lemma:second torus}, an action
extending to two tori is both an $S^1(1,b;r)$ and an $S^1(1,b';n)$, with
$n=|2b-r|$, and the vertex carrying the ball changes label accordingly.
Lemma~\ref{lemma:invariance} says that nothing computed from the pair depends on
that choice: the capacities, the number of toric extensions before and after the
blow-up, and hence all the homotopy types deduced from them, are invariants of
the equivalence class of $(\alpha,\fixedcomp)$ and of $c$.
\end{Remark}

\begin{Remark}
The three capacities of Lemma~\ref{lemma:capacity basics} are genuinely
different. For $r=4$ on $(\SSS,\om_{2.1})$ and $b=1$, the negative section of
$\T_4$ has area $0.1$ while the vertex $Q$ corresponds, under
Lemma~\ref{lemma:dictionary}, to a vertex of the positive section of $\T_2$. 
Hence $c^{\T_4}(M,Q)=0.1$ while $c^{S^1}(M,Q)=c^{\T_2}(M,Q)=1$, by
Corollary~\ref{cor:capacity is max}. An equivariant
ball may therefore be much larger than some invariant spheres through its
centre, and the toric action whose boundary are those invariant spheres need not survive the blow-up.
\end{Remark}

%%%%%%%%%%%%%%%%%%%%%%%%%%%%%%%%%%%%%%%%%%%%%%%%%%%%%%%%%%%%%%%%%%%%%%%%%%%%%%%%
\subsection{Stabilizers of invariant balls}
\label{subsec:stabilizers}
%%%%%%%%%%%%%%%%%%%%%%%%%%%%%%%%%%%%%%%%%%%%%%%%%%%%%%%%%%%%%%%%%%%%%%%%%%%%%%%%

As discussed above, in order to understand the homotopy type of space of equivariant embeddings we need to understand the fibration 
\[\Stab^G(\iota) \rightarrow \Symp^G(M,\oml) \twoheadrightarrow \IEmb(M,c,\fixedcomp)\]
where $\Stab^G(\iota)$ denotes the space of equivariant symplectomorphisms that send a given ball embedding to itself. In this section we relate the homotopy type of the stabiler with that of the space of equivariant symplectomorphisms of the blow up $\bX_2$. 

\begin{lemma}\label{lemma:indecomposable exceptional class}
Let a compact group $G$ act symplectically on a closed symplectic $4$-manifold
$(M,\om)$ and let $Z$ be an invariant exceptional sphere. Write
$\mathcal{S}^G_{[Z]}$ for the space of invariant embedded symplectic spheres in
the class $[Z]$, and $\mathcal{J}^G_\om([Z])\subseteq\mathcal{J}^G_\om$ for the
subspace of those invariant compatible structures for which $[Z]$ has an
embedded holomorphic representative. Then the map sending $J$ to that
representative is a fibration
\[
\mathcal{J}^G_\om([Z])\ \longrightarrow\ \mathcal{S}^G_{[Z]}
\]
with contractible fibres, and $\Symp^G(M,\om)$ acts on $\mathcal{S}^G_{[Z]}$
with local sections, the stabilizer of $Z$ being $\Symp^G(M,Z,\om)$.
Consequently:
\begin{enumerate}
\item if $[Z]$ is equivariantly indecomposable, then the inclusion
$\Symp^G(M,Z,\om)\into\Symp^G(M,\om)$ is a homotopy equivalence;
\item if $\mathcal{J}^G_\om([Z])$ is homotopy equivalent to a single orbit of $\Symp^G_h(M,\om)$,
presented by the orbit map as $\Symp^G_h(M,\om)/K$ with $K$ the stabilizer of
some $J_0\in\mathcal{J}^G_\om([Z])$, then the inclusion
\[
K\ \into\ \Symp^G_h(M,Z,\om):=\Symp^G(M,Z,\om)\cap\Symp^G_h(M,\om)
\]
is a homotopy equivalence.
\end{enumerate}
\end{lemma}
\begin{proof}
A class of square $-1$ has at most one embedded holomorphic representative, by
positivity of intersections, and that representative is invariant when $J$ is. 
The map of the statement is therefore well defined. It is onto, every invariant
embedded symplectic sphere being holomorphic for some invariant compatible
structure, and its fibre over $Z'$ is the space of invariant compatible
structures for which $Z'$ is holomorphic, contractible
by~\cite[Lemma A.2]{ChPin-Memoirs}. Being a fibration, it is a homotopy
equivalence. That the action on $\mathcal{S}^G_{[Z]}$ admits local sections is
the equivariant displacement theorem of~\cite[Lemma A.10]{ChPin-Memoirs}.

(1) Equivariant indecomposability means that
$\mathcal{J}^G_\om([Z])=\mathcal{J}^G_\om$, which is contractible, so that
$\mathcal{S}^G_{[Z]}$ is contractible as well. Since it is connected and the
orbits of $\Symp^G(M,\om)$ are open, the action is transitive and
\[
\Symp^G(M,Z,\om)\to\Symp^G(M,\om)\to\mathcal{S}^G_{[Z]}
\]
is a fibration with contractible base.

(2) The map of the statement is equivariant, so $\Symp^G_h(M,\om)$ acts
transitively on $\mathcal{S}^G_{[Z]}$ as well, with stabilizer
$\Symp^G_h(M,Z,\om)$ at the representative $Z_0$ of $J_0$. The two orbit fibrations therefore form a ladder
\[
\begin{tikzcd}
K \arrow[r] \arrow[d,hookrightarrow] & \Symp^G_h(M,\om) \arrow[r] \arrow[d,equal] & \mathcal{J}^G_\om([Z]) \arrow[d] \\
\Symp^G_h(M,Z,\om) \arrow[r] & \Symp^G_h(M,\om) \arrow[r] & \mathcal{S}^G_{[Z]}
\end{tikzcd}
\]
whose right-hand vertical map is the homotopy equivalence above, and so is
the left-hand one.
\end{proof}

\begin{thm}\label{thm:stabilizer}
Let $G$ be $S^1$ or $\T$ acting on a rational ruled surface $(M,\oml)$, let
$\iota\in\Emb^G(M,\fixedcomp,c)$ and let $(\bX_2,\tom)$ be the corresponding
equivariant blow-up, with exceptional divisor $Z$. Then
\[
\Stab^G(\iota)\ \simeq\ \Symp^G(\bX_2,\tom)
\]
in each of the following two cases:
\begin{enumerate}
\item $M=\SSS$ and $[Z]=E$;
\item $M=\CCC$, $c\leq\lambda-1$, and $[Z]=F-E$;
\end{enumerate}
\end{thm}
\begin{proof}
By Proposition~\ref{prop:homotopy equivalence Stab(B) and Stab(E)}, $\Stab^G(\iota)$ is homotopy equivalent to
$\Symp^G(\bX_2,Z,\tom)$, so it suffices to see that $[Z]$ is equivariantly
indecomposable and to apply part~(1) of
Lemma~\ref{lemma:indecomposable exceptional class}. In case (1) the class is
$E$ and in case (2) it is $F-E$, by Lemma~\ref{lemma:two blow-downs}; both are
indecomposable by Lemma~\ref{lemma:indecomposable}.
\end{proof}

\section{Centralizers of circle actions on $\bX_2$}\label{sec:Centralizers on bX2}

The results of the previous section reduce the determination of the stabilizer of an invariant embedded ball in $\SSS$ to the computation of the centralizer of the blown-up action. The same holds for balls in $\CCC$ if the capacity satisfies $c\leq \lambda-1$. When $M=\CCC$ and $c>\lambda-1$, the stabilizer is computed in~\S\ref{subsubsec:decomposable case} from the open stratum alone.

We compute the homotopy type of the centralizers following the framework of~\cite{ChPin-Memoirs}. From Corollary~\ref{cor:blow down}, it is sufficient to consider actions obtained from equivariantly blowing-up actions on $\SSS$ at a ball of capacity $c$. As before, we denote the fibre class by $F$, the class of a section of self-intersection $0$ by $B$, and the class of the exceptional divisor by $E=[\exdivisor]$. 

\subsection{Stratification of $\jomlc$ and toric extensions}\label{section:Strataofinvariantacs}

By Lemma~\ref{lemma:indecomposable}, the exceptional classes $E$ and $F-E$ are indecomposable. For $\lambda>1$, the class $B-E$ can be represented by reduced $J$-holomorphic curves. The way $B-E$ decomposes defines a stratification of the space $\jomlc$ of almost-complex structure tamed by $\omlc$. Define the homology classes
\[
D_{2k-1}:=B-kF \quad\text{and}\quad D_{2k}:=B-kF-E.
\]
Write $\lambda = \ell+\nu$ where $\ell$ is the largest integer strictly smaller than $\lambda$, and where $\nu\in (0,1]$.

\begin{lemma}[{\cite[Lemma~2.5]{P-com}}]\label{Lemma:NegativeselfintersectioncurvesinSSSC}
The space $\jomlc$ of tamed almost complex structures on $(\bX_2,\,\omlc)$
for which the class $D_0 = B-E$ is represented by an embedded $J$-holomorphic sphere is open and dense in $\jomlc$. If for a given $J$ there is no such sphere, then there is a unique integer $1\leq k \leq \ell$ such that either the class $D_{2k-1} = B-kF$
or the class $D_{2k} = B-kF-E$ is represented by a unique embedded $J$-holomorphic sphere. Moreover, the other $J$-holomorphic spheres with negative self-intersections represent the classes $E$ and $F-E$.\qed
\end{lemma}
Since the classes $D_i$ intersect negatively, it follows that the space $\jomlc$ of $\omlc$-compatible almost complex structures admits a finite partition
\begin{equation}\label{eq:Stratification}
\jomlc=\tilde{U}_0\sqcup \tilde{U}_1\sqcup\cdots\sqcup \tilde{U}_{N}
\end{equation}
where
\[
N=
\begin{cases}
2\ell   & \text{if~} c<\nu\\
2\ell-1 & \text{if~} c\geq\nu.
\end{cases}
\]
where each subspace $\tilde{U}_i$ is characterized by the existence of embedded $J$-holomorphic spheres representing the class $D_{i}$, that is, 
\[
\tilde{U}_{i}=\left\{J\in \jomlc~|~\text{there exists an embedded $J$-holomorphic sphere $C$ in class $D_i$}\right\}.
\]

\begin{prop}[{\cite[Proposition 2.8]{P-com}}]\label{prop:StructureStratification} Let $1\leq i\leq N$. Then
\begin{enumerate}
\item The subspace $\tilde{U}_i$ is a smooth, co-oriented, codimension $2i$ submanifold whose closure is the union $\bigsqcup_{i\leq m\leq N}\tilde{U}_{m}$.
\item The strata are in bijection with the equivalence classes of toric actions on the symplectic blow-up $(\bX_2,\,\omlc)$. 
\item Each stratum $\tilde{U}_{i}$ contains a complex structure $\tilde{J}_{i}$ coming from the blow-up of the Hirzebruch surface $W_{i+1}$ at the fixed point $P$ (away from the zero section $s_{0}$).
\item The number of strata, and hence the number of inequivalent $T^{2}$-actions on $(\bX_2,\,\omlc)$, changes exactly when $\lambda$ crosses an integer or when the capacity $c$ reaches the critical value $\nu$.
\end{enumerate}\qed
\end{prop}

%%%%%%%%%%

\begin{thm}\label{Theorem:ToricExtension=intersectionwithstrata}
The circle action $S^1(a,b;r;\fixedcomp)$ on $(\bX_2,\,\omlc)$ extends  to the toric action $\ttT_i$ on $\bX_2$ iff $\jsomc \cap \tilde{U}_i $ is non-empty.
\end{thm}
\begin{proof}
When the $S^1(a,b;r;\fixedcomp)$ extends to $\ttT_i$ then by definition the canonical $\ttT_i$ invariant integrable almost complex structure belongs to the intersection $\jsomc \cap \tilde{U}_i$. Conversely, suppose $ J \in \jsomc \cap \tilde{U}_i $, then there exists an invariant $J$-holomorphic curve $\overline{D}$ in class $D_i$. Furthermore, as $\overline{D}$ is invariant under the $S^1$ action, it has at least two fixed points $p_0$ and $p_1$ on it. Then we have either of the two following possibilities
\begin{enumerate}
    \item If $i=2k$ is even, then there exists two invariant $J$-holomorphic curves $F_{p_0}$ and $E_{p_1}$ in classes $F$ and $E$ respectively passing through $p_0$ and $p_1$. Furthermore, there is a curve $(F-E)_{p}$ in class $F-E$ that is disjoint from the curve in class $\overline{D}$.
    \item If $i=2k-1$ is odd, there exists two invariant $J$-holomorphic curves $F_{p_0}$ and $(F-E)_{p_1}$ in classes $F$ and $F-E$ respectively passing through $p_0$ and $p_1$. Furthermore, there is a curve $(E)_{p}$ in class $E$ that is disjoint from the curve in class $\overline{D}$.
\end{enumerate}
In either case, we have a curve $C$ of self-intersection $-1$ disjoint from our given curve $\overline{D}$. We can then $S^1$-equivariantly blow down along $C$ to obtain either the manifold $\SSS$ or the manifold $\CCC$, endowed with an $S^1$ action admitting an invariant curve of self-intersection $-(i+1)$. By~\cite[Proposition 3.22]{ChPin-Memoirs}, that action extends to the toric action $\T_{i+1}$, which lives on $\SSS$ when $i$ is odd and on $\CCC$ when $i$ is even. Blowing up again gives the extension of $S^1(a,b;r;\fixedcomp)$ to $\ttT_i$, by Lemma~\ref{lemma:dictionary}.
\end{proof}

\begin{figure}[ht]
\centering
\begin{subfigure}[b]{0.48\textwidth}
\centering
\begin{tikzpicture}[scale=0.85, every node/.style={scale=0.85}]
\coordinate (P) at (0,0);      \coordinate (S) at (1,1.95);
\coordinate (R) at (0,3.25);   \coordinate (Q) at (-1,2.6);
\draw (P) -- node[left]{$b$} (Q);
\draw [dashed, blue] (Q) -- (R);
\draw (R) -- node[right]{$b-1$} (S);
\draw [dashed, blue] (S) -- (P);
\foreach \v in {P,S,R,Q} \draw [fill] (\v) ellipse (0.08cm and 0.08cm);
\node[below] at (P) {$P,\ \mu=0$};
\node[right] at (S) {$S,\ \mu=\lambda$};
\node[above] at (R) {$R,\ \mu=\lambda-1+b$};
\node[left] at (Q) {$Q,\ \mu=b$};
\end{tikzpicture}
\caption{$\alpha$, for $b\geq2$}
\end{subfigure}
\hfill
\begin{subfigure}[b]{0.48\textwidth}
\centering
\begin{tikzpicture}[scale=0.85, every node/.style={scale=0.85}]
\coordinate (P) at (0,0);        \coordinate (S) at (1,1.95);
\coordinate (R) at (0,3.25);     \coordinate (q1) at (-0.9,3.75);
\coordinate (q2) at (-1.1,0.52);
\draw [dotted] (-1.6,3.25) -- (1.6,3.25);
\draw (P) -- node[below left=-2pt]{$b$} (q2);
\draw (q2) -- node[left]{$b+1$} (q1);
\draw [dashed, blue] (q1) -- (R);
\draw (R) -- node[right]{$b-1$} (S);
\draw [dashed, blue] (S) -- (P);
\foreach \v in {P,S,R,q1,q2} \draw [fill] (\v) ellipse (0.08cm and 0.08cm);
\node[below] at (P) {$P$};
\node[right] at (S) {$S$};
\node[below left] at (R) {$R$};
\node[above] at (q1) {$q_1,\ \mu=b+c,\ (1,-1-b)$};
\node[left] at (q2) {$q_2,\ \mu=b(1-c)$};
\end{tikzpicture}
\caption{blow-up at $Q$, for $b\geq2$}
\end{subfigure}

\bigskip
\begin{subfigure}[b]{0.48\textwidth}
\centering
\begin{tikzpicture}[scale=0.85, every node/.style={scale=0.85}]
\coordinate (Q) at (0,-2.6);   \coordinate (R) at (1,-1.95);
\coordinate (P) at (-1,0);     \coordinate (S) at (0,1.95);
\draw (Q) -- node[left]{$-b$} (P);
\draw [dashed, blue] (Q) -- (R);
\draw (R) -- node[right]{$1-b$} (S);
\draw [dashed, blue] (S) -- (P);
\foreach \v in {P,S,R,Q} \draw [fill] (\v) ellipse (0.08cm and 0.08cm);
\node[below] at (Q) {$Q,\ \mu=b$};
\node[right] at (R) {$R,\ \mu=\lambda-1+b$};
\node[left] at (P) {$P,\ \mu=0$};
\node[above] at (S) {$S,\ \mu=\lambda$};
\end{tikzpicture}
\caption{$\alpha$, for $b\leq-2$}
\end{subfigure}
\hfill
\begin{subfigure}[b]{0.48\textwidth}
\centering
\begin{tikzpicture}[scale=0.8, every node/.style={scale=0.8}]
\coordinate (S) at (0.4,1.95);  \coordinate (P) at (-1.2,0);
\coordinate (q2) at (-1.8,-0.8); \coordinate (q1) at (-0.6,-1.6);
\coordinate (R) at (1.2,-2.2);
\draw [dotted] (-2.2,-2.2) -- (2.2,-2.2);
\draw (P) -- node[right]{$-b$} (q2);
\draw (q2) -- node[above right=-2pt]{$-1-b$} (q1);
\draw [dashed, blue] (q1) -- (R);
\draw (R) -- node[right]{$1-b$} (S);
\draw [dashed, blue] (S) -- (P);
\foreach \v in {P,S,R,q1,q2} \draw [fill] (\v) ellipse (0.08cm and 0.08cm);
\node[left] at (P) {$P$};
\node[above] at (S) {$S$};
\node[below left=-1pt, align=right] at (R) {$R,\ \mu=\lambda-1+b$\\$(-1,1-b)$};
\node[below left=-1pt] at (q1) {$q_1,\ \mu=b+c$};
\node[left] at (q2) {$q_2,\ \mu=b(1-c)$};
\end{tikzpicture}
\caption{blow-up at $Q$, for $b\leq-2$}
\end{subfigure}
\caption{The graph of $\alpha=S^1(1,b;1)$ on $(\CCC,\oml)$ and the graph that its blow-up at $Q$ of size $c>\lambda-1$ would have, the height of a vertex being its moment value $\mu$. Solid edges are labelled by their isotropy weight and dashed blue edges are free spheres. The dotted line marks the height of $R$. In~(b) the highest vertex $q_1$ has weights of opposite signs, and in~(d) so does the lowest vertex $R$. For $b=1$, the vertices $R$ and $S$ form a fixed fibre at height $\lambda$, and for $b=-1$, the vertices $q_1$ and $q_2$ merge into the exceptional divisor, fixed at height $c-1$.}
\label{fig:corner graphs}
\end{figure}
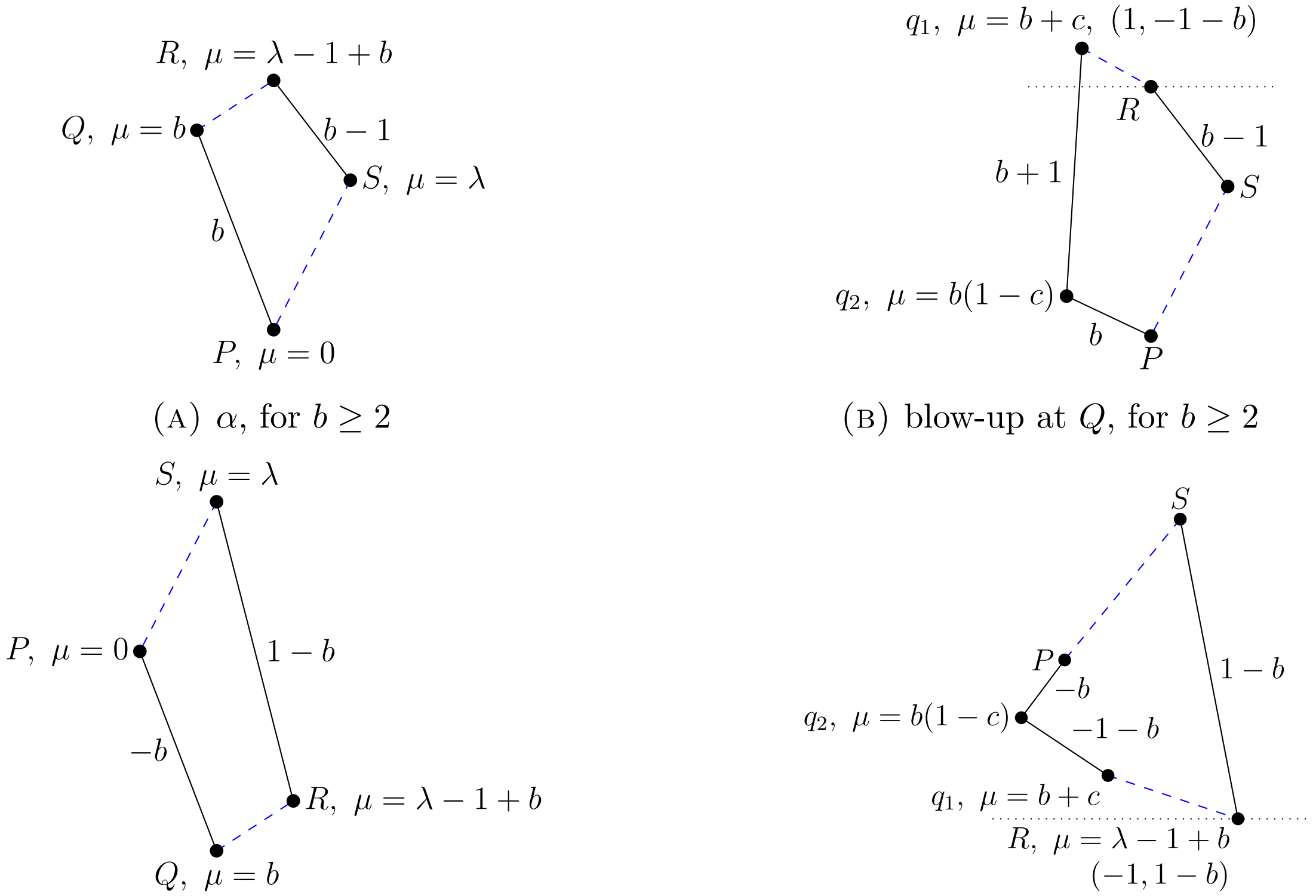

\begin{cor}\label{cor:capacity is max}
Let $\alpha$ be a Hamiltonian circle action on a rational ruled surface
$(M,\oml)$ and let $\fixedcomp$ be a component of its fixed point set. Then
\[
c^{S^1}(M,\fixedcomp)=c^{\T}_{\max}(\fixedcomp):=\max\ \bigl\{\,c^{\T}(M,p)\ :\ \T\in\Ext(\alpha),\
p\in\fixedcomp\ \text{a fixed point of }\T\,\bigr\}.
\] 
\end{cor}
\begin{proof}
A $\T$-equivariant ball is $S^1$-equivariant, which gives the inequality
$ c^{S^1}(M,\fixedcomp) \geq c^{\T}_{\max}(\fixedcomp)$. Conversely, let $c$ be the capacity of an $S^1$-equivariant ball centred on $\fixedcomp$. It is enough to find $\T\in\Ext(\alpha)$ and a fixed point $p\in\fixedcomp$ of $\T$ such that $c\leq c^{\T}(M,p)$.

Assume first that $M=\SSS$, or that $M=\CCC$ and $c\leq\lambda-1$. The
$S^1$-equivariant blow-up $\widehat\alpha_{\fixedcomp,c}$ exists and $\jsomc$ meets at least one of the strata, so that the set
$\Ext(\widehat\alpha_{\fixedcomp,c})$ is non-empty by Theorem~\ref{Theorem:ToricExtension=intersectionwithstrata}. Condition~(2) of Proposition~\ref{prop:capacity rule} holds, so that  $\Ext(\widehat\alpha_{\fixedcomp,c})$ consists of blow-ups $\T(p)$ with $\T\in\Ext(\alpha)$, $p\in\fixedcomp$. Consequently, $c^{S^1}<c^{\T}(M,p)$.

Assume now that $M=\CCC$ and $c>\lambda-1$. Since the Gromov width of $(\CCC,\oml)$ is $\oml(F)=1$, we have $1<\lambda<1+c<2$, so that $\T_1$ is the only toric action on $(\CCC,\oml)$ by Lemma~\ref{lemma:list of tori}.
Therefore, $\Ext(\alpha)=\{\T_1\}$, and we may assume that  $\alpha=S^1(a,b;1)$ with $a\geq0$ (and with $b=1$ when $a=0$). The moment polygon $\Delta_1$ has vertices $P=(0,0)$, $Q=(0,1)$, $R=(\lambda-1,1)$ and $S=(\lambda,0)$. The two edges meeting at $P$ or at $S$ have lengths $\lambda$ and $1$, and the two edges meeting at $Q$ or at $R$ have lengths $1$ and $\lambda-1$. Lemma~\ref{lemma:capacity basics} then gives
\[
c^{\T_1}(M,P)=c^{\T_1}(M,S)=1
\qquad\text{and}\qquad
c^{\T_1}(M,Q)=c^{\T_1}(M,R)=\lambda-1 .
\]
If $\fixedcomp$ contains $P$ or $S$, then $c^{\T}_{\max}(\fixedcomp)=1\geq c^{S^1}(M,\fixedcomp)$ and we are done. Otherwise $\fixedcomp$ is contained in the negative section $QR$, and we show that this contradicts the assumption $c>\lambda-1$. Let $x\in\fixedcomp$ be the centre of the ball and let $\Sigma$ be the invariant sphere over $QR$, whose class $B_{-1}$ has area $\lambda-1$.

Suppose first that $a=0$ or $a\geq2$. Then the stabilizer $\Gamma\subset S^1$ of the points of $\Sigma$ other than $Q$ and $R$ is, respectively, the whole circle or its subgroup $\Z_a$. Then $\Gamma$ fixes also the proper transform of $\Sigma$ in the blow-up, which lies in the class $B_{-1}-E$. By Lemma~\ref{LR}, it is $J$-holomorphic for every invariant almost complex structure $J$ tamed by the blown-up form, so that its area $\lambda-1-c$ is positive, contradicts the inequality $c>\lambda-1$.

Suppose finally that $a=1$. Then $QR$ is not fixed by the action, while $Q$ and $R$ are isolated fixed points exactly when $b\neq0$ and $b\neq1$ respectively, so that $\fixedcomp$ is either $Q$ with $b\neq0$ or $R$ with $b\neq1$. The affine symmetry $(x,y)\mapsto(\lambda-x-y,y)$ of $\Delta_1$ exchanges $P$ with $S$ and $Q$ with $R$, and replaces $S^1(a,b;1)$ by $S^1(a,a-b;1)$. Here it replaces $(1,b)$ by $(1,1-b)$, so that by Lemma~\ref{lemma:invariance} we may assume $\fixedcomp=Q$ and $b\neq0$. We argue on the graphs of $\alpha$ and of its blow-up in the sense of~\cite[\S2.1]{karshon}, drawn so that the height of a vertex is its moment value, see Figure~\ref{fig:corner graphs}. By the proof of Lemma~\ref{lemma:second torus}, the vertices $P,Q,R,S$ have heights $0$, $b$, $\lambda-1+b$ and $\lambda$, and the weights at $Q$ and $R$ are $(-b,1)$ and $(-1,1-b)$. In particular the free sphere $QR$ rises from $Q$ to $R$ by its area $\lambda-1$. The blow-up of size $c$ at $Q$ replaces the vertex $Q$ by two vertices and leaves the others unchanged~\cite[\S6.1]{karshon}. These are $q_1$, at height $b+c$ with weights $(1,-1-b)$, on the proper transform of $QR$, and $q_2$, at height $b(1-c)$ with weights $(-b,1+b)$, on the proper
transform of $QP$. When $b=-1$ they merge into a single fixed vertex at height $c-1$, the exceptional divisor. Since $c>\lambda-1$, the vertex $q_1$ lies above $R$. Now the highest vertex corresponds to the maximum of the moment map, at which both weights are non-positive, and the lowest vertex to the minimum, at which both weights are non-negative~\cite[Lemma~2.1 and \S2.1]{karshon}. If $b\geq1$, the vertex $q_1$ is the highest one, since the heights $0$, $\lambda$, $\lambda-1+b$ and $b(1-c)$ of the other vertices are smaller than $b+c$, but its weights have opposite signs. If $b\leq-1$, the vertex $R$ is the lowest one, since $\lambda-1+b$ is negative and smaller than $b+c\leq b(1-c)$, but its weights also have opposite signs. Both situations are impossible.
\end{proof}

The following  lemmata are useful to understand which strata the space of invariant almost complex structures intersect.

\begin{lemma}[Lemma 5.4 in~\cite{karshon}]\label{weight}
Let $S^1$  act on $S^2$ by rotating it $k$ times while fixing the
north and south poles. Suppose that the action lifts to a complex line bundle $E$ over $S^2$. Then $S^1$ acts linearly on the fibres over the north and south poles. Let $w(n)$ and $w(s)$ be the weights for these actions. Then
\[
w(n) - w(s) = -ek
\]
where $e$ is the self intersection of the zero section.\qed
\end{lemma}

\begin{lemma}\label{Lemma:SymplecticArea}
Suppose $S^1$ acts on $S^2$ by rotating it $k$ times while fixing the
north and south poles. Let $\mu(n)$ and $\mu(s)$ denote the values of the moment map at the north and south poles. Then the symplectic area of the sphere is $(\mu(n) - \mu(s))/k$.\qed
\end{lemma}

\begin{thm}\label{Theorem:IntersectionwithStrataanot1}
Given the circle action $S^1(a,b;r)$ on $W_r$, let $p\in\fixedcomp$ be a fixed point of the toric action $\T_r$. If $a \neq \pm 1$ then the action $S^1(a,b;r;\fixedcomp)$ on $\bX_2$ only extends to the toric action $\tilde{T}_r(p)$ obtained by equivariantly blowing-up $W_r$ at $p$.
\end{thm}
\begin{proof}
The circle action $S^1(a,b;r)$ restricted to the invariant curve $\overline{D}_r$ in class $D_r$ has stabilizer $\Z_a\subset S^1(a,b;r)$ of order $|a|\geq 2$. After the equivariant blow-up at $p$, the subgroup $\Z_a$ still acts trivially on the proper transform $\overline{D}'$ of $\overline{D}_r$. By Lemma~\ref{LR}, the invariant curve $\overline{D}'$ is $J$-holomorphic with respect to any invariant structure $\tilde{J}\in\jsomc$. This shows that the entire space $\jsomc$ belong to the stratum characterised by the class $[\overline{D}']$. Theorem~\ref{Theorem:ToricExtension=intersectionwithstrata} implies the result. 
\end{proof}

\begin{lemma}\label{Lemma:Atmost2tori}
An action of the form $S^1(\pm1,b;r;\fixedcomp)$ extends to at most two inequivalent tori. 
\end{lemma}
\begin{proof}
We prove this for the actions of the form $S^1(1,b;r;\fixedcomp)$ and $r$ even, the other cases are similar. Identify $\bX_2$ with $\bX_2$, and assume that there are at least 3 inequivalent tori $\ttT_{i_1},\ttT_{i_2},\ttT_{i_3}$ extending the action $S^1(\pm1,b;r;\fixedcomp)$. Associated to each of these toric actions are $S^1$ invariant curves $\overline{D}_{i_j}$ representing distinct classes $D_{i_j}$ (as defined in Section~\ref{section:Strataofinvariantacs}) and which are $J_{i_j}$-holomorphic for some $J_{i_j}\in\jsomc$. Since the class $E$ is indecomposable, it is represented by an invariant $J_{i_j}$-holomorphic curve. We can thus perform the $S^1$-equivariant blow down along invariant representatives of $E$ to obtain three $S^1$ invariant curves $\overline{D}'_{i_j}$ on $\SSS$ representing classes of the form $B-k_{i_j}F$ for $k_{i_j}\geq 0$, no longer necessarily distinct. But by~\cite[Proposition 3.22]{ChPin-Memoirs}, these classes characterize toric extensions of circle actions on $\SSS$, and by \cite[Proposition 2.21]{ChPin-Memoirs}, we know that the blow down circle action can admit at most 2 inequivalent toric extensions. On the other hand, looking at the Delzant polytopes of Figure~\ref{fig:toric blow-ups}, we see that in order to obtain 3 different classes $D_{i_j}$ by blowing up a single circle action on $\SSS$ with a fixed size $c$, one must perform the toric blow-up of at least 2 inequivalent tori. This shows that the blown down $S^1$ action extends to exactly two toric actions, which must be of the forms $S^1(1,b;m)$ and $S^1(1,b';m')$ for some parameters satisfying $b\neq0$, $b'\neq0$, $m=2k>0$, and $m'=2k'\geq 0$. Moreover, two of the blow-ups must be done on the same toric action, say $S^1(1,b;m)$, at two different fixed points, one on the negative curve in class $B-kF$, and one on a curve in class $B+kF$.
These two blow-ups of the same size $c$ must yield the original action $S^1(\pm 1, b;r;\fixedcomp)$. By~\cite[Proposition 2.21]{ChPin-Memoirs}, the possible graphs for the $S^1(1,b;m)$ action on $\SSS$ are the ones depicted in Figure~\ref{fig:PossibleGraphsToricExtensions}. In each case, blowing-up at the maximum and at the minimum give different moment map labels (up to translation). Similarly, blowing up at the two interior fixed points yields a new edge (representing the exceptional divisor) with different moment map labels. Finally, blowing up at an extremum and at an interior fixed point yields different sets of edges. This shows that it is not possible to obtain the same circle action by performing an equivariant blow-up of $S^1(1,b;m)$ of size $c$ at two different fixed points. 
\end{proof}

%%%%%%%%%% Figure of graphs of possible toric extensions
\begin{figure}
\centering
~
% a>0, b>0, am-b>0
\subcaptionbox*{}{} % In order to create a subfigure counter
\subcaptionbox{When $a=1,~b>0$, and $m-b>0$\label{fig:a=1|b>0|m-b>0}}
[.43\linewidth]
{\begin{tikzpicture}[scale=0.9, every node/.style={scale=0.9}]
\draw [fill] (0,0) ellipse (0.1cm and 0.1cm);
\draw [fill] (1,-1.5) ellipse (0.1cm and 0.1cm);
\draw [fill] (0,-4) ellipse (0.1cm and 0.1cm);
\draw [fill] (-1,-2.5) ellipse (0.1cm and 0.1cm);
\draw (0,0) -- (1,-1.5);
\draw [dashed, blue] (0,0) -- (-1,-2.5);
\draw [dashed, blue] (1,-1.5) -- (0,-4);
\draw (-1,-2.5) -- (0,-4);
\node[left] at (-0.6,-3.2) {$b$};
\node[right] at (0.7,-0.5) {$m-b$};
\node[above] at (0,0.2) {$\mu = (\lambda + k)$};
\node[right] at (-0.6,-2.5) {$\mu = b$};
\node[right] at (1,-1.8) {$\mu = b + (\lambda - k) $};
\node[below] at (0, -4.2) {$\mu = 0$};
\node[left] at (-1,-2.5) {Q};
\node[left] at (0.9,-1.4) {R};
\node[right] at (0,0) {S};
\node[right] at (0,-4) {P};
\end{tikzpicture}}
% a>0, b>0, am-b<0
\subcaptionbox{When $a=1,~b>0$ and $m-b<0$
}
[.52\linewidth]
{\begin{tikzpicture}[scale=0.9, every node/.style={scale=0.9}]
\draw [fill] (0,0) ellipse (0.1cm and 0.1cm); 
\draw [fill] (1,-2) ellipse (0.1cm and 0.1cm); 
\draw [fill] (0,-4) ellipse (0.1cm and 0.1cm); 
\draw [fill] (-1,-3) ellipse (0.1cm and 0.1cm); 
\draw [dashed, blue] (0,0) .. controls (-4,-1.5) and (-2,-3.5) .. (0,-4);
\draw [dashed, blue] (1,-2) -- (-1,-3);
\draw (0,0) -- (1,-2);
\draw (-1,-3) -- (0,-4);
\node[right] at (0,-4) {P};
\node[left] at (-1,-3) {Q};
\node[left] at (0.9,-2) {S};
\node[right] at (0,0) {R};
\node[above] at (0,0.2) {$\mu = b +  (\lambda -k)$};
\node[right] at (1.2,-2) {$\mu = (\lambda + k) $};
\node[right] at (-0.7,-3) {$\mu = b$};
\node[right] at (-0.4,-3.5) {$b$};
\node[right] at (0.7,-1) {$b-m$};
\node[below] at (0,-4.2) {$\mu = 0$};
\end{tikzpicture}
}
% a>0, b<0
\subcaptionbox{When $a=1,~b<0$\label{fig:a=1|b>0}}
[.5\linewidth]
{\begin{tikzpicture}[scale=0.9, every node/.style={scale=0.9}]
\draw [fill] (0,0) ellipse (0.1cm and 0.1cm);
\draw [fill] (1,-3) ellipse (0.1cm and 0.1cm);
\draw [fill] (0,-4) ellipse (0.1cm and 0.1cm);
\draw [fill] (-1,-2) ellipse (0.1cm and 0.1cm);
\draw [dashed, blue] (0,0) .. controls (-4,-1.5) and (-2,-3.5) .. (0,-4);
\draw [dashed, blue] (-1,-2) -- (1,-3);
\draw (0,0) -- (-1,-2) ;
\draw (1,-3) -- (0,-4);
\node[above] at (0,1) {\rule{0em}{1em}}; % added to make a larger vertical separation between the graphs
\node[above] at (0,0.2) {$\mu = (\lambda +k)-b$};
\node[right] at (1.4,-3) {$\mu = -b $};
\node[right] at (-0.8,-1.9) {$\mu = (\lambda-k)$};
\node[below] at (0,-4.2) {$\mu = 0$};
\node[right] at (0.5, -3.6) {$-b$};
\node[right] at (-0.4,-1.0) {$m-b$};
\node[right] at (1,-3) {P};
\node[left] at (-1,-2) {R};
\node[right] at (0,0) {S};
\node[right] at (0,-4) {Q};
\end{tikzpicture}}
\caption{Graphs of $S^{1}(1,b;m)$ with possible locations of invariant curves in classes $B\pm k_{i_k}F$ and with moment map labels $\mu$.} 
\label{fig:PossibleGraphsToricExtensions}
\end{figure}
%%%%%%%%%% end of figure

\begin{lemma}\label{Lemma:Fixedsurface_toricextension}
Consider an action of the form $S^1(\pm1,b;r;\fixedcomp)$ on $(\bX_2,\,\omlc)$ with at least one fixed surface. 
\begin{enumerate}
\item If the circle action has two fixed surfaces, then it only extends to a single toric action.
\item If the circle action has a single fixed surface whose homology class is neither $E$ nor $F-E$, then it only extends to a single toric action.
\item If the circle action has a single fixed surface whose homology class is either $E$ or $F-E$, then  we have the following 2 cases: 
\begin{enumerate}
\item When $2\lambda < |2b-r| $ then it extends to only one torus.
\item When $2\lambda > |2b-r| $ then it extends to two tori. 
\end{enumerate}
\end{enumerate}
\end{lemma}
\begin{proof}
\textit{Case 1:} The two fixed surfaces are extrema of the moment map and, therefore, cannot both intersect any given invariant sphere $C$ in class $E$ (a fortiori if one fixed surface represents $E$ itself). By blowing-down the action along $C$, we get an action on $\SSS$ that contains a fixed surface. In particular, it is different from the types of circle actions depicted in Figure~\ref{fig:PossibleGraphsToricExtensions}. The same argument as in the proof of Lemma~\ref{Lemma:Atmost2tori} establishes the first claim.\\

\textit{Case 2:} Since the fixed surface does not represent $E$, the  argument used in Case 1 applies exactly as before.\\

\textit{Case 3a: $2\lambda > |2b-r|$}  By Lemma~\ref{Lemma:Atmost2tori}, we know the circle action extends to at most 2 tori. We shall now explicitly produce the 2 extensions, thus proving the claim. We present the extensions for the case $S^1(1,1,r,\fixedcomp)$ with $\fixedcomp$ is the isolated fixed point $P$, $r=2k$ in figure \ref{fig:even hirz}. The other cases are similar and we leave the details to the reader. 

\begin{figure}[H] % Graphs of possible toric extensions
\centering
~
% a>0, b>0, am-b>0
\subcaptionbox*{}{} % In order to create a subfigure counter
\setcounter{subfigure}{0} % To set the subfigure counter to 6
\subcaptionbox{When $a=1$, $r>1>0$ and $\fixedcomp = P$\label{fig:a=1|b=1|m-b>0}}
[.43\linewidth]
{\begin{tikzpicture}[scale=0.9, every node/.style={scale=0.9}]
\draw [fill] (0,0) ellipse (0.1cm and 0.1cm);
\draw [fill] (1,-1) ellipse (0.1cm and 0.1cm);
\draw [fill] (0,-4) ellipse (0.3cm and 0.1cm);
\draw [fill] (-1,-3) ellipse (0.1cm and 0.1cm);
\draw [blue](0,0) -- (1,-1);
\draw [dashed] (0,0) -- (-1,-3);
\draw [dashed] (1,-1) -- (0,-4);
\draw [dashed] (-1,-3)--(0,-4);
\node[right] at (0.7,-0.2) {$r-1$};
\node[above] at (0,0.0) {$\mu = \lambda +k-c$};
\node[right] at (1.1, -1) {$\mu= 1+\lambda-k-c$};
\node[left] at (-1.1,-3) {$\mu = 1-c$};
\node[below] at (0, -4.2) {$\mu = 0$};
\end{tikzpicture}}
\vspace{1mm}
% a>0, b>0, am-b>0
\subcaptionbox*{}{} % In order to create a subfigure counter
\setcounter{subfigure}{1} % To set the subfigure counter to 6
\subcaptionbox{When $a=1$, $r>1>0$ and $\fixedcomp = P$\label{fig:a=1|b>0|m-b>0-B}}
[.43\linewidth]
{\begin{tikzpicture}[scale=0.9, every node/.style={scale=0.9}]
\draw [fill] (0,0) ellipse (0.1cm and 0.1cm);
\draw [fill] (1,-1) ellipse (0.1cm and 0.1cm);
\draw [fill] (0,-4) ellipse (0.3cm and 0.1cm);
\draw [fill] (-1,-3) ellipse (0.1cm and 0.1cm);
\draw (0,0) -- (1,-1);
\draw [dashed] (0,0) -- (0,-4);
\draw [dashed] (-1,-3) -- (1,-1);
\draw [dashed](-1,-3)--(0,-4);
\node[right] at (0.7,-0.2) {$r-1$};
\node[above] at (0,0.0) {$\mu = \lambda +k-c$};
\node[right] at (1.1, -1) {$\mu= 1+\lambda-k-c$};
\node[left] at (-1.1,-3) {$\mu = 1-c$};
\node[below] at (0, -4.2) {$\mu = 0$};
\end{tikzpicture}}
\end{figure}

By Corollary 2.22 and 2.23 in \cite{ChPin-Memoirs}, The circle action $S^1(1,1;r)$ extends to two tori. The above two extended graphs arise by considering the blow-up at the minima of these two toric extensions, thus completing the proof for Case:3a.  
\\

\textit{Case 3b: $2\lambda < |2b-r|$}. Suppose on the contrary, that the circle action extended to two tori in this case. Without loss of generality denote them by $\ttT_{i_1}$, $\ttT_{i_2}$. Furthermore, we see that the invariant curve $D_{|2b-r|}$ joining the 2 interior fixed points has negative area and hence the circle action cannot extend to $\ttT_{|2b-r|}$. Increasing $\lambda$ such that $2\lambda > |2b-r|$, by Case 3a, we see that the circle action now extends to three tori $\ttT_{i_1}$, $\ttT_{i_2}$ and $\ttT_{|2b-r|}$ thus contradicting Lemma~\ref{Lemma:Atmost2tori}. Hence when $2\lambda < |2b-r|$, the circle action can extend to at most 2 tori. 
\end{proof}

We now determine the toric extensions of a blown-up circle action. By
Theorem~\ref{Theorem:ToricExtension=intersectionwithstrata} this amounts to determining which strata the space $\jsomc$ meets, and by
Proposition~\ref{prop:capacity rule} it is a comparison of the capacity of the ball with the equivariant capacities of the toric extensions of the action on the base.

\begin{thm}\label{Thm:Intersectionwithstrata_a=1}
Let $\alpha=S^1(a,b;r)$ act on $(M,\oml)$, let $\fixedcomp$ be a component of its fixed point set and let $\widehat\alpha=S^1(a,b;r;\fixedcomp)$ be its equivariant blow-up at a ball of capacity $c$. Assume $M=\SSS$, or $M=\CCC$ with $c\leq\lambda-1$, so that Proposition~\ref{prop:capacity rule} applies as an equality. Then:
\begin{enumerate}
\item if $a\neq\pm1$, or if $a=1$ and $b\in\{0,r\}$ with $r=0$, then
$\widehat\alpha$ has a single toric extension;
\item if $a=1$ and $b\in\{0,r\}$ with $r\geq1$, so that $\fixedcomp$ is a fixed fibre, then $\widehat\alpha$ extends to $\ttT_{r-1}$, and also to $\ttT_{r}$ when $2\lambda>r+\epsilon_r+2c$;
\item if $a=1$ and $b\notin\{0,r\}$, so that $\fixedcomp$ is one of the four vertices $P,Q,R,S$ of $\Delta_r$, then $\widehat\alpha$ extends to the toric actions listed in Table~\ref{table:NumberOfStrata}, each subject to the condition displayed there.
\end{enumerate}
\end{thm}
\begin{proof}
By Proposition~\ref{prop:capacity rule}, the toric extensions of $\widehat\alpha$ are the blow-ups $\T(p)$ with $\T\in\Ext(\alpha)$, $p$ a fixed point of $\T$ lying in $\fixedcomp$ and $c<c^{\T}(M,p)$, and by Lemma~\ref{lemma:capacity basics} that capacity is the smaller of the areas of the two invariant spheres through $p$: a fibre of area $1>c$, and a section. By Lemma~\ref{lemma:dictionary}, the blow-up of $\T_m$ gives $\ttT_{m-1}$ at a vertex of its positive section and $\ttT_{m}$ at a vertex of its negative section, the latter being admissible exactly when $2\lambda>m+\epsilon_m+2c$, by~\eqref{eq:negative section}.

For $a\neq\pm1$, Proposition~\ref{prop:two tori} gives $\Ext(\alpha)=\{\T_r\}$, so that $\Ext\bigl(\widehat\alpha_{\fixedcomp,c}\bigr)$ is a singleton, which is (1). If $\fixedcomp$ is a fixed fibre, its two ends are a vertex of the positive section and a vertex of the negative section of $\T_r$, and both are fixed points of $\T_r$ lying in $\fixedcomp$. They contribute the two distinct toric actions $\ttT_{r-1}$ and $\ttT_r$, the first always admissible and the second exactly when $c$ is smaller than the negative section. This gives (2). For $r=0$ the two toric blow-ups are exchanged by the symmetry $y\mapsto1-y$ of the rectangle, which fixes the circle, so only one stratum is met.

For (3), $\Ext(\alpha)=\{\T_r,\T_n\}$ with $n=|2b-r|$ by Proposition~\ref{prop:two tori}, the second being present only when
$2\lambda>n+\epsilon_n$. It remains to say which vertex of $\Delta_n$
corresponds to $\fixedcomp$, and that is the content of Lemma~\ref{lemma:second torus}, where the correspondence depends only on the sign of $2b-r$. Combining it with the previous paragraph gives
Table~\ref{table:NumberOfStrata}.
\end{proof}

\begin{table}[H]
\begin{center}
\begin{tabular}{|c|c||c|c|}
\hline
$2b-r$ & $\fixedcomp$ & first extension & second extension \\
\hline\hline
\multirow{4}{*}{$<0$}
 & $P$ & $\ttT_{r-1}$, \ $2\lambda>r+\epsilon_r$
       & $\ttT_{n}$, \ $2\lambda>n+\epsilon_n+2c$ \\\cline{2-4}
 & $Q$ & $\ttT_{r}$, \ $2\lambda>r+\epsilon_r+2c$
       & $\ttT_{n-1}$, \ $2\lambda>n+\epsilon_n$ \\\cline{2-4}
 & $R$ & $\ttT_{r}$, \ $2\lambda>r+\epsilon_r+2c$
       & $\ttT_{n}$, \ $2\lambda>n+\epsilon_n+2c$ \\\cline{2-4}
 & $S$ & $\ttT_{r-1}$, \ $2\lambda>r+\epsilon_r$
       & $\ttT_{n-1}$, \ $2\lambda>n+\epsilon_n$ \\
\hline
\multirow{4}{*}{$>0$}
 & $P$ & $\ttT_{r-1}$, \ $2\lambda>r+\epsilon_r$
       & $\ttT_{n-1}$, \ $2\lambda>n+\epsilon_n$ \\\cline{2-4}
 & $Q$ & $\ttT_{r}$, \ $2\lambda>r+\epsilon_r+2c$
       & $\ttT_{n}$, \ $2\lambda>n+\epsilon_n+2c$ \\\cline{2-4}
 & $R$ & $\ttT_{r}$, \ $2\lambda>r+\epsilon_r+2c$
       & $\ttT_{n-1}$, \ $2\lambda>n+\epsilon_n$ \\\cline{2-4}
 & $S$ & $\ttT_{r-1}$, \ $2\lambda>r+\epsilon_r$
       & $\ttT_{n}$, \ $2\lambda>n+\epsilon_n+2c$ \\
\hline
\end{tabular}
\end{center}
\caption{Toric extensions of $S^1(1,b;r;\fixedcomp)$ on $(\bX_2,\omlc)$, for $b\notin\{0,r\}$ and $n=|2b-r|$, under the standing assumption of
Theorem~\ref{Thm:Intersectionwithstrata_a=1}, which excludes the decomposable regime $M=\CCC$, $c>\lambda-1$, where the count may differ (Remark~\ref{rem:corner counterexample}). Each extension is present exactly when the displayed condition holds; the first one always is, that condition being also the one for the ball to exist at $\fixedcomp$ inside $\T_r$. The type of a vertex depends on $b$: for $0<b<r$ the minimum is $P$, the maximum is $S$ and $Q,R$ are interior; for $b>r$ the maximum is $R$ and $Q,S$ are interior; for $b<0$ the minimum is $Q$ and $P,R$ are interior.}
\label{table:NumberOfStrata}
\end{table}

\begin{cor}\label{cor:Next}
Under the assumptions of Theorem~\ref{Thm:Intersectionwithstrata_a=1}, the number $N_{\mathrm{ext}}$ of toric extensions of $\widehat\alpha$ is $1$ or $2$. It equals $2$ exactly when $\Ext(\alpha)$ has two elements whose capacities at $\fixedcomp$ both exceed $c$, or when $\fixedcomp$ is a fixed fibre with $r\geq1$ and $2\lambda>r+\epsilon_r+2c$. When $M=\CCC$ and $c>\lambda-1$, the count may differ, see Remark~\ref{rem:corner counterexample}.\qed
\end{cor}

\subsection{Homotopical description of strata}
The results of this subsection are the transpositions to $\bX_2$ of Chapter 3 of~\cite{ChPin-Memoirs}, where the same chain of fibrations is established for $\SSS$ and $\CCC$. We indicate for each statement the one it transposes. The arguments carry over verbatim, all the constructions involved (the spaces of invariant curves and configurations, the gauge groups, the normal bundles) being equivariant by definition, and the transitivity statements resting on the equivariant displacement theorem of~\cite[Lemma A.10]{ChPin-Memoirs}, which holds for any closed symplectic $4$-manifold. There is one exception. When the curve $\overline{D}$ of negative self-intersection is fixed pointwise by the circle, it meets an invariant exceptional curve at a point that is not determined by the action, and this curve has to be taken into account, see part~(2) of Proposition~\ref{cor:homotopy type of Fix}.

Consider an action $S^1(a,b;r;\fixedcomp)$ on $\bX_2$ and assume that the space of invariant almost complex structures $\jsomc$ meets the stratum $U_{k}$. By Theorem~\ref{Theorem:ToricExtension=intersectionwithstrata}, the circle action then extends to the toric action $\ttT_k$.

As the action of $\Symp^{S^1}_h(\bX_2,\omlc)$ on $\jsomc \cap U_k$ preserves integrability, it is not transitive. However, as explained in Section~\ref{section:Strataofinvariantacs}, for every $J\in \jsomc \cap U_k$, there exists a $J$-holomorphic invariant curve in class $D_k$. Moreover, as explained in~\cite[Lemma 3.41]{ChPin-Memoirs}), the space $\mathcal{S}^{S^1}_{D_{s}}$ of invariant curves in class $D_k$ is homotopy equivalent to $\jsomc \cap U_k$.
Furthermore, the argument used in the proof of~\cite[Lemma 3.37]{ChPin-Memoirs} applies mutatis mutandis to show that the action of 
$\Symp^{S^1}_h(\bX_2,\omlc)$ on $\mathcal{S}^{S^1}_{D_{k}}$ is  transitive with fibre $\Stab^{S^1}(\overline{D})$ consisting of equivariant symplectomorphisms that send a fixed curve $\overline{D}$ in class $D_k$ to itself. 
\begin{prop}[{compare~\cite[Corollary 3.38 and Lemma 3.39]{ChPin-Memoirs}}]\label{prop:Transitivityofcurves}
Assume that the space of invariant almost complex structures $\jsomc$ has a non-empty intersection with the stratum $U_{k}$. The action of $\Symp^{S^1}_h(\bX_2,\omlc)$ on $\mathcal{S}^{S^1}_{D_{k}}$ induces a Serre fibration
\[
\Stab^{S^1}(\overline{D}) \longrightarrow \Symp^{S^1}_{h}(\bX_2,\omlc) \longtwoheadrightarrow {\mathcal{S}^{S^1}_{D_{k}}} \mathbin{\xrightarrow{\text{~~~$\simeq$~~~}}} \jsomc \cap U_{k}\rule{0em}{2em}
\]\qed
\end{prop}

The stabilizer $\Stab^{S^1}(\overline{D})$ acts on $\overline{D}$ inducing a restriction map $\Stab^{S^1}(\overline{D}) \to  \Symp^{S^1}(\overline{D})$.
\begin{prop}[{compare~\cite[Lemma 3.42]{ChPin-Memoirs}}]\label{prop:first restriction is a fibration}
The restriction map
\[
\Fix^{S^1}(\overline{D}) \to \Stab^{S^1}(\overline{D}) \longtwoheadrightarrow \Symp^{S^1}(\overline{D}) \mathbin{\xrightarrow{\text{~~~$\simeq$~~~}}} S^1 ~\text{or}~ \SO(3)
\]
is a fibration with fibre $\Fix^{S^1}(\overline{D})$ consisting of all equivariant symplectomorphisms that act as the identity on the curve $\overline{D}$. The group $\Symp^{S^1}(\overline{D})$ is homotopy equivalent to $S^1$ whenever the circle action restricts to a non-trivial action on $\overline{D}$, while it is homotopy equivalent to $\SO(3)$ when $\overline{D}$ is pointwise fixed. 
\end{prop}
\begin{proof}
The proof of  $\Symp^{S^1}(\overline{D}) \simeq S^1 ~\text{or}~ \SO(3)$ is the content of~\cite[Lemma 3.43]{ChPin-Memoirs}. The transitivity of the action follows from \cite[Lemma 3.42]{ChPin-Memoirs}.
\end{proof}

\begin{prop}[{compare~\cite[Lemma 3.44]{ChPin-Memoirs}}]\label{prop:Fibration to gauge}
The restriction of the derivative to the symplectic normal bundle $N(\overline{D})$ of $\overline{D}$ induces a fibration
\[
\Fix^{S^1} (N(\overline{D})) \to \Fix^{S^1}(\overline{D}) \longtwoheadrightarrow  \Gauge1^{S^1}(N(\overline{D})) \mathbin{\xrightarrow{\text{~~~$\simeq$~~~}}} S^1\rule{0em}{2em}
\]
where $\Gauge1^{S^1}(N(\overline{D}))$ is the equivariant automorphisms group of $N(\overline{D})$ and where the fibre $\Fix^{S^1} (N(\overline{D}))$ is homotopy equivalent to the group $\Fix^{S^1}_{\near}(\overline{D})$ of equivariant symplectomorphisms that are identity in a neighbourhood of $\overline{D}$. 
\end{prop}
\begin{proof}
The proof follows mutatis mutandis the proof of~\cite[Lemma 3.44]{ChPin-Memoirs}.
\end{proof}

Assume, until part~(1) of Proposition~\ref{cor:homotopy type of Fix}, that the circle acts non-trivially on $\overline{D}$, which then contains exactly two fixed points. Given any invariant almost complex structure $J \in \jsomc \cap U_k$ with corresponding $J$-holomorphic curve $\overline{D}$ in class $D_k$, these two fixed points $p_0, p_1\in \overline{D}$ can be labelled so that,
\begin{enumerate}
\item if $k$ is even, then there exists two invariant $J$-holomorphic curves $F_{p_0}$ and $C_{p_1}$ in classes $F$ and $E$ respectively passing through $p_0$ and $p_1$,
\item if $k$ is odd, there exists two invariant $J$-holomorphic curves $F_{p_0}$ and $C_{p_1}$ in classes $F$ and $F-E$ respectively passing through $p_0$ and $p_1$.
\end{enumerate}
Let $\mathcal{S}^{S^1}_{\overline{D}\vee [F] \vee [C]}$ denote the  space of unparametrized $S^1$ invariant wedge of three spheres $S_1 \vee S_2 \vee S_3$ such that 
\begin{itemize}
\item $S_1 = \overline{D}$,
\item $S_2$ is a sphere in class $F$ which is equal to the invariant curve ${{F}_{p_0}}$ near $p_0$, 
\item $S_3$ is a sphere in class $[C_{p_1}]$ which is equal to the invariant curve $C_{p_1}$ near $p_1$.
\end{itemize}

\begin{prop}[{compare~\cite[Lemma 3.37]{ChPin-Memoirs}}]\label{prop:Fix near is contractible}
The group $\Fix^{S^1}_{\near}(\overline{D})$ acts transitively on the space of configurations $\mathcal{S}^{S^1}_{\overline{D}\vee [F] \vee [C]}$ and there is a fibration
\[
\Stab^{S^1}({F}_{p_0} \vee C_{p_1}) \cap \Fix^{S^1}_{\near}(\overline{D}) \to \Fix^{S^1}_{\near}(\overline{D}) \longtwoheadrightarrow   \mathcal{S}^{S^1}_{\overline{D}\vee [F] \vee [C]}\mathbin{\xrightarrow{\text{~~~$\simeq$~~~}}} \mathcal{J}_{\omlc}^{S^1}(\overline{D})\simeq \{*\}.
\]
Moreover, the space $\mathcal{S}^{S^1}_{\overline{D}\vee [F] \vee [C]}$ is homotopy equivalent to the contractible space $\mathcal{J}_{\omlc}^{S^1}(\overline{D})$ of compatible invariant almost complex structures on $\bX_2$ for which the curve $\overline{D}$ is holomorphic. Finally, the fibre is the space $\Stab^{S^1}({F}_{p_0}\vee C_{p_1}) \cap \Fix^{S^1}_{\near}(\overline{D})$ of equivariant symplectomorphisms that that are identity in a neighbourhood of $\overline{D}$ and send the fixed curves ${F}_{p_0} \sqcup C_{p_1}$ to themselves.
\end{prop}
\begin{proof}
The proof is similar to that of~\cite[Lemma 3.45]{ChPin-Memoirs}.
\end{proof} 

\begin{prop}[{compare~\cite[Lemma 3.47]{ChPin-Memoirs}}]\label{prop:restriction two is an equivalence}
The group $\Stab^{S^1}({F}_{p_0} \vee C_{p_1}) \cap \Fix^{S^1}_{\near}(\overline{D})$ acts on the wedge $F_{p_0}\vee C_{p_1}$ and the restriction map defines a fibration 
\[ 
H  \to \Stab^{S^1}({F}_{p_0} \vee C_{p_1}) \cap \Fix^{S^1}_{\near}(\overline{D}) \longtwoheadrightarrow  \Symp^{S^1}(F_{p_0} \vee C_{p_1}, N(p_0,p_1)) \mathbin{\xrightarrow{\text{~~~$\simeq$~~~}}} \{*\}
\]
where $\Symp^{S^1}(F_{p_0} \vee C_{p_1}, N(p_0,p_1))$ is the group of equivariant symplectomorphisms of $F_{p_0}$ and $C_{p_1}$ which are identity in a neighbourhood of $p_0$
and $p_1$, and where the fibre
\[
H:= \Fix^{S^1}(F_{p_0} \vee C_{p_1})\cap \Fix^{S^1}_{\near}(\overline{D})
\]
consists of equivariant symplectomorphisms that pointwise fix a neighbourhood of $\overline{D}$ and fix the curve $F_{p_0}$ and $E_{p_1}$. Furthermore, $\Symp^{S^1}(F_{p_0} \vee C_{p_1}, N(p_0,p_1))$ is contractible. 
\end{prop}
\begin{proof}
The proof of the transitivity of action is similar to that of~\cite[Lemma 3.46]{ChPin-Memoirs}.     
\end{proof}

\begin{prop}[{compare~\cite[Lemma 3.47]{ChPin-Memoirs}}]\label{prop:H is contractible}
The restriction map
\[
\begin{tikzcd}[row sep=tiny, column sep=small]
\{*\} \simeq H'\ar[r] & H \ar[r,two heads] & \Gauge1^{S^1}(N(\overline{D}\vee F_{p_0} \vee C_{p_1})) \simeq \{*\} \\
 & \phi \ar[r,maps to] & \left.d\phi\right|_{N(\overline{D}\vee F_{p_0} \vee C_{p_1})}
\end{tikzcd}
\]
onto the group $\Gauge1^{S^1}(N(\overline{D}\vee F_{p_0} \vee C_{p_1}))$ of equivariant gauge transformations of the normal bundle of $\overline{D}\vee F_{p_0} \vee C_{p_1}$ which is identity in a neighbourhood of $\overline{D}$ is a fibration. Its fibre $H'$ is homotopy equivalent to the group
\[
\Symp^{S^1}(\id\near \overline{D}\vee F_{p_0} \vee C_{p_1})
\]
of equivariant symplectomorphisms that are the identity in a neighbourhood of $\overline{D}\vee F_{p_0} \vee C_{p_1}$. Furthermore, both $\Gauge1^{S^1}(N(\overline{D}\vee F_{p_0} \vee C_{p_1}))$ and $\Symp^{S^1}(\id\near \overline{D}\vee F_{p_0} \vee C_{p_1})$ are contractible. It follows that $H$ is also contractible.
\end{prop}
\begin{proof}
The transitivity of the action and the fact that $\Gauge1^{S^1}(N(\overline{D}\vee F_{p_0} \vee C_{p_1}))$ is contractible follow from arguments similar to the ones appearing in the proof of \cite[Lemma 3.47]{ChPin-Memoirs}. The contractibility of $\Symp^{S^1}(\id\near \overline{D}\vee F_{p_0} \vee C_{p_1})$ follows directly from Proposition~\ref{prop:equivariant Gromov}. 
\end{proof}

When $\overline{D}$ is pointwise fixed, the configurations used above must be modified. Let $C$ be the class $E$ if $k$ is even and $F-E$ if $k$ is odd, so that $C\cdot D_k=1$. By Lemma~\ref{lemma:indecomposable}, the class $C$ is indecomposable. For every invariant compatible $J$ for which $\overline{D}$ is holomorphic, it is therefore represented by a unique embedded $J$-holomorphic sphere $C_J$, which is invariant and meets $\overline{D}$ transversally at a single point. When the circle acts non-trivially on $\overline{D}$ that point is the fixed point $p_1$, but when $\overline{D}$ is pointwise fixed it depends on $J$, so that the structures whose curve $C_J$ passes through a prescribed point no longer form a contractible space. We therefore first fix the curve in the class $C$.

\begin{prop}\label{cor:homotopy type of Fix}
\begin{enumerate}
\item If the circle acts non-trivially on $\overline{D}$, the restriction of the derivative to the symplectic normal bundle
\[
\Fix^{S^1}(\overline{D}) \longtwoheadrightarrow  \Gauge1^{S^1}(N(\overline{D})) \mathbin{\xrightarrow{\text{~~~$\simeq$~~~}}} S^1\rule{0em}{2em}
\]
is a homotopy equivalence.
\item Suppose that $\overline{D}=\overline{D}_k$ is pointwise fixed, let $\overline{C}$ be the $J_k$-holomorphic sphere in the class $C$, let $p_1$ be its intersection point with $\overline{D}$, and write $\Stab^{S^1}(\overline{D}\cup\overline{C})$ for the group of equivariant symplectomorphisms that send each of the curves $\overline{D}$ and $\overline{C}$ to itself. Then
\begin{enumerate}
\item the inclusion $\Stab^{S^1}(\overline{D}\cup\overline{C})\into\Stab^{S^1}(\overline{D})$ is a homotopy equivalence;
\item the restriction to $\overline{D}$ defines a fibration
\[
\Fix^{S^1}(\overline{D})\cap\Stab^{S^1}(\overline{C})\to\Stab^{S^1}(\overline{D}\cup\overline{C})\longtwoheadrightarrow\Symp(\overline{D},p_1)\mathbin{\xrightarrow{\text{~~~$\simeq$~~~}}} S^1,
\]
where $\Symp(\overline{D},p_1)$ is the group of symplectomorphisms of $\overline{D}$ fixing $p_1$;
\item the restriction of the derivative to the symplectic normal bundle induces a homotopy equivalence
\[
\Fix^{S^1}(\overline{D})\cap\Stab^{S^1}(\overline{C})\longtwoheadrightarrow\Gauge1^{S^1}(N(\overline{D}))\mathbin{\xrightarrow{\text{~~~$\simeq$~~~}}} S^1 .
\]
\end{enumerate}
\end{enumerate}
\end{prop}
\begin{proof}
(1) Working backward, Propositions~\ref{prop:H is contractible}, \ref{prop:restriction two is an equivalence}, and \ref{prop:Fix near is contractible} show that the group $\Fix^{S^1}_{\near}(\overline{D})$ is contractible. The result then follows from Proposition~\ref{prop:Fibration to gauge}.

(2)(a) Let $\mathcal{S}^{S^1}_{C}(\overline{D})$ be the space of invariant embedded symplectic spheres in the class $C$ that are $J$-holomorphic for some invariant compatible $J$ for which $\overline{D}$ is holomorphic. As in Lemma~\ref{lemma:indecomposable exceptional class}, the map $J\mapsto C_J$ from the contractible space $\mathcal{J}^{S^1}_{\omlc}(\overline{D})$ onto $\mathcal{S}^{S^1}_{C}(\overline{D})$ is a fibration, whose fibre over $C'$ is the contractible space of invariant compatible structures for which $\overline{D}$ and $C'$ are both holomorphic. Hence $\mathcal{S}^{S^1}_{C}(\overline{D})$ is contractible. The group $\Stab^{S^1}(\overline{D})$ acts on it with local sections, as in Proposition~\ref{prop:Transitivityofcurves}, hence transitively, and the stabilizer of $\overline{C}$ is $\Stab^{S^1}(\overline{D}\cup\overline{C})$. The orbit fibration therefore has a contractible base.

(2)(b) Every element of $\Stab^{S^1}(\overline{D}\cup\overline{C})$ fixes $p_1=\overline{D}\cap\overline{C}$. Since the circle acts trivially on $\overline{D}$, the restriction to $\overline{D}$ is subject to no other condition. That it is a fibration onto $\Symp(\overline{D},p_1)$ is proved as in Proposition~\ref{prop:first restriction is a fibration}, the extensions being chosen to preserve $\overline{C}$ near $p_1$. Finally $\Symp(\overline{D},p_1)$ retracts onto the circle of rotations of $\overline{D}\cong S^2$ about $p_1$.

(2)(c) As in Proposition~\ref{prop:Fibration to gauge}, this map is a fibration whose fibre is homotopy equivalent to $\Fix^{S^1}_{\near}(\overline{D})\cap\Stab^{S^1}(\overline{C})$. Choose a point $p_0\in\overline{D}$ different from $p_1$. For every invariant $J$ for which $\overline{D}$ and $\overline{C}$ are holomorphic, the fibre class $F$ is represented through $p_0$ by a unique embedded $J$-holomorphic sphere $F_{p_0}$, which is invariant. Indeed, the only singular representative of $F$ is the union of $\overline{C}$ with a curve in the class $F-C$ (see the proof of~\cite[Lemma 2.3]{P-com}), and the latter is disjoint from $\overline{D}$ since $(F-C)\cdot D_k=0$. The group $\Fix^{S^1}_{\near}(\overline{D})\cap\Stab^{S^1}(\overline{C})$ acts on the space of configurations $\overline{D}\vee S_2\vee\overline{C}$, with $S_2$ in the class $F$ and equal to $F_{p_0}$ near $p_0$. As in Proposition~\ref{prop:Fix near is contractible}, this action is transitive, the space of configurations is homotopy equivalent to the contractible space of invariant compatible structures for which $\overline{D}$ and $\overline{C}$ are holomorphic, and the stabilizer of a configuration is the group $\Stab^{S^1}(F_{p_0}\vee C_{p_1})\cap\Fix^{S^1}_{\near}(\overline{D})$ of that proposition, with $C_{p_1}=\overline{C}$. That group is contractible by Propositions~\ref{prop:restriction two is an equivalence} and~\ref{prop:H is contractible}. Hence $\Fix^{S^1}_{\near}(\overline{D})\cap\Stab^{S^1}(\overline{C})$ is contractible.
\end{proof}

Finally, combining the previous results, we are in a position to show that the invariant stratum $\jsomc \cap U_k$ is homotopy equivalent to a homogeneous space under the action of the equivariant symplectomorphism group. 

\begin{thm}[{compare~\cite[Theorem 3.48]{ChPin-Memoirs}}]\label{Thm:Homogenousstrata}
Consider a Hamiltonian $S^1(a,b;r;\fixedcomp)$ action on $(\bX_2,\,\omlc)$ and assume that $\jsomc$ meets the stratum $U_k$. Then $\jsomc \cap U_k$ is homotopy equivalent to a homogeneous space under the action of the equivariant symplectomorphism group. More precisely, the stabilizer of the complex structure $J_k \in \jsomc \cap U_k$ is the group $\Isom(\omlc, J_k)$ of equivariant K\"ahler isometries of the pair $(\omlc, J_k)$ acting trivially on homology, this group is the torus $\ttT_k$, and there is a homotopy equivalence
\[
\jsomc \cap U_k \simeq \frac{\Symp_h^{S^1}(\bX_2,\omlc)}{\ttT_k}.
\]
\end{thm}
\begin{proof}
We first identify $\Isom(\omlc,J_k)$, which contains $\ttT_k$. When $k\geq1$, the stabilizer of $J_k$ in the group of all symplectomorphisms of $(\bX_2,\omlc)$ acting trivially on homology is the torus $\ttT_k$, by~\cite[Proposition 2.8(iii)]{P-com}, so that $\Isom(\omlc,J_k)=\ttT_k$. When $k=0$, the elements of $\Isom(\omlc,J_0)$ are biholomorphisms acting trivially on homology. They preserve the $J_0$-holomorphic spheres in the exceptional classes $F-E$ and $B-E$, which are disjoint since $(F-E)\cdot(B-E)=0$. Blowing these spheres down identifies $(\bX_2,J_0)$ with the blow-up of $\CP^2$ at two points, which we may take to be $[0:1:0]$ and $[0:0:1]$, and identifies $\Isom(\omlc,J_0)$ with a compact subgroup of the group $G$ of projective transformations fixing these two points. The group $G$ consists of the classes of the invertible matrices
\[
\begin{pmatrix} u&0&0\\ x&v&0\\ y&0&w\end{pmatrix},
\]
and the homomorphism from $G$ to $(\C^*)^3/\C^*$ that records the diagonal has a kernel isomorphic to $\C^2$, which contains no non-trivial compact subgroup. Thus $\Isom(\omlc,J_0)$ maps injectively into the maximal compact subgroup of $(\C^*)^3/\C^*$, which is a $2$-torus. As the $2$-torus $\ttT_0$ maps onto it, $\Isom(\omlc,J_0)=\ttT_0$.

By Proposition~\ref{prop:Transitivityofcurves}, it remains to show that the inclusion $\ttT_k\into\Stab^{S^1}(\overline{D}_k)$ is a homotopy equivalence, where $\overline{D}_k$ is the $J_k$-holomorphic sphere in the class $D_k$. The restriction of $\ttT_k$ to $\overline{D}_k$ has for kernel the circle of $\ttT_k$ that fixes $\overline{D}_k$ pointwise, which acts on the normal bundle of $\overline{D}_k$ by rotations of weight $\pm1$, and for image a circle of rotations of $\overline{D}_k$. When $\overline{D}_k$ is not pointwise fixed by the circle action, this gives a ladder of fibrations
\[
\begin{tikzcd}
&\Fix^{S^1}(\overline{D}_k) \arrow[r] &\Stab^{S^1}(\overline{D}_k) \arrow[r,twoheadrightarrow] &\Symp^{S^1}(\overline{D}_k) \\
&S^1 \arrow[u,hookrightarrow] \arrow[r] &\ttT_k \arrow[u,hookrightarrow] \arrow[r] &S^1 \arrow[u,hookrightarrow]
\end{tikzcd}
\]
in which the first vertical arrow is an equivalence by part~(1) of Proposition~\ref{cor:homotopy type of Fix}, and the third one by Proposition~\ref{prop:first restriction is a fibration}. When $\overline{D}_k$ is pointwise fixed, the circle acting on $\bX_2$ is the kernel above, and the $J_k$-holomorphic sphere $\overline{C}$ in the class $C$ is $\ttT_k$-invariant. By part~(2)(a) of Proposition~\ref{cor:homotopy type of Fix}, it suffices to show that the inclusion $\ttT_k\into\Stab^{S^1}(\overline{D}_k\cup\overline{C})$ is a homotopy equivalence, and there is a ladder of fibrations
\[
\begin{tikzcd}
&\Fix^{S^1}(\overline{D}_k)\cap\Stab^{S^1}(\overline{C}) \arrow[r] &\Stab^{S^1}(\overline{D}_k\cup\overline{C}) \arrow[r,twoheadrightarrow] &\Symp(\overline{D}_k,p_1) \\
&S^1 \arrow[u,hookrightarrow] \arrow[r] &\ttT_k \arrow[u,hookrightarrow] \arrow[r] &S^1 \arrow[u,hookrightarrow]
\end{tikzcd}
\]
in which the first and third vertical arrows are equivalences by parts~(2)(c) and~(2)(b) of Proposition~\ref{cor:homotopy type of Fix}. In both cases, the inclusion of $\ttT_k$ is therefore a homotopy equivalence.
\end{proof}

\begin{Remark}
The curve $\overline{D}_k$ is pointwise fixed by the circle action iff the $S^1$ action on the blow up is of the form $S^1(0,\pm 1;k;\fixedcomp)$.
\end{Remark}

%%%%%%%%%%%%%%%%%%%%%%%%%%%%%%%%%%%%%%%%%%%%%%%%%%%%%%%%%%%%%%%%%%%%%%%%%%%%%%%%
\subsection{Calculation of codimension}\label{section:codimension}
%%%%%%%%%%%%%%%%%%%%%%%%%%%%%%%%%%%%%%%%%%%%%%%%%%%%%%%%%%%%%%%%%%%%%%%%%%%%%%%%

By Theorem~\ref{Thm:Intersectionwithstrata_a=1}, we know $\jsomc$ intersects with 2 strata for certain circle actions of the form $S^1(1,b;r,\fixedcomp)$. Furthermore, by Theorem \ref{Thm:Homogenousstrata}, each stratum is homotopy equivalent to a homogenous space under the action of the group of equivariant symplectomorphisms. In this section, we calculate the codimension of a given stratum in $\jsomc$ and we determine the isotropy representation at a standard complex structure. We follow the approach of~\cite[Sections 2 and 3]{AGK} and~\cite[Section 2]{Anj-Lal-Pin} which consists in restricting the problem to integrable compatible complex structures, and then in analysing their inclusion in the space of all compatible almost complex structures. We adapt the arguments to the equivariant setting as in~\cite[Chapter 5]{ChPin-Memoirs} where similar results are obtained for circle actions on $\SSS$ and $\CCC$.
Again, as the arguments apply mutatis mutandis to the one-point blow-up, we only explain the main steps and refer the reader to these references for further details.

Let $\bX_2=(\bX_2, \omlc)$ equipped with an action $S^1(1,b;r,\fixedcomp)$. Denote by $\iomlc\subset\jomlc$ the subset of integrable compatible almost complex structures, and by $\isomlc\subset\jsomc$ the invariant ones. Given an integrable $J$, we write $H^{0,q}_{J}(\bX_2)$ for the $q^{\text{th}}$ Dolbeault cohomology group with coefficients in the sheaf of germs of holomorphic functions, and $H^{0,q}_{J}(T\bX_2)$ for the $q^{\text{th}}$ Dolbeault cohomology group with coefficients in the sheaf of germs of holomorphic vector fields. We then write $H^{0,q}_{J}(\bX_2)^{S^1}$ and $H^{0,q}_{J}(T\bX_2)^{S^1}$ for the subspaces invariant under the induced $S^1$ action. As explained in~\cite{Cathelineau}, these latter groups describe infinitesimal equivariant deformations of $J$.
\begin{lemma}[{\cite[Lemma A.5]{Anj-Lal-Pin}}]
For any $J\in \iomlc$, the cohomology groups $H^{0,2}_{J}(\bX_2)$ and $H^{0,2}_{J}(T\bX_2)$ are zero.\qed
\end{lemma}

\begin{cor}[{\cite[Theorem 5.4]{ChPin-Memoirs}}]\label{cor:thmA}
The space $\iomlc^{S^1}$ is a submanifold of $\jomlc^{S^1}$. Moreover, the tangent space of $\iomlc^{S^1}$ at $J$ is naturally identified with $T_{J}((\Diff(\bX_2) \cdot J)\cap \iomlc^{S^1})\oplus H^{0,1}_{J}(T\bX_2)^{S^1}$. Here, $H^{0,1}_{J}(T\bX_2)^{S^1}$ represents the infinitesimal equivariant compatible deformations of $J$ in $\iomlc^{S^1}$ and coincides with the infinitesimal equivariant deformations of $J$ in the set of all invariant integrable structures $\Ii^{S^1}$ on $\bX_2$.\qed
\end{cor}
The next step is to show that the inclusion $\isomlc\into\jsomc$ is transverse to the inclusion of each non-empty stratum $U_i\cap\jsomc \into\jsomc$, which implies that the stratification of $\jsomc$ induces a stratification of $\isomlc$ with strata whose normal bundles can be described using equivariant deformation theory as in Corollary~\ref{cor:thmA}.

\begin{prop}\label{prop:thmC}
Suppose $J_i\in\isomlc$ is an integrable complex structure in the stratum $\isomlc\cap U_i$. Let $u:S^2\to \bX_2$ be a $J_i$-holomorphic curve in class $D_i$. Then,
\begin{enumerate}
\item the map $u^{*}:H^{0,1}_{J_i}(T\bX_2)^{S^1}\to H^{0,1}(u^{*}(T\bX_2))^{S^1}$ is an isomorphism,
\item the inclusion $\isomlc\into\jsomc$ is transversal at $J_i$ to the inclusion of the stratum\\ $\jsomc\cap U_i\into\jsomc$,
\item the infinitesimal complement to the stratum $\jsomc\cap U_i\subset \jsomc$ at $J_i$ can be identified with the space of infinitesimal equivariant deformations $H^{0,1}_{J_i}(T\bX_2)^{S^1}$.
\end{enumerate}
\end{prop}
\begin{proof}
The first statement follows from~\cite[Lemma A.8]{Anj-Lal-Pin}. The second and third statements is a consequence of the first following the argument of~\cite[Theorem 5.15]{ChPin-Memoirs}.
\end{proof}

\begin{cor}\label{cor:codim given by deformations}
The complex codimension of a stratum $\jsomc \cap U_i$ containing an integrable complex structure $J_i$ is equal to the dimension of the space of $S^1$ invariant infinitesimal deformations of $J_i$, that this, to the rank of the space $H^{0,1}_{J_i}(T\bX_2)^{S^1}$.\qed
\end{cor}

In order to compute the invariant subspace $H^{0,1}_{J_i}(T\bX_2)^{S^1}$, we now recall how the torus $\T_i$ associated to the complex structure $J_i\in\jsomc \cap U_i$ acts on infinitesimal deformations $H^{0,1}_{J_i}(T\bX_2)^{S^1}$.
As in Proposition~\ref{prop:StructureStratification}, we denote by $\ttT_i$ the action obtained from blowing-up the standard toric actions on the Hirzebruch surface $W_{i+1}$ at the fixed point $P$, and we write $x:=(1,0)$ and $y:=(0,1)$ for the corresponding subcircles in $\ttT_i$.

\begin{prop}\label{prop:character}
The K\"ahler isometry group of $J_i$ is the torus $\ttT_i$. In the basis $\{x,y\}$, the character of the representation of $\ttT_i$, $i\geq 1$, on $H^{0,1}_{J_i}(T\bX_2)$, written multiplicatively, is given by 
\begin{equation}\label{eq:character}
\chi_i(x,y)=
\frac{1}{y}\left(x + x^2 + \cdots + x^{i}\right)
\end{equation}
In particular, the character of the action of a subcircle of $\ttT_i$ with direction $(a,b)$ in that basis is
\[
\chi_{a,b;i}(e^{\sqrt{-1}t}) = \chi_i(e^{\sqrt{-1}at}, e^{\sqrt{-1}bt}) = \sum_{j=1}^i e^{\sqrt{-1}t(aj-b)}
\]
\end{prop}
\begin{proof}
The character of the representation can be computed using the Atiyah-Bott fixed point formula as in \cite[Theorem 4.2]{AGK} and~{\cite[Proposition B.1]{Anj-Lal-Pin}}.
\end{proof}
\begin{cor}\label{Cor:Codimensionofstratum}
Suppose the circle action extends to $\ttT_i$, let $\overline D_i$ be the invariant curve of self-intersection $-(i+1)$ of that toric action, and let $q$ be one of the two fixed points lying on $\overline D_i$. Write $m$ for the weight of the circle action on the tangent line to $\overline D_i$ at $q$, and $m'$ for its weight on the other tangent line at $q$. Then the complex codimension of $\jsomc\cap U_i$ in $\jsomc$ is
\[
\begin{cases}
1 & \text{if } m\neq0 \text{ and } -m'/m\in\{1,\ldots,i\},\\
0 & \text{otherwise.}
\end{cases}
\]
In particular every non-empty stratum has complex codimension $0$ or $1$, and if the action extends to two toric actions, one of the two strata is open and the other has complex codimension $1$.
\end{cor}
\begin{proof}
In the basis $\{x,y\}$ of Proposition~\ref{prop:character} the circle has a direction $(a,b)$, and the dimension of the invariant subspace of $H^{0,1}_{J_i}(T\bX_2)$ is
\[
\int_{S^1} \chi_{a,b;i}(g)\, dg
= \frac{1}{2\pi} \sum_{j=1}^{i}\int_{0}^{2\pi} e^{\sqrt{-1}t(aj-b)} \,dt
= \#\{\,j\in\{1,\ldots,i\}\ :\ aj=b\,\},
\]
which is $0$ or $1$, the integers $j$ being distinct. It remains to read the condition $b=aj$ off the weights at $q$. In that basis the two fixed points of $\overline D_i$ have weight vectors $\{(1,0),(0,-1)\}$ and $\{(-1,0),(i+1,-1)\}$, the first vector of each pair pointing along $\overline D_i$. Pairing with $(a,b)$ gives $(m,m')=(a,-b)$ at the first point and $(m,m')=(-a,\,a(i+1)-b)$ at the second. In both cases $m\neq0$ means $a\neq0$, and the condition $-m'/m\in\{1,\ldots,i\}$ reads $b=aj$ with $j\in\{1,\ldots,i\}$ at the first point and $b=a(i+1-j)$ at the second. Since the two conditions define the same set, the criterion does not depend on the chosen end, nor on the
presentation. When $m=0$ the circle acts trivially on $\overline D_i$, so $a=0$, and then $b=\pm1$ by effectiveness and no term contributes.

By Corollary~\ref{cor:codim given by deformations} this dimension is the complex codimension of the stratum at $J_i$, which is homogeneous under $\Symp^{S^1}_h(\bX_2,\omlc)$ by Theorem~\ref{Thm:Homogenousstrata}. Finally, $\jsomc$ being connected, if it decompose into two strata, one must be open while the other has complex codimension~$1$.
\end{proof}

\begin{Remarks}\label{remark:open stratum = less negative}
(1) Corollary~\ref{Cor:Codimensionofstratum} shows that an effective Hamiltonian circle action on $\bX_2$ extends to more than one inequivalent toric action only if it acts on the curve $\overline D_i$ with weight $m=\pm1$. Indeed, effectiveness gives $\gcd(m,m')=1$, while the correspondence of Theorem~\ref{Theorem:ToricExtension=intersectionwithstrata} forces at least one stratum to have positive codimension, since $\jsomc$ is connected. The criterion then gives $-m'=jm$ for some $j$, so that $m$ divides $m'$ and $m=\pm1$.

(2) When a circle action $S^1(1,b;i)\subset\T_i$ on $\bX_2$ extends to exactly two toric actions, the open stratum contains almost complex structures that admit curves with higher (negative) self intersection, while the closed stratum contains almost complex structures that admit curves with strictly lower (negative) self intersection.
\end{Remarks}

\subsubsection{Genericity of almost complex structures}\label{section:genericity}

The results collected in this subsection and in the next one are not used elsewhere in this section. They serve in Section~\ref{sec:deformation to isotopy}, where paths of invariant pairs are deformed. A reader interested only in the homotopy types of \S\ref{section:homotopy type SSSC} may pass to it directly.

\begin{defn}\label{Defn:Genericstrata}
An almost complex structure in the codimension 0 stratum of $\jsomc$ is said to be \emph{generic} (in the Fredholm sense).
\end{defn}
In \cite{karshon}, a different notion of genericity was defined for invariant almost complex structures. Recall that  an invariant sphere is \emph{free} if the stablizier at its points, other than the poles, is trivial. A fixed point is \emph{interior} w.r.t.\ a moment map if it is not a maximum nor a minimum of the map. Recall also that any $S^1$-invariant compatible pair $(\omega,J)$  defines a {\bf compatible} Riemannian metric $\langle \cdot, \cdot \rangle=\omega(\cdot,J\cdot)$. The gradient vector field of the moment map $\Phi$ is then $-J \xi_{S^1}$, where $\xi_{S^1}$ is the vector field that generates the $S^1$-action (See equation 3.2 in~\cite{karshon}). The vector fields $\xi_{S^1}$ and $J \xi_{S^1}$ generate a $\C^{\times}=(S^1)^\C$-action. The closure of a non-trivial $\C^{\times}$-orbit is a sphere, called a {\bf gradient sphere}. On a gradient sphere, $S^1$ acts by rotation with two fixed points at the north and south poles. All other points on the sphere have the same stabilizer.

\begin{defn}\label{Defn:genericKarshon}
An invariant compatible almost complex structure $J$ is said to be generic in the sense of Karshon if the metric $g(\cdot,\cdot):= \om(\cdot,J\cdot)$ admits no free gradient sphere between two interior fixed points.
\end{defn}

\begin{lemma}[{\cite[Corollary 3.8]{karshon}}]\label{lem:generic1} For a {\it fixed} invariant symplectic form $\omega$, the subset of generic almost complex structures in the sense of Karshon is open and dense in $\mathcal{J}^{S^1}(\omega)$.\qed
\end{lemma}

\begin{Remark}
Karshon states this, and the results of~\cite[\S3]{karshon} that precede it, in
terms of invariant metrics compatible with $\om$ rather than of invariant
$\om$-compatible almost complex structures. The two formulations are equivalent:
the assignment $J\mapsto\om(\cdot,J\cdot)$ is a homeomorphism from the space of
invariant $\om$-compatible almost complex structures onto the space of invariant
$\om$-compatible metrics, with inverse given by polar decomposition, and it
carries the gradient flow of the moment map for one to the gradient flow for the
other. We work throughout with almost complex structures, and quote Karshon's
statements in that form.
\end{Remark}

We now proceed to show that for any Hamiltonian circle action on
$\bX_2$, the two notions of genericity coincide. To this end, we first relate gradient spheres to $J$-holomorphic spheres.

\begin{lemma}\label{cl:jhol=gradient} 
Given any invariant compatible pair $(\omega,J)$, a non-fixed sphere is an embedded $J$-holomorphic sphere if and only if it is a smooth gradient sphere w.r.t.\ the metric $\omega(\cdot,J \cdot)$.
\end{lemma}
\begin{proof}
By definition, if $C$ is a smooth gradient sphere w.r.t.\ the metric $\omega(\cdot,J \cdot)$ then the tangent space $T_p C$ at $p \in C$ is spanned by the vectors $\xi_{S^1}$ and $J\xi_{S^1}$. Moreover, a smooth non-fixed sphere with this property is a gradient sphere  w.r.t.\  $\omega(\cdot,J \cdot)$. 
\end{proof}
We will henceforth use the terms gradient sphere and holomorphic sphere interchangeably.

\begin{prop}\label{prop:Fredholm genericity implies Karshon}
Given any Hamiltonian circle action on $\bX_2 = (\bX_2,\omlc)$, an invariant compatible almost complex structure $J$ belongs to the open stratum of $\jsomc$ iff it is generic in the sense of Karshon.
\end{prop}
\begin{proof}
Choose $J\in\jsomc$ and suppose there is a free gradient sphere $C$ between two interior fixed points. Then, by Lemma~\ref{weight}, its self-intersection is strictly negative. By Lemma~\ref{Lemma:NegativeselfintersectioncurvesinSSSC}, the curve $C$ represents one of the classes $F-E$, $E$, or $D_i$ for some index $i\geq 0$. Again, Lemma~\ref{weight} shows that $C\cdot C=-1$ is not possible as the normal weights at an interior fixed point cannot be zero. Hence, $C$ represents the class $D_i$ for some $i\geq 1$ and, by the correspondence between curves in classes $D_i$ and strata, $J$ belongs to the stratum $J_i\in \jsomc\cap U_i$. If this is the open stratum, then representatives of the class $D_i$ persist under small deformations of $J_i$ inside $\jsomc$. However, by~\cite[Lemma 3.6]{karshon}, generic perturbations of the metric destroys free gradient spheres connecting interior fixed points. Therefore, $\jsomc\cap U_i$ cannot be open in $\jsomc$ so that $J$ is not Fredholm generic.

To show the converse, first note that if $\jsomc$ consists of a single (open) stratum, then $J$ is automatically Fredholm generic. So, suppose $\jsomc$ is made of two strata that correspond to two possible toric extensions. Looking at the extended graphs of these extensions, exactly one contains an interior free gradient sphere in some class $D_j$, while the other contains a free gradient sphere in class $D_i$ which connects one extremum. By a case-by-case analysis, we always have $D_i^2 > D_j^2$.  By Remark~\ref{remark:open stratum = less negative}, this shows that Karshon genericity implies Fredholm genericity.
\end{proof} 

\begin{cor}\label{cor:the set of generic pairs is connected}
Given any Hamiltonian circle action on $\bX_2$, the subset of generic almost complex structures in the sense of Karshon is open, dense, and connected.\qed
\end{cor}

\subsubsection{Genericity of paths of invariant symplectic pairs}\label{section:genericity of pairs}

Given a symplectic form $\om$ and an effective Hamiltonian $S^1$-action on $\bX_2$, let $\Omega^{S^1}_{\om}$ denote the space of invariant symplectic forms on $\bX_2$ that are equivariantly deformation equivalent to $\om$. Since a symplectic circle action on a $4$-manifold is Hamiltonian iff it has fixed points~\cite{McDuff-MomentMap}, the given circle action is Hamiltonian with respect to any form $\om'\in\Omega^{S^1}_{\om}$ and admits a moment map w.r.t. $\om'$.
Let $\mathcal{J}^{S^1}$ denote the space of invariant almost complex structures. We endow $\Omega^{S^1}_{\om}$ and $\mathcal{J}^{S^1}$ with the $C^\infty$ topology. Define the set 
\begin{equation}\label{eq:P}
\mathcal{P}_{\om}:= \left\{ (\om',J) \in \Omega^{S^1}_{\om}\times \mathcal{J}^{S^1}~|~ \text{$J$ is compatible with $\om'$} \right\}.
\end{equation}

\begin{prop}\label{prop:genericity space of pairs}
The set of generic pairs is an open, dense, and connected subset of $\mathcal{P}_{\om}\subset\Omega^{S^1}_{\om}\times \mathcal{J}^{S^1}$.
\end{prop}
\begin{proof}
We first observe that since $\Omega_\om^{S^1}$ is connected, and since the projection map $p:\mathcal{P}_\om\to\Omega_\om^{S^1}$ has contractible fibres, $\mathcal{P}_{\om}$ is itself connected. Now, if the circle action extends to only one toric action for every invariant symplectic form $\om'\in\Omega_\om^{S^1}$, then every pair $(\om',J)\in\mathcal{P}_\om$ is generic, in which case the claim holds trivially. So, suppose the circle action admits two toric extensions for some symplectic form $\om'$. We claim that the set of non-generic pairs is then a co-oriented submanifold of $\mathcal{P}_\om$ of complex codimension $1$. To see this, let $D_i$ be the class characterizing the closed stratum in $\mathcal{J}_{\om'}^{S^1}$ (i.e the homology class of the free gradient between two interior fixed points). Consider the universal moduli space 
\[
\mathcal{M}(D_i, \mathcal{P}_\om) :=\left\{ (u,\om'',J)~|~u:S^2\to\bX_2,~[u]=D_i,~(\om'',J)\in\mathcal{P}_\om,~\delbar_Ju=0\right\}.
\]
This is an ILB-manifold, and the image of its projection $\pi:\mathcal{M}(D_i, \mathcal{P}_\om)\to \mathcal{P}_\om$ is the set of pairs admitting a curve in the class $D_i$, that is the set of non-generic pairs. For a fixed symplectic form, the image is known to be a submanifold of complex codimension $1$, see~\cite[Corollary 5.13]{ChPin-Memoirs}. The same analysis carries over in the case the symplectic form varies in $\Om_\om^{S^1}$ as any symplectic curve persist under small perturbations of symplectic forms.
\end{proof}

\begin{cor}\label{cor:perturb}
Given a path $(\om_t,J_t) \in \mathcal{P}_\om$ connecting generic elements $(\om_0,J_0)$ and $(\om_1,J_1)$, there is a perturbation $(\tilde{\om}_t, \tilde{J}_t)$  of the path in $\mathcal{P}_\om$ with end points fixed such that $(\tilde{\om}_t, \tilde{J}_t)$ is generic for all $t$. \qed
\end{cor}

We now establish that invariant $J$-holomorphic curves (i.e. gradient spheres) persist along deformations of generic pairs. We first note that the definition of an interior fixed point does not depend on the element of a path-connected component of  $\Omega^{S^1}_\om(M)$.

\begin{lemma}\label{lemma:interior}
Let $\om \in \Omega^{S^1}_{\om}(M)$ and let $p$ be an interior fixed point for the moment map $\Phi_{\om}$ associated to $\om$. Then $p$ is an interior fixed point with respect to the moment map $\Phi_{\om'}$ for any $\om'$ in the same path-connected component of  $\om$ in $\Omega^{S^1}_{\om}(M)$.
\end{lemma}
\begin{proof}
To show that a fixed point is an interior fixed point it suffices to show that one of the weights at the fixed point is positive and the other is negative.  Given a symplectic form in $\Omega^{S^1}_{\om}(M)$, and a fixed point $p$, we have a continuous map 
$\Wt_p: \Omega^{S^1}_{\om}(M) \rightarrow \left(\Z \times \Z\right)/\Z_2$
that associates to each $\om \in \Omega^{S^1}_{\om}(M)$ the weight associated to the fixed point $p$. By continuity and the fact that the target is a discrete space, we get that the map is constant on path-connected components. 
\end{proof}
We also recall the criterion given by Ahara and Hattori for when an $(\omega,J)$-gradient sphere is smooth.
\begin{lemma}\label{lemma:Smmothness of gradient spheres}
{\cite[Lemma 4.9]{AhHa}.} A gradient sphere is smooth at its poles except in the following situation: if the gradient sphere is free and the pole in question is an isolated minimum (or maximum) of the moment map $\Phi$ with both isotropy weights $>1$ (or $<-1$).\qed
\end{lemma} 

\begin{prop}\label{prop:persist}
Let $(\om_t,J_t)$ be a path of generic pairs in the set $\mathcal{P}_\om$. Let $A$ be a homology class represented by an embedded sphere $\overline{A}_0$ between fixed points $p$ and $q$  that is a gradient sphere for the generic metric $g_0(\cdot,\cdot)= \om_0(\cdot, J_0\cdot)$. Further, in case $\overline{A}_0$ is a free gradient sphere from the maximum to the minimum, we assume the weights at $p$ to be $(-1,-a)$ and at $q$ to be  $(1,b)$, with $a,b \geq 0$. Then,  for all $t$, there is an embedded representative $\overline{A}_t$ that is a gradient sphere for the generic metric $g_t(\cdot,\cdot)= \om_t(\cdot, J_t\cdot)$  in class $A$ between the connected components of the fixed point set containing $p$ and $q$. 
\end{prop}

\begin{proof} The proof breaks into two cases depending on whether the sphere $\overline{A}_0$ is free or not.\\

\noindent\textbf{Case 1.} \textit{$\overline{A}_0$ is a free gradient sphere:} Recall that the moment map $\Phi_t$ associated to the pair $(\om_t, J_t)$ is a Morse-Bott function. There are two subcases to consider.\\

\textbf{1.a.} \textit{When either $p$ or $q$ is an interior fixed point:}  WLOG assume that $p$ is an interior point and q is a minimum of the moment map $\Phi_0$ with respect to $\om_0$. We show that there is a  gradient sphere between the interior fixed point $p$ and the connected component of $q$ in $M^{S^1}$ for any generic metric.  As the index of an interior fixed point is $2$, there is a $2$-dimensional family of gradient lines flowing out of the the critical point $p$. By genericity of the metric, this family needs to flow into the minimum $q$. Hence for any generic metric $g_t$ there always exists a gradient sphere between an interior fixed point $p$ and the connected component of $q$. This sphere is smooth by Lemma~\ref{lemma:Smmothness of gradient spheres}. By continuity, the  homology class of this sphere is the same as that of $\overline{A}_0$ for all $t$.\\

\textbf{1.b.} \textit{When $p$ is a maximum of the moment map and $q$ is a minimum:}  As above we show that there is always a  gradient sphere between the maximum and minimum of the moment map for any generic metric  $g_t(\cdot,\cdot)=\omega_t(\cdot,J_t \cdot)$ with $\omega_t \in \Omega^{S^1}_{\om}(M)$. Note that the union of unstable manifolds for the moment map $\Phi_t$ from all interior fixed points is a finite union of dimension $2$ manifolds. Hence the complement is open and dense. As the dimension of the stable manifold at the minimum is $4$,  there is a gradient trajectory between the minimum and the maximum. This sphere is smooth by Lemma~\ref{lemma:Smmothness of gradient spheres} as we assume that in case $p$ is a maximum and $q$ is a minimum the weights are $(-1,-a)$ at $p$ and $(1,b)$ at $q$ with $a,b \geq 0$. By continuity, the homology class of this sphere is the same as that of $\overline{A}_0$.\\

\noindent\textbf{Case 2.} \textit{$\overline{A}_0$ is a $\Z_k$-sphere:}  Consider the space 
\[
\begin{aligned}
T:= \{t \in [0,1]~|&~ \text{$A$ is represented by a smooth gradient sphere for the generic metric $g_t$} \\
&\text{passing through the connected components of $M^{S^1}$ containing $p$ and $q$}\}.
\end{aligned}
\]
As $T$ is non-empty (since $0 \in T$), it is enough to show that the set $T$ is connected.\\  

To show that $T$ is open, we note the following. Let $C$ be a $\Z_k$-sphere for $(\om,J)$, then there is a small neighbourhood $U$ of $(\om,J) \in \mathcal{P}$ such that $C$ is a $\Z_k$-sphere for all $(\tilde{\om},\tilde{J}) \in U$. This follows from the following observation. By Lemma~\ref{LR}, the curve $C$ is $\tilde{J}$-holomorphic (and hence $\tilde{\omega}$ symplectic) for all $(\tilde{\omega},\tilde{J}) \in U$ and is a $\Z_k$-sphere.\\

To show that $T$ is closed, let $t_1,\cdots, t_n \rightarrow t_{\infty}$. Let $\overline{A}_1, \cdots, \overline{A}_n $ be the  $\Z_{k}$ spheres corresponding to $g_{t_i}$. By Lemma~\ref{LR}, a $\Z_k$-sphere being a component of the fixed point set of the subgroup $\Z_k\subset S^1$ and therefore holomorphic for every invariant almost complex structure, we can take $\overline{A}_i = \overline{A}_1$ for all $i$. Hence $\overline{A}_1$ is a $\Z_k$-sphere for $g_{t_\infty}$.
\end{proof}

\begin{Remark}\label{Remark:GenericforSSSandCCC} 
The analogues of Propositions~\ref{prop:Fredholm genericity implies Karshon}, \ref{prop:genericity space of pairs}, and~\ref{prop:persist}, as well as of Corollary~\ref{cor:perturb} hold for Hamiltonian circle actions on $\SSS$ and $\CCC$.
\end{Remark}

%%%%%%%%%%%%%%%%%%%%%%%%%%%%%%%%%%%%%%%%%%%%%%%%%%%%%%%%%%%%%%%%%%%%%%%%%%%%%%%%
\subsection{Homotopy type of equivariant symplectomorphisms of $\bX_2$}\label{section:homotopy type SSSC}
%%%%%%%%%%%%%%%%%%%%%%%%%%%%%%%%%%%%%%%%%%%%%%%%%%%%%%%%%%%%%%%%%%%%%%%%%%%%%%%%

In this section, we compute the homotopy type of equivariant symplectomorphisms of $(\bX_2,\omlc)$. We first treat the case of a single toric extension.

\begin{thm}\label{Thm:Singlestratumhomotopytype}
Suppose the blown-up action $S^1(a,b;r;\fixedcomp)$ on $(\bX_2,\omlc)$ has a
single toric extension $\ttT$. Then the inclusion
$\ttT\into\Symp^{S^1}_h(\bX_2,\omlc)$ is a homotopy equivalence.
\end{thm}
\begin{proof}
Since there is a single toric extension, $\jsomc$ meets a single stratum, by
Theorem~\ref{Theorem:ToricExtension=intersectionwithstrata}. That stratum is
therefore all of $\jsomc$, which is contractible. By
Theorem~\ref{Thm:Homogenousstrata} the stratum is homotopy equivalent to
$\Symp^{S^1}_h(\bX_2,\omlc)/\ttT$, which gives the result.
\end{proof}

We now present the homotopy type of equivariant symplectomorphisms when the circle action extends to two toric actions. To do so, we need to first calculate the homology of the group of equivariant symplectomorphisms. 
\subsubsection{Homological injectivity}
In this section, we show that the $\ttT_k \hookrightarrow \Symp^{S^1}(\bX_2,\omlc)$ induce injective maps in homology and use it to calculate the homotopy type of equivariant symplectomorphisms.  
\begin{lemma}\label{Lemma:injectiveinhomology}
    The inclusion of $\ttT_k \hookrightarrow \Symp^{S^1}(\bX_2,\omlc)$ induces a map that is injective in homology. 
\end{lemma}
\begin{proof}
By~\cite[Proposition 3.15]{P-com}, the inclusion of $\ttT_k$ in
$\Symp(\bX_2,\omlc)$ is injective in homology, over any coefficient field. It
factors through $\Symp^{S^1}(\bX_2,\omlc)$, so the first map of that
factorisation is injective as well.
\end{proof}

Recall that the contractible space $\jsomc$ of invariant compatible almost
complex structures decomposes as the disjoint union
\[
\jsomc=(\jsomc\cap U_i)\sqcup(\jsomc\cap U_j)=:U^{S^1}_i\sqcup U^{S^1}_j,
\qquad i<j,
\]
where, by \S\ref{section:codimension}, both pieces are submanifolds of $\jsomc$,
the first one open and the second of real codimension $2$
(Corollary~\ref{Cor:Codimensionofstratum}). In particular $U^{S^1}_i=\jsomc-
U^{S^1}_j$ is connected. The closed stratum $U^{S^1}_j$ is co-oriented, being an
almost complex submanifold, so that the Alexander--Eells duality in the
contractible manifold $\jsomc$, as explained in Appendix D
of~\cite{ChPin-Memoirs}, induces isomorphisms
\begin{equation}\label{eq:AlexanderEellsIsomorphismHomology}
\lambda_{*}:H_{p}(U^{S^1}_j;k)\to H_{p+1}(U^{S^1}_i;k).
\end{equation}

By Theorem~\ref{Thm:Homogenousstrata} we also have the two fibrations
\begin{equation}\label{eq:TheTwoMainFibrations}
\begin{aligned}
\ttT_i&\to\Symp^{S^1}_h(\bX_2,\omlc)\xrightarrow{\ p_i\ }
\Symp^{S^1}_h(\bX_2,\omlc)/\ttT_i\simeq U^{S^1}_i,\\
\ttT_j&\to\Symp^{S^1}_h(\bX_2,\omlc)\xrightarrow{\ p_j\ }
\Symp^{S^1}_h(\bX_2,\omlc)/\ttT_j\simeq U^{S^1}_j .
\end{aligned}
\end{equation}
The first one shows that $\Symp^{S^1}_h(\bX_2,\omlc)$ is connected, its fibre and
its base being connected, and the second one then shows that $U^{S^1}_j$ is
connected as well. By Lemma~\ref{Lemma:injectiveinhomology} the two inclusions
$\ttT_i\into\Symp^{S^1}_h(\bX_2,\omlc)$ and
$\ttT_j\into\Symp^{S^1}_h(\bX_2,\omlc)$ are injective in homology over any 
field,. The tori having finite-dimensional homology, the restrictions in
cohomology are therefore surjective, so that the Leray--Hirsch theorem applies to
both fibrations and gives isomorphisms of graded modules (not of algebras) 
\begin{equation}\label{eq:SplittingCohomology}
H^*(\Symp^{S^1}_h(\bX_2,\omlc);k)\cong H^*(U^{S^1}_i;k)\otimes
H^*(\ttT_i;k)\cong H^*(U^{S^1}_j;k)\otimes H^*(\ttT_j;k),
\end{equation}
and, dualising over $k$, the corresponding splittings in homology
\begin{equation}\label{eq:SplittingHomology}
H_*(\Symp^{S^1}_h(\bX_2,\omlc);k)\cong H_*(U^{S^1}_i;k)\otimes
H_*(\ttT_i;k)\cong H_*(U^{S^1}_j;k)\otimes H_*(\ttT_j;k).
\end{equation}

\begin{thm}\label{Thm:Homologyofcentralizer}
Let $S^1(a,b;r;\fixedcomp)$ be a circle action that extends to two toric actions,
see Table~\ref{table:NumberOfStrata}. Then, for any field $k$,
\[
H^p\bigl(\Symp^{S^1}_h(\bX_2,\omlc); k\bigr) \simeq \begin{cases}
k^4 & p \geq 2,\\
k^3 & p =1,\\
k & p=0.
\end{cases}
\]
In particular the topological group $\Symp^{S^1}_h(\bX_2,\omlc)$ is of finite
type. The same holds for $\Symp^{S^1}(\bX_2,\omlc)$, which by
Lemma~\ref{lemma:Symp=Symph} coincides with $\Symp^{S^1}_h(\bX_2,\omlc)$ up to
homotopy, except for the actions $S^1(1,1;0,\Sigma)$ with $\lambda=1$.
\end{thm}
\begin{proof}
The proof is identical to \cite[Theorem 4.3]{ChPin-Memoirs}.
\end{proof}

\begin{thm}\label{Thm:pushout centralizer SSS}
Let an $S^1$-action on $(\bX_2,\,\omlc)$ be such that $\jsomc$ meets two strata,
$\jsomc\cap U_i$ and $\jsomc\cap U_j$ with $i<j$, the first one open and dense
and the second of real codimension $2$. Let $\widetilde P$ be the homotopy
pushout, in the category of topological groups, of the two inclusions
$\ttT_i\hookleftarrow S^1\hookrightarrow\ttT_j$. Then the canonical map
\[
\widetilde\Upsilon:\widetilde P\longrightarrow\Symp^{S^1}_h(\bX_2,\omlc),
\]
induced by the inclusions of the two tori, is a homotopy equivalence. In
particular
\[
\Symp^{S^1}(\bX_2,\omlc) \simeq \Omega S^3 \times \mathbb{T}^3 .
\]
\end{thm}
\begin{proof}
That exactly one of the two strata is open, and that the other has real
codimension $2$, is Corollary~\ref{Cor:Codimensionofstratum}. That the closed one
is the deeper of the two is Remark~\ref{remark:open stratum = less negative}. The
argument is then the gluing of~\cite[\S\S4.2.2 and 4.2.3]{ChPin-Memoirs}, which
applies verbatim to complete the proof.
\end{proof}

\begin{thm}
The following table describes the homotopy type of the space of $S^1$
equivariant symplectomorphisms on $(\bX_2,\,\omlc)$ in terms of the number
$N_{\mathrm{ext}}$ of toric extensions of the action, which
Theorem~\ref{Thm:Intersectionwithstrata_a=1} and
Corollary~\ref{cor:Next} determine.
\end{thm}

\begin{table}[H]
\begin{center}
\begin{tabular}{|l|c|c|}
\hline
action & $N_{\mathrm{ext}}$ & $\Symp_h^{S^1}(\bX_2,\omlc)$\\
\hline\hline
single toric extension & $1$ & $\simeq\ttT$\\
\hline
$a=\pm1$, two toric extensions & $2$ & $\simeq\Omega S^3\times\mathbb{T}^3$\\
\hline
\end{tabular}
\end{center}
\caption{Homotopy type of $\Symp_h^{S^1}(\bX_2,\omlc)$ according to the number
of toric extensions of the circle action. The value of $N_{\mathrm{ext}}$ is
read off Table~\ref{table:NumberOfStrata}: it equals $2$ exactly when the
action on the base extends to two tori whose equivariant capacities at
$\fixedcomp$ both exceed $c$, or when $\fixedcomp$ is a fixed fibre with
$r\geq1$ and $2\lambda>r+\epsilon_r+2c$.}
\label{table:CentralizersCircleSSSChomotopytype}
\end{table}

\begin{proof}
    The proof is a consequence of Theorems~\ref{Thm:Singlestratumhomotopytype} and \ref{Thm:pushout centralizer SSS}.
\end{proof}

\begin{lemma}\label{lemma:Symp=Symph}
Let $M=\bX_2$. For every circle action $S^1(a,b;r;\fixedcomp)$ on $(M,\omlc)$,
\[
\Symp^{S^1}(M,\omlc)\simeq
\begin{cases}
\Symp_h^{S^1}(M,\omlc)\times\Z_2 & \text{for the actions } S^1(1,1;0;\fixedcomp)
\text{ with }\lambda=1,\\[2pt]
\Symp_h^{S^1}(M,\omlc) & \text{otherwise.}
\end{cases}
\]
\end{lemma}
\begin{proof}
We have the fibration
\[
\Symp^{S^1}_{h}(M,\omlc) \rightarrow \Symp^{S^1}(M,\omlc)
\twoheadrightarrow G,
\]
where $G\subseteq\Aut_{c_1,\om}(H^2(M,\R))$ is the image of $\Symp^{S^1}(M,\omlc)$
in the group of automorphisms of $H^2(M,\R)$ preserving the classes of $c_1(TM)$
and of $\omlc$.

We first determine that group. The only classes of square $-2$ in $H_2(M;\Z)$
are $\pm(B-F)$, so the isometries of the intersection form preserving $c_1$ form
the group of order two generated by the reflection in $B-F$, which exchanges $B$
and $F$ and fixes $E$, equivalently exchanges the two disjoint exceptional
classes $B-E$ and $F-E$. That reflection preserves the class of $\omlc$ if and
only if $\lambda=1$. Hence $\Aut_{c_1,\om}(H^2(M,\R))\cong\Z_2$ when $\lambda=1$,
and is trivial otherwise.

It remains to say when the generator is realised by an equivariant
symplectomorphism. Blowing down the exceptional divisor, it is realised exactly
when the exchange of the two factors of $\SSS$ can be taken equivariant. That
exchange conjugates $S^1(a,b;0)$ to $S^1(b,a;0)$, and the two are the same
subgroup only for $a=b=1$: for instance it conjugates $S^1(1,-1;0)$ to its own
inverse, so it is not equivariant there. For the diagonal action
$S^1(1,1;0;\fixedcomp)$ with $\lambda=1$ it is equivariant, and the blow-ups at
the two fixed points it exchanges are equivariantly symplectomorphic, so the
generator is realised whichever component $\fixedcomp$ carries the ball. Hence
$G=\Z_2$ in that case and $G$ is trivial in all others.

Finally, $G$ being discrete, $\Symp^{S^1}(M,\omlc)$ is the union of $|G|$
components, any two of which are homeomorphic by translation by an element
mapping to the corresponding class. The stated homotopy types follow.
\end{proof}

%%%%%%%%%%%%%%%%%%%%%%%%%%%%%%%%%%%%%%%%%%%%%%%%%%%%%%%%%%%%%%%%%%%%%%%%%%%%%%%%
\section{Homotopy type of equivariant symplectomorphisms of $\CP^2$}
\label{section:Centralizers on CP2}
%%%%%%%%%%%%%%%%%%%%%%%%%%%%%%%%%%%%%%%%%%%%%%%%%%%%%%%%%%%%%%%%%%%%%%%%%%%%%%%%

Let $(M,\om)=(\CP^2,\om_{FS})$, with the Fubini--Study form normalized so that a
line has area one. Every symplectomorphism of $\CP^2$ acts trivially on homology,
so that $\Symp^{S^1}_h(\CP^2,\om)=\Symp^{S^1}(\CP^2,\om)$, and we drop the
subscript throughout. Recall also that $\Symp(\CP^2,\om)\simeq\PU(3)$, and that
$\CP^2$ carries a unique Hamiltonian toric action up to conjugation,
\[
(e^{i\theta_1},e^{i\theta_2})\cdot[z_0:z_1:z_2]
=[e^{i\theta_1}z_0:e^{i\theta_2}z_1:z_2],
\]
which we denote by $\T$, with moment polygon the standard triangle. By
Theorem~\ref{thm:extoric} every Hamiltonian circle action on $(\CP^2,\om)$ is
conjugate to a subcircle $S^1(a,b)\subset\T$, that is to
$e^{i\theta}\mapsto(e^{ia\theta},e^{ib\theta})$. We normalise $a\geq0$, and
$b=1$ when $a=0$.

\begin{lemma}\label{lemma:CP2 degenerate list}
The circle $S^1(a,b)$ fixes a line pointwise if and only if
$(a,b)\in\{(1,0),(0,1),(1,1)\}$. In that case the fixed point set consists of
that line together with one isolated fixed point, at which the two weights of the
circle are equal. In all other cases the fixed point set consists of the three
vertices of the triangle, at which the moment map takes three distinct values.
\end{lemma}
\begin{proof}
The circle fixes pointwise the sphere over a facet exactly when $(a,b)$ is
parallel to the inward normal of that facet, and the three inward normals of the
standard triangle are $(1,0)$, $(0,1)$ and $(-1,-1)$. With the normalisation
$a\geq0$ this gives the three pairs listed. At the vertex opposite to such a
facet the two edges leaving it are the two other facets, whose primitive
directions pair with $(a,b)$ to the same integer, so the weights are equal. In
general the moment values at the three vertices are $0$, $a$ and $b$, which are
pairwise distinct unless $a=b$, $a=0$ or $b=0$, that is unless the action is one
of the three above.
\end{proof}

We call the three circle actions of Lemma~\ref{lemma:CP2 degenerate list} that fix a line pointwise \emph{degenerate}. Equivalently, a circle action on $\CP^2$ is degenerate when its two weights at some fixed point are equal.%

\begin{cor}\label{cor:CP2 fixes the vertices}
If $S^1(a,b)$ has isolated fixed points, every equivariant symplectomorphism of
$(\CP^2,\om)$ fixes each of them.
\end{cor}
\begin{proof}
An equivariant symplectomorphism preserves the moment map, hence permutes the
fixed points preserving their moment values, which are pairwise distinct by
Lemma~\ref{lemma:CP2 degenerate list}.
\end{proof}

The next three lemmas concern an arbitrary invariant line, whether or not it is
pointwise fixed, and are used in both cases below.

\begin{lemma}\label{lemma:CP2 restriction fibration}
The restriction map $\Stab^{S^1}(\overline L)\to\Symp^{S^1}(\overline L)$,
$\phi\mapsto\phi|_{\overline L}$, is a Serre fibration with fibre
$\Fix^{S^1}(\overline L)$, the group of equivariant symplectomorphisms fixing
$\overline L$ pointwise.
\end{lemma}
\begin{proof}
As in~\cite[Lemma 3.42]{ChPin-Memoirs}.
\end{proof}

\begin{lemma}[{\cite[Lemma 3.43]{ChPin-Memoirs}}]\label{lemma:CP2 Symp of the line}
For a Hamiltonian circle action on $\overline L\cong\CP^1$,
\[
\Symp^{S^1}(\overline L)\simeq
\begin{cases}
S^1 & \text{if $\overline L$ is not pointwise fixed,}\\
\PU(2)\cong\SO(3) & \text{if it is.}
\end{cases}
\]\qed
\end{lemma}

\begin{lemma}\label{lemma:CP2 gauge}
The group $\Fix^{S^1}(\overline L)$ acts transitively on the space
$\Aut^{S^1}(N(\overline L))$ of equivariant symplectic automorphisms of the
normal bundle, with evaluation fibration
\[
\Fix^{S^1}(N(\overline L))\to\Fix^{S^1}(\overline L)\to
\Aut^{S^1}(N(\overline L)),
\]
and $\Aut^{S^1}(N(\overline L))\simeq\U(1)$, while the fibre
$\Fix^{S^1}(N(\overline L))$, the group of equivariant symplectomorphisms that
are the identity near $\overline L$, is contractible. Hence
$\Fix^{S^1}(\overline L)\simeq S^1$.
\end{lemma}
\begin{proof}
Transitivity and the fibration are as in~\cite[Lemma 3.44]{ChPin-Memoirs}, which
also computes $\Aut^{S^1}(N(\overline L))\simeq S^1$. For the fibre, note that an
element of $\Fix^{S^1}(N(\overline L))$ is supported in the complement of a
neighbourhood of $\overline L$. That complement is a star-shaped domain of
$\R^4$: the complement of a line in $(\CP^2,\om_{FS})$ is an open ball of
capacity one, on which the toric action is linear in affine coordinates, so the
circle acts linearly there, in particular near the boundary. Proposition~\ref{prop:equivariant Gromov} therefore
applies and gives $\Fix^{S^1}(N(\overline L))\simeq\{*\}$.
\end{proof}

%%%%%%%%%%%%%%%%%%%%%%%%%%%%%%%%%%%%%%%%%%%%%%%%%%%%%%%%%%%%%%%%%%%%%%%%%%%%%%%%
\subsection{Isolated fixed points}
\label{subsec:CP2 isolated}
%%%%%%%%%%%%%%%%%%%%%%%%%%%%%%%%%%%%%%%%%%%%%%%%%%%%%%%%%%%%%%%%%%%%%%%%%%%%%%%%

Throughout this subsection the fixed point set consists of three points, that
is $(a,b)\notin\{(1,0),(0,1),(1,1)\}$. Let
$p_0$ and $p_1$ be the minimum and the maximum of the moment map, and let
$\mathcal{S}^{S^1}_{L,p_0,p_1}$ be the space of invariant symplectic spheres in
the class $L=[\CP^1]$ through $p_0$ and $p_1$. By
Corollary~\ref{cor:CP2 fixes the vertices}, $\Symp^{S^1}(\CP^2,\om)$ acts on
that space.

\begin{lemma}\label{lemma:CP2 transitive}
The group $\Symp^{S^1}(\CP^2,\om)$ acts transitively on
$\mathcal{S}^{S^1}_{L,p_0,p_1}$, and that space is contractible.
\end{lemma}
\begin{proof}
Given $\overline L_0,\overline L_1$ in the space, choose invariant compatible
almost complex structures $J_0,J_1$ making them holomorphic, join them by a path
in the contractible space $\jj^{S^1}(\om)$, and apply the equivariant
displacement theorem of~\cite[Lemma A.10]{ChPin-Memoirs}, which gives an
equivariant symplectomorphism carrying $\overline L_0$ to $\overline L_1$.

For the contractibility, note that for every $J\in\jj^{S^1}(\om)$ there is a
unique $J$-holomorphic curve in the class $L$ through $p_0$ and $p_1$: existence
is Gromov's theorem, and uniqueness follows from positivity of intersections,
two distinct curves in the class $L$ meeting in a single point. That curve is
invariant, by uniqueness, and the assignment $J\mapsto\overline L$ is a fibration
\[
\jj^{S^1}(\om,\overline L)\to\jj^{S^1}(\om)\to
\mathcal{S}^{S^1}_{L,p_0,p_1}
\]
whose total space and fibre are contractible, the fibre being the set of
invariant structures for which a given curve is holomorphic.
\end{proof}

We write $\Stab^{S^1}(\overline L)$ for the stabilizer of a point of that space.

\begin{prop}\label{prop:CP2 isolated}
If $S^1(a,b)$ has isolated fixed points, then
$\Symp^{S^1}(\CP^2,\om)\simeq\T$.
\end{prop}
\begin{proof}
By Lemma~\ref{lemma:CP2 transitive}, the evaluation fibration
\[
\Stab^{S^1}(\overline L)\to\Symp^{S^1}(\CP^2,\om)\to
\mathcal{S}^{S^1}_{L,p_0,p_1}\simeq\{*\}
\]
gives $\Symp^{S^1}(\CP^2,\om)\simeq\Stab^{S^1}(\overline L)$. The curve
$\overline L$ is not pointwise fixed, by
Lemma~\ref{lemma:CP2 degenerate list}, so
Lemmas~\ref{lemma:CP2 restriction fibration}
to~\ref{lemma:CP2 gauge} give the fibration
\[
\Fix^{S^1}(\overline L)\simeq S^1\to\Stab^{S^1}(\overline L)\to
\Symp^{S^1}(\overline L)\simeq S^1 .
\]
Now $\T$ is contained in $\Stab^{S^1}(\overline L)$, the line $\overline L$
being invariant under the whole torus, and the inclusion is compatible with that
fibration: the subcircle of $\T$ whose direction is normal to the facet
corresponding to $\overline L$ acts trivially on $\overline L$, hence lands in the
fibre $\Fix^{S^1}(\overline L)$, where it induces an equivalence, and the quotient
circle maps to $\Symp^{S^1}(\overline L)$ by an equivalence as well, being the
circle that acts on $\overline L$. The five lemma, applied to the long exact
sequences of homotopy groups of the two fibrations, gives
$\T\simeq\Stab^{S^1}(\overline L)$.
\end{proof}

%%%%%%%%%%%%%%%%%%%%%%%%%%%%%%%%%%%%%%%%%%%%%%%%%%%%%%%%%%%%%%%%%%%%%%%%%%%%%%%%
\subsection{A fixed surface}
\label{subsec:CP2 surface}
%%%%%%%%%%%%%%%%%%%%%%%%%%%%%%%%%%%%%%%%%%%%%%%%%%%%%%%%%%%%%%%%%%%%%%%%%%%%%%%%

Suppose now that $(a,b)\in\{(1,0),(0,1),(1,1)\}$, so that the fixed point set is
a line $\Sigma$ together with one isolated point. No space of curves is needed
here: an equivariant symplectomorphism preserves the fixed point set, hence
preserves $\Sigma$, which is its only two-dimensional component. In other words
\[
\Stab^{S^1}(\Sigma)=\Symp^{S^1}(\CP^2,\om).
\]

\begin{prop}\label{prop:CP2 surface}
If $S^1(a,b)$ fixes a line pointwise, then $\Symp^{S^1}(\CP^2,\om)\simeq\U(2)$.
\end{prop}
\begin{proof}
By Lemma~\ref{lemma:CP2 restriction fibration} applied to $\Sigma$, and by
Lemmas~\ref{lemma:CP2 Symp of the line} and~\ref{lemma:CP2 gauge}, we have the
fibration
\[
\Fix^{S^1}(\Sigma)\simeq S^1\to\Symp^{S^1}(\CP^2,\om)\to
\Symp^{S^1}(\Sigma)=\Symp(\Sigma)\simeq\SO(3),
\]
the curve being pointwise fixed. The group $\U(2)$ is contained in $\Symp^{S^1}(\CP^2,\om)$ here, as the
centralizer of the circle in $\PU(3)$, and the inclusion is compatible with that
fibration: its centre acts trivially on $\Sigma$, hence lands in the fibre
$\Fix^{S^1}(\Sigma)$, where it induces an equivalence, and the quotient
$\PU(2)\cong\SO(3)$ maps to $\Symp(\Sigma)$ by an equivalence. The five lemma
again gives $\U(2)\simeq\Symp^{S^1}(\CP^2,\om)$.
\end{proof}

\begin{cor}\label{Cor:HomotopytypeofSympCP2}
Let $S^1(a,b)\subset\T$ be an effective Hamiltonian circle action on
$(\CP^2,\om)$, normalized so that $a\geq0$. Then
\[
\Symp^{S^1}(\CP^2,\om)\simeq
\begin{cases}
\U(2) & \text{if }(a,b)\in\{(1,0),(0,1),(1,1)\},\\
\T & \text{otherwise.}
\end{cases}
\]
In particular $\Symp^{S^1}(\CP^2,\om)$ is connected in either case.\qed
\end{cor}

These two homotopy types, and the connectedness, are what
Section~\ref{section:Homotopy types of embedding spaces} uses to describe the
spaces of equivariant ball embeddings in $\CP^2$.

\section{From equivariant deformation to isotopy}\label{sec:deformation to isotopy}

This section proves the deformation-to-isotopy lemma required
in~\S\ref{section:role of 3.12 and 3.13} for Hamiltonian circle actions on
$\CP^2$, $\SSS$ and $\CCC$, which completes the proof of path-connectedness. As
in~\S\ref{sec:toric}, the isotopies are obtained by equivariant inflation
(Lemma~\ref{lem:inf}). The invariant curves along which we inflate now come from
a toric extension of the blown-up circle action, and they are followed along a
generic path of invariant almost complex structures by the results
of~\S\ref{section:genericity} and~\S\ref{section:genericity of pairs}. The
inflations themselves are carried out in families in
Proposition~\ref{prop:special}, for the configurations of curves that occur.

\begin{prop}\label{prop:special}
Let a compact Lie group $G$ act homologically trivially on a closed four-manifold
$X$, and let $\mathcal{N}$ be a compact manifold, possibly with boundary. Let
$(\omega_{s,\nu},J_{s,\nu})$, $(s,\nu)\in[0,1]\times\mathcal{N}$, be a smooth
family of invariant symplectic forms and of invariant $\omega_{s,\nu}$-compatible
almost complex structures, and let $\alpha\colon[0,1]\times\mathcal{N}\to[0,1]$
and $c\colon\mathcal{N}\to(0,1)$ be smooth functions such that
$\alpha(0,\nu)=0=\alpha(1,\nu)$. Assume that one of the following holds.
\begin{itemize}
\item[(i)] $X=\CP^2\#\overline{\CP^2}$, the couplings of $\omega_{s,\nu}$ with
$L$ and $E$ are $1$ and $(1-\alpha(s,\nu))\,c(\nu)$, and the class $L-E$ is
represented by an embedded invariant $J_{s,\nu}$-holomorphic sphere.
\item[(ii)] $X$ is $(\SSS)\#\overline{\CP^2}$ or $(\CCC)\#\overline{\CP^2}$, the
couplings of $\omega_{s,\nu}$ with $B^*$, $F$ and $E$ are $\lambda\geq1$, $1$
and $(1-\alpha(s,\nu))\,c(\nu)$, and, for some integer $k\geq0$, each class of
one of the configurations
\begin{equation}\label{eq:b1}
B^*+kF,\quad F-E,\quad B^{**}-kF-E,
\end{equation}
\begin{equation}\label{eq:b2}
B^*+kF-E,\quad F,\quad B^{**}-kF,
\end{equation}
is represented by an embedded invariant $J_{s,\nu}$-holomorphic sphere, where
$k\geq1$ in configuration~\eqref{eq:b2} if $X=(\SSS)\#\overline{\CP^2}$. Here
$B^*=B^{**}=B$ if $X=(\SSS)\#\overline{\CP^2}$, while $B^*=B_{+1}$ and
$B^{**}=B_{+1}-F$ if $X=(\CCC)\#\overline{\CP^2}$.
\end{itemize}
Assume moreover that these spheres depend smoothly on $(s,\nu)$. Then there is a
smooth family of invariant symplectic forms $\omega_{s,t,\nu}$, $0\leq t\leq 1$,
such that $\omega_{s,t,\nu}$ tames $J_{s,\nu}$, $\omega_{s,0,\nu}=\omega_{s,\nu}$,
$\omega_{s,1,\nu}$ lies in the class $[\omega_{0,\nu}]$, and
$\omega_{s,t,\nu}=\omega_{s,\nu}$ for all $t$ whenever $\alpha(s,\nu)=0$.
\end{prop}

\begin{proof}
Each step below is an inflation along one of the given families of spheres,
performed with Lemma~\ref{lem:inf} for all values of $(s,\nu)$ at once. This is
possible because the construction of~\cite{CPP} adapts to smooth families of
data, as noted after Lemma~\ref{lem:inf}, and because the size of each inflation
is an explicit smooth function of $\alpha$, $c$ and the couplings of
$\omega_{s,\nu}$. These sizes vanish where $\alpha$ vanishes. At such points every
inflation therefore has size zero and is performed with an almost complex
structure compatible with the current form, so that it is constant, as in the
proof of Lemma~\ref{lem:toric correction}. This gives the last assertion. We
drop $\nu$ from the notation, and write $\gamma_s=(1-\alpha(s))c$ for the coupling
of $\omega_s$ with $E$.

In case (i), set $\delta_s=(c-\gamma_s)/(1-c)\geq0$. As $L-E$ has
self-intersection $0$, we may inflate $\omega_s$ along it by $\delta_s$. The
couplings of the resulting form with $L$ and $E$ are $1+\delta_s$ and
$\gamma_s+\delta_s=c(1+\delta_s)$, so that dividing the forms of the path by
$1+t\delta_s$, $0\leq t\leq1$, gives the required family.

In case (ii), let $\lambda^*=\omega_s(B^{**})$, which is $\lambda$ if
$X=(\SSS)\#\overline{\CP^2}$ and $\lambda-1$ if $X=(\CCC)\#\overline{\CP^2}$. For
configuration~\eqref{eq:b1}, set $A_0=B^*+kF$, $A_1=F-E$ and
$\delta_0^s=(c-\gamma_s)/(\lambda^*-k-c)$. For configuration~\eqref{eq:b2}, set
$A_0=B^*+kF-E$, $A_1=F$ and $\delta_0^s=(c-\gamma_s)/(1-c)$. In both cases, let
$\delta_1^s=(\lambda^*-k)\,\delta_0^s$. These sizes are non-negative, since
$c-\gamma_s=\alpha(s)c\geq0$, and since $\lambda^*-k-c$ and $\lambda^*-k$ are the
$\omega_0$-areas of the $J_0$-holomorphic spheres in the classes $B^{**}-kF-E$
and $B^{**}-kF$. In configuration~\eqref{eq:b1}, the relation
$(\lambda^*-k-c)\,\delta_0^s=c-\gamma_s$ gives
$\delta_1^s=c(1+\delta_0^s)-\gamma_s$.

We first inflate $\omega_s$ along $A_0$ by $\delta_0^s$. This is allowed because
$A_0\cdot A_0$ is non-negative: it equals $2k$ or $2k-1$, with $k\geq1$ in the
second case, on $(\SSS)\#\overline{\CP^2}$, and $2k+1$ or $2k$ on
$(\CCC)\#\overline{\CP^2}$. Let $\widetilde\omega_s$ be the resulting form. We
then inflate $\widetilde\omega_s$ along $A_1$ by $\delta_1^s$. If $A_1=F$, its
self-intersection is $0$. If $A_1=F-E$, its self-intersection is $-1$, and the
inflation is allowed because $c<1$ gives
\[
\delta_1^s=c(1+\delta_0^s)-\gamma_s<1+\delta_0^s-\gamma_s=\widetilde\omega_s(F-E).
\]
Since $A_0\cdot B^{**}=k$, $A_1\cdot B^{**}=1$, $A_0\cdot F=1$ and $A_1\cdot F=0$,
the final form has couplings $\lambda^*+k\delta_0^s+\delta_1^s=\lambda^*(1+\delta_0^s)$
with $B^{**}$ and $1+\delta_0^s$ with $F$. Its coupling with $E$ is
$\gamma_s+\delta_1^s$ in configuration~\eqref{eq:b1} and $\gamma_s+\delta_0^s$ in
configuration~\eqref{eq:b2}, which is $c(1+\delta_0^s)$ in both cases. As the
classes $B^{**}$, $F$ and $E$ form a basis of $H_2(X;\mathbb{Q})$, dividing the
forms of the concatenated path by $1+t\delta_0^s$, $0\leq t\leq1$, gives the
required family.
\end{proof}

\begin{Remark}\label{Remark:Epreservedthroughinflation}
All the classes along which we inflate in the proof of
Proposition~\ref{prop:special} have non-negative intersection with $E$.
Furthermore, if an invariant sphere in the class $E$ is symplectic for all the
forms $\omega_{s,\nu}$, then it remains symplectic for the forms
$\omega_{s,1,\nu}$.
\end{Remark}

The following theorem gives the path-connectedness asserted in
Theorem~\ref{thm:main-a} and in part~(1) of Theorem~\ref{thm:main-b}.

\begin{thm}\label{thm:main0}
Let $(M,\omega)$ be $(\CP^2,\lambda\omega_{\operatorname{FS}})$, $(\SSS,\oml)$ or
$(\CCC,\oml)$, endowed with a Hamiltonian circle action, and let $\fixedcomp$ be
a connected component of its fixed point set. Then $\EmbS$ is path-connected.
\end{thm}

\begin{proof}
Let $G=S^1$, let $g_0,g_1\in\EmbGiota$, and let $\wt\sigma_t$ be a deformation
given by Proposition~\ref{prop:deformation from embeddings}, with $\delta$ as
there. By Proposition~\ref{prop:isotopy implies connectedness}, it suffices to
show that the truncations $\wt\sigma^a$, $\delta\leq a\leq c$, satisfy the three
conditions listed in~\S\ref{section:role of 3.12 and 3.13}.

\emph{The truncations.} By the proof of
Proposition~\ref{prop:deformation from embeddings}, the truncation $\wt\sigma^a$
is, up to reparametrization, the concatenation of the forms $\wt\sigma_{g_0,a'}$
of Lemma~\ref{lemma:rescaled blow-up}, as $a'$ decreases from $a$ to $\delta$,
with the forms $\wt\sigma_{g_1,a'}$, as $a'$ increases back to $a$. The two paths
agree at $a'=\delta$ by part~(5) of that lemma. This concatenation can be
reparametrized by $s\in[0,1]$, smoothly at the junction as in the proof of
Proposition~\ref{prop:deformation from embeddings}, so as to obtain a smooth
family $\wt\sigma^a_s$, $(s,a)\in[0,1]\times[\delta,c]$. By part~(3) of
Lemma~\ref{lemma:rescaled blow-up}, its couplings with the classes coming from
$M$ are those of $\omega$, and its coupling with $E=[\exdivisor]$ is
$(1-\alpha(s,a))\,a$, where the function $\alpha$ has the following properties.
\begin{itemize}
\item[(a)] It is smooth, with values in $[0,1)$.
\item[(b)] It vanishes at $s=0$ and $s=1$, since the endpoints of
$\wt\sigma^a$ are the forms $\wt\sigma_{g_0,a}$ and $\wt\sigma_{g_1,a}$.
\item[(c)] It vanishes identically for $a=\delta$, since $\wt\sigma^\delta$ is
constant.
\end{itemize}

\emph{Almost complex structures and curves.} By Lemma~\ref{lem:generic1},
Proposition~\ref{prop:genericity space of pairs} and
Corollary~\ref{cor:the set of generic pairs is connected}, together with
Remark~\ref{Remark:GenericforSSSandCCC} when $M=\CP^2$, whose blow-up is $\CCC$,
there are invariant almost complex structures $J^i_{a'}$ compatible with
$\wt\sigma_{g_i,a'}$, depending smoothly on $a'\in[\delta,c]$, with
$J^0_\delta=J^1_\delta$, such that all the pairs
$(\wt\sigma_{g_i,a'},J^i_{a'})$ are generic in the sense of Karshon and $J^0_c$
is any prescribed generic structure. We choose $J^0_c$ as follows. By
Theorem~\ref{thm:extoric}, the circle action on $(\wM,\wt\sigma_{g_0,c})$ extends
to a Hamiltonian $\T$-action, and by
Theorems~\ref{Theorem:ToricExtension=intersectionwithstrata}
and~\ref{Theorem:IntersectionwithStrataanot1}, Lemma~\ref{Lemma:Atmost2tori},
Remark~\ref{remark:open stratum = less negative}(2) and
Proposition~\ref{prop:Fredholm genericity implies Karshon}, this extension can
be chosen so that it preserves a $\wt\sigma_{g_0,c}$-compatible complex structure
$J$ for which the pair $(\wt\sigma_{g_0,c},J)$ is generic. We take $J^0_c=J$. The
invariant spheres of the toric action, described in~\S\ref{subsec:ruled}
and~\S\ref{subsec:blowup}, are fixed by circle subgroups of $\T$, hence
$J$-holomorphic by Lemma~\ref{LR}. They include a sphere in the class $L-E$ if
$M=\CP^2$, and spheres in the three classes of one of the
configurations~\eqref{eq:b1} and~\eqref{eq:b2} if $M$ is $\SSS$ or $\CCC$. By
Proposition~\ref{prop:persist}, these classes are represented by embedded
invariant $J^i_{a'}$-holomorphic spheres for all $i$ and $a'$.

\emph{Inflation.} Composing the families $J^i_{a'}$ and the corresponding
spheres with the reparametrization above gives invariant almost complex
structures $J_{s,a}$ compatible with $\wt\sigma^a_s$, and invariant
$J_{s,a}$-holomorphic spheres, that depend smoothly on $(s,a)$.
Proposition~\ref{prop:special} applies to the family
$(\wt\sigma^a_s,J_{s,a})$, with $\mathcal{N}=[\delta,c]$ and the function $\alpha$
above. If $M=\CP^2$, case~(i) applies to the forms divided by $\lambda$, with
$c(a)=a/\lambda$. If $M$ is $\SSS$ or $\CCC$, case~(ii) applies with $c(a)=a$. We
obtain invariant forms $\sigma_{s,u,a}$ taming $J_{s,a}$, with
$\sigma_{s,0,a}=\wt\sigma^a_s$ and $\sigma_{s,1,a}$ in the class of
$\wt\sigma^a_0$, which by properties~(b) and~(c) coincide with
$\wt\sigma^a_s$ for all $u$ when $s=0$, $s=1$ or $a=\delta$. The path $u\mapsto\sigma_{\cdot,u,a}$ is therefore a
homotopy relative to the endpoints from $\wt\sigma^a$ to the isotopy
$\wt\tau^a=\sigma_{\cdot,1,a}$, which depends continuously on $a$ and is constant
for $a=\delta$. By Remark~\ref{Remark:Epreservedthroughinflation}, the
exceptional divisor $\exdivisor$, which is symplectic for all the forms
$\wt\sigma^a_s$, remains symplectic along $\wt\tau^a$. The three conditions
of~\S\ref{section:role of 3.12 and 3.13} are thus satisfied.
\end{proof}

%%%%%%%%%%%%%%%%%%%%%%%%%%%%%%%%%%%%%%%%%%%%%%%%%%%%%%%%%%%%%%%%%%%%%%%%%%%%%%%%
\section{Homotopy type of space of equivariant embeddings}\label{section:Homotopy types of embedding spaces}\label{section:embedding spaces}
%%%%%%%%%%%%%%%%%%%%%%%%%%%%%%%%%%%%%%%%%%%%%%%%%%%%%%%%%%%%%%%%%%%%%%%%%%%%%%%%

We now use the results from the previous sections to describe the homotopy type of the space of equivariant embeddings centered on a connected component of a fixed point set $\fixedcomp$. 

%%%%%%%%%%%%%%%%%%%%%%%%%%%%%%%%%%%%%%%%%%%%%%%%%%%%%%%%%%%%%%%%%%%%%%%%%%%%%%%%
\subsection{Embeddings in $\CP^2$}
\label{subsec:CP2 embeddings}
%%%%%%%%%%%%%%%%%%%%%%%%%%%%%%%%%%%%%%%%%%%%%%%%%%%%%%%%%%%%%%%%%%%%%%%%%%%%%%%%

Let $(\CP^2,\om)$ be normalized so that a line has area one and let $\alpha=S^1(a,b)\subset\T$ with $a\geq0$. Recall that the fixed point set consists in either three isolated fixed points or, for any of the three degenerate actions $S^1(1,0)$, $S^1(0,1)$, $S^1(1,1)$ of Lemma~\ref{lemma:CP2 degenerate list}, one isolated point together with a fixed line.
Let $\fixedcomp$ be a component of the fixed point set. Recall also that the centralizer of the circle in the
K\"ahler isometry group $\PU(3)$ is $\T$ in the isolated case, and $\U(2)$ is the degenerate case. Write $K$ for this centralizer. By Corollary~\ref{Cor:HomotopytypeofSympCP2} the inclusion $K\into\Symp^{S^1}(\CP^2,\om)$ is a homotopy equivalence.

We now make three observations. First, the moment triangle has vertices $(0,0)$, $(1,0)$, $(0,1)$ so that  $c^{\T}(\CP^2,p)=1$ at every vertex and every $c<1$ is admissible. Second, computations similar to Lemma~\ref{lemma:divisor label} show that, after blow-up, the label of the exceptional divisor in the Karshon graph is  
\[
|b-a| \ \text{ at } (0,0), \qquad |b| \ \text{ at } (1,0), \qquad
|a| \ \text{ at } (0,1).
\]
In particular, for a  degenerate action, the label is $0$ if $\fixedcomp$ is the isolated fixed point, and it is $1$ if $\fixedcomp$ is the fixed line. Third, blowing up a ball of capacity $c$ gives $\CCC$ with $\om(L)=1$ and $\om(E)=c$. The fibre class is $F:=L-E$, of area $1-c$, and the exceptional divisor is a section of self-intersection $-1$. In the normalisation of~\S\ref{subsec:ruled}, the symplectic size of the section with self-intersection $+1$ is then
\[
\hat\lambda=\frac{1}{1-c}>1.
\]
The class $E$ is the only exceptional class of $\CCC$, and its possible decompositions are of the form $E=E-kF$ for $1\leq k \leq\hat{\ell}$. Hence, $E$ is indecomposable if and only if $c\leq\frac12$.

\begin{prop}\label{prop:CP2 embeddings}
For every $c\in(0,1)$, $\IEmb^{S^1}(\CP^2,\fixedcomp,c)\ \simeq\ \fixedcomp$.
\end{prop}
\begin{proof}
\emph{The stabilizer on the blow-up.} By Proposition~\ref{prop:homotopy equivalence Stab(B) and Stab(E)}, there is a homotopy equivalence $\Stab^{S^1}(\iota)\simeq\Symp^{S^1}(\CCC,Z,\tom)$. The blow-up of $\T$ at a vertex belonging to $\fixedcomp$ is a toric extension $\ttT$ of $\widehat\alpha_{\fixedcomp,c}$. Its moment image is the polygon of $\T_1$, for every $c<1$. Its integrable structure $J_1$ is invariant, and $\exdivisor$, being an edge, is the unique embedded $J_1$-holomorphic representative of the exceptional class $[\exdivisor]=E=B-F=D_1$. By the description of the strata of $\jj_{\om_{\hat\lambda}}$ in \cite[Theorem 3.18]{ChPin-Memoirs}, the stratum $U_1$ consists precisely of those structures for which $D_1$ has an embedded holomorphic representative. Hence, for any $0<c<1$,
\[
\mathcal{J}^{S^1}_{\tom}([\exdivisor])=\jj^{S^1}_{\tom}\cap U_1.
\]
That stratum contains $J_1$, and by \cite[Theorem 6.1]{ChPin-Memoirs} it is a single orbit of $\Symp^{S^1}_h(\CCC,\tom)$. The stabilizer of $J_1$ is $\ttT$ itself unless the divisor $\exdivisor$ is
pointwise fixed (i.e. for the actions $S^1(0,\pm1;m)$) in which case the stabilizer is $\U(2)$. Write $\tilde{K}$ for that stabilizer. By part~(2) of Lemma~\ref{lemma:indecomposable exceptional class} the inclusion $\tilde{K}\into\Symp^{S^1}(\CCC,\exdivisor,\tom)$ is a homotopy equivalence. 
Thus $\Stab^{S^1}(\iota)\simeq \tilde{K}$, for every $c\in(0,1)$.

\emph{Forming the quotient.}  Blowing down maps $\tilde{K}$ into $K$, and the square formed by the two inclusions and the blow-down commutes. There are two cases to consider.

If $\fixedcomp$ is an isolated fixed point, then $\Symp^{S^1}(\CP^2,\fixedcomp,\om)=\Symp^{S^1}(\CP^2,\om)$ by
Corollary~\ref{cor:CP2 fixes the vertices}, and $\tilde{K}\to K$ is an isomorphism: both stabilizers are two-tori when $\alpha$ is not degenerate, and both are $\U(2)$ when it is. Hence the inclusion $\Stab^{S^1}(\iota)\into\Symp^{S^1}(\CP^2,\fixedcomp,\om)$ is a homotopy equivalence, and the evaluation fibration makes
$\IEmb^{S^1}(\CP^2,\fixedcomp,c)$ weakly contractible, hence contractible by Remark~\ref{rem:CW}.

If $\fixedcomp$ is a fixed line $\Sigma$, then $\alpha$ is degenerate, so $K=\U(2)$. The stabilizer $\tilde{K}$ is a two-torus, mapped by the blow-down onto a maximal torus of $\T\subset K$. Moreover
$\Symp^{S^1}(\CP^2,\Sigma,\om)=\Symp^{S^1}(\CP^2,\om)$, every equivariant symplectomorphism preserving the fixed point set and hence its only two-dimensional component. The evaluation fibration therefore has the homotopy type of $\T\to\U(2)\to\U(2)/\T$, whose base is $\CP^1\simeq \fixedcomp$.
\end{proof}

\begin{Remark}\label{rem:CP2 large ball}
Proposition~\ref{prop:CP2 embeddings} shows that the homotopy type of $\IEmb^{S^1}(\CP^2,\fixedcomp,c)$ is independent of the number of toric extensions of $\widehat\alpha_{\fixedcomp,c}$. It is one whenever $c\leq\frac12$ since $\hat\lambda\leq2$, and by Lemma~\ref{lemma:list of tori} the only toric action on $(\CCC,\om_{\hat\lambda})$ is $\T_1$. It can also be two, as the following example shows. Take $\alpha=S^1(1,2)$, and blow-up at the vertex $(0,0)$ with capacity $c=\tfrac9{10}$ so that $\hat\lambda=10$. The label of $\exdivisor$ is $|b-a|=1$. The blown-up action has the four fixed points
\[
p'\ (\tfrac{9}{10};\{1,1\}),\quad S\ (1;\{-1,1\}),\quad p''\ (\tfrac{9}{5};\{2,-1\}),\quad
Q\ (2;\{-2,-1\}),
\]
and a single $\Z_k$-sphere with $k\geq2$, namely $p''Q$, of area $0.1$. The Delzant quadrilateral with vertices
\[
(0,0),\quad (\tfrac{1}{10},\tfrac{1}{10}),\quad (\tfrac{9}{10},\tfrac{1}{10}),\quad (\tfrac{11}{10},0)
\]
has the same graph for its subcircle $\mu=x$, so it extends the same action. Its self-intersections are
$(3,0,-3,0)$, so it is $\T_3$, which acts on $(\CCC,\om_{\hat\lambda})$ because $2\hat\lambda>4$. 
\end{Remark}

%%%%%%%%%%%%%%%%%%%%%%%%%%%%%%%%%%%%%%%%%%%%%%%%%%%%%%%%%%%%%%%%%%%%%%%%%%%%%%%%
\subsection{Embeddings in rational ruled surfaces}
\label{subsec:ruled embeddings}
%%%%%%%%%%%%%%%%%%%%%%%%%%%%%%%%%%%%%%%%%%%%%%%%%%%%%%%%%%%%%%%%%%%%%%%%%%%%%%%%

Let $(M,\oml)$ be a rational ruled surface, $\alpha$ a Hamiltonian circle action on it with an isolated fixed point $p$, let $c$ be the capacity of an equivariant ball centred at $p$ and let $\iota\in\Emb^{S^1}(M,p,c)$. We write
\[
N=\bigl|\Ext(\alpha)\bigr|,
\qquad
\widehat N=\bigl|\Ext(\widehat\alpha_{p,c})\bigr| .
\]
for the numbers of toric extensions.

\begin{thm}\label{thm:embeddings by extensions}
When $\fixedcomp=\{p\}$ is an isolated fixed point in a rational ruled surface,
\[
\IEmb^{S^1}(M,p,c)\simeq
\begin{cases}
\Omega S^3\times S^1 & \text{if }(N,\widehat N)=(2,1),\\
\{*\} & \text{otherwise}.
\end{cases}
\]
Moreover $\Emb^{S^1}(M,p,c)\to\IEmb^{S^1}(M,p,c)$ is a principal bundle whose structure group $\Symp^{S^1}(B^4(c))\simeq\U^{S^1}(2)$ is $\T$ when the two weights of the circle at $p$ are distinct, and $\U(2)$ when they are equal.
\end{thm}
The proof splits into three cases.

\subsubsection{Proof of Theorem~\ref{thm:embeddings by extensions} in the indecomposable case}
Assume that $M=\SSS$, or that $M=\CCC$ with $c\leq\lambda-1$, so that the class of the exceptional divisor is indecomposable and Proposition~\ref{prop:capacity rule} applies as an equality. We exclude the action $S^1(1,1;0)$ with $\lambda=1$, treated in the second case below. Under these assumptions, $\Symp^{S^1}(M,p,\oml)=\Symp^{S^1}_h(M,\oml)$, the centralizer being connected and $p$ isolated. 

By Theorem~\ref{thm:stabilizer} the stabilizer $\Stab^{S^1}(\iota)$ is homotopy equivalent to $\Symp^{S^1}_h(\bX_2,\omlc)$, and by Theorems~\ref{Thm:Singlestratumhomotopytype} and~\ref{Thm:pushout centralizer SSS} the two centralizers are
\[
\Symp^{S^1}_h(M,\oml)\simeq
\begin{cases}\T & \text{if } N=1,\\ \Omega S^3\times\mathbb{T}^3 & \text{if } N=2,\end{cases}
\qquad
\Symp^{S^1}_h(\bX_2,\omlc)\simeq
\begin{cases}\T & \text{if } \widehat N=1,\\ \Omega S^3\times\mathbb{T}^3 & \text{if } \widehat N=2.\end{cases}
\]
The inequality $\widehat N\leq N$ holds because, under the assumptions, Proposition~\ref{prop:capacity rule} is an equality and each torus of $\Ext(\alpha)$ contributes at most one extension when $p$ is isolated. The pair
$(N,\widehat N)=(1,2)$ is therefore excluded.

Consider now the evaluation fibration
\[
\Stab^{S^1}(\iota)\longrightarrow\Symp^{S^1}(M,p,\oml)\longrightarrow\IEmb^{S^1}(M,p,c).
\]
For the pair $(2,1)$ the fibre is a two-torus and the total space is $\Omega S^3\times\mathbb{T}^3$. Consider the commutative diagram
\[
\begin{tikzcd}[column sep=normal]
\{*\}\ar{d} \ar{r}& \Symp^{S^1}_h(\bX_2,\omlc)\simeq\T \ar{d} \ar{r}& \T \ar{d}\\
\Omega S^3 \ar{d} \ar{r}& \Symp^{S^1}_h(M,\oml) \ar{r}{f_P} \ar{d}& \T\times S^1 \ar{d}\\
\Omega S^3 \ar{r}& \IEmb^{S^1}(M,p,c) \ar{r}& (\T\times S^1)/\T\simeq S^1
\end{tikzcd}
\]
whose first column is the fibre of the evaluation fibration, and whose second column is that fibration. By~\cite[\S4.2.3]{ChPin-Memoirs}, $f_P$ is the map obtained by viewing $\Symp^{S^1}_h(M,\oml)$ as the homotopy pushout of the two tori extending $\alpha$. This map $f_P$ admits a section, by the same reference, and that section descends to the quotient, so that the bottom row is a trivial fibration and $\IEmb^{S^1}(M,p,c)\simeq\Omega S^3\times S^1$.

For the pairs $(1,1)$ and $(2,2)$ the quotient is contractible because the inclusion $\Stab^{S^1}(\iota)\into\Symp^{S^1}(M,p,\oml)$ is itself a homotopy equivalence. Indeed, 
when $N=\widehat N=1$, the unique toric extension $\ttT$ of $\widehat\alpha$ blows down to the unique toric extension $\T$ of $\alpha$, by Proposition~\ref{prop:capacity rule}, and the blow-down is an isomorphism of tori, so
that the square
\[
\begin{array}{ccc}
\ttT & \lhook\joinrel\longrightarrow & \Stab^{S^1}(\iota) \simeq \Symp^{S^1}_h(\bX_2,\omlc)\\[2pt]
\cong\big\downarrow & & \big\downarrow\\[2pt]
\T & \lhook\joinrel\longrightarrow & \Symp^{S^1}(M,p,\oml)
\end{array}
\]
commutes. The two horizontal maps are homotopy equivalences, by Theorem~\ref{Thm:Singlestratumhomotopytype} on the top and \cite[Theorems 6.7 and 4.8]{ChPin-Memoirs} for the bottom inclusion. 

When $N=\widehat N=2$, the two toric extensions of $\widehat\alpha$ blow down to the two toric extensions of $\alpha$, again by Proposition~\ref{prop:capacity rule}, so that the inclusion carries the diagram $\ttT_i\hookleftarrow S^1\hookrightarrow\ttT_j$ to the diagram $\T_i\hookleftarrow S^1\hookrightarrow\T_j$ by isomorphisms. Consequently, it induces a homotopy equivalence of the homotopy pushouts. Since the canonical maps out of those pushouts are homotopy equivalences, by Theorem~\ref{Thm:pushout centralizer SSS} on the left and
by~\cite[Theorem 4.7]{ChPin-Memoirs} on the right, the inclusion is one as well.

Finally $\Emb^{S^1}(M,p,c)\to\IEmb^{S^1}(M,p,c)$ is a principal bundle with structure group the group of equivariant symplectomorphisms of the ball, which is $\U^{S^1}(2)$: the diagonal torus when the two weights of the isotropy
representation at $p$ are distinct, and $\U(2)$ when they are equal.

\subsubsection{Proof of Theorem~\ref{thm:embeddings by extensions} for $S^1(1,1;0)$ on $(\SSS,\om_1)$}
\label{subsubsec:exceptional action}

For the action $S^1(1,1;0)$ on $(\SSS,\om_1)$ the centralizer $\Symp^{S^1}(M,\oml)$ is not connected, so that it differs from $\Symp^{S^1}(M,p,\oml)$. The latter is homotopy equivalent to $\T$, by the proof of~\cite[Theorem 4.1]{ChPin-Memoirs}, and so is $\Symp^{S^1}_h(\bX_2,\omlc)$. Proceeding as before, the evaluation fibration makes $\IEmb^{S^1}(M,p,c)$ contractible in that case as well.

\subsubsection{Proof of Theorem~\ref{thm:embeddings by extensions} in the decomposable case}
\label{subsubsec:decomposable case}

Assume finally that $M=\CCC$ and $c>\lambda-1$. By the proof of Corollary~\ref{cor:capacity is max}, this forces $1<\lambda<2$, so that $\Ext(\alpha)=\{\T_1\}$ and $N=1$. As in that proof, we may assume that $\alpha=S^1(a,b;1)$ with $a\geq0$, and the fixed point $p$ is then $P$ or $S$. The symmetry of $\Delta_1$ used in that proof exchanges $P$ and $S$, so that by Lemma~\ref{lemma:invariance} we may assume $p=P$. The weights of $\alpha$ at $P$ are $a$ and $b$, both non-zero since $P$ is isolated, hence $a\geq1$. 

We first observe that, in this regime, the two rulings of the blow-up exchange their roles with respect to Lemma~\ref{lemma:two blow-downs}. 
\begin{lemma}\label{lemma:decomposable normalisation}
Let $\lambda-1<c<1$. In the normalisation of~\S\ref{subsec:blowup}, the blow-up of $(\CCC,\oml)$ along a ball of capacity $c$ is, after rescaling by $1/(\lambda-c)$, the manifold $(\bX_2,\om_{\lambda',c'})$ with
\[
\lambda'=\frac{1}{\lambda-c}>1
\qquad\text{and}\qquad
c'=\frac{1-c}{\lambda-c}\in(0,1),
\]
and its exceptional divisor $\exdivisor$ represents the class $D_0=B-E$ of~\eqref{eq:Di}.
\end{lemma}
\begin{proof}
In the basis $L,E_1,E_2$ of~\S\ref{subsec:blowup}, where $E_2$ is the class of the exceptional divisor $\exdivisor$ and $E_1=B_{-1}$ is the class of the section of self-intersection $-1$ of $\CCC$, the areas of $L$, $E_1$ and $E_2$ are
$\lambda$, $\lambda-1$ and $c$. The classes
\[
B':=L-E_1,\qquad F':=L-E_2,\qquad E':=E_0=L-E_1-E_2
\]
have areas $1$, $\lambda-c$ and $1-c$. Since $\lambda-c<1$, the class of smaller area among $B'$ and $F'$ is $F'$, so that the normalisation of~\S\ref{subsec:blowup} has fibre class $F'$, section class $B'$ and exceptional class $E'$. Dividing the areas by $\lambda-c$ gives $\lambda'$ and $c'$, and $c'<1$ because $\lambda>1$. Finally $E_2=B'-E'$.
\end{proof}

In particular, the class $[\exdivisor]$ of the exceptional divisor is decomposable, and the space of almost complex structures on $(\bX_2,\om_{\lambda',c'})$ admits several strata, characterised by the possible degenerations of $[\exdivisor]$. Another consequence is that the blown-up action $\widehat\alpha$ may have two toric extensions, see Remark~\ref{rem:corner counterexample}. Consequently, we no longer have in general the equivalence $\Symp^{S^1}_h(\bX_2,\exdivisor,\tom)\simeq \Symp^{S^1}_h(\bX_2,\tom)$ of the stabilizer with the whole centralizer.

We choose the ball as follows. By Lemma~\ref{lemma:capacity basics}, $c^{\T_1}(M,P)=1>c$, so that cutting the corner $P$ of $\Delta_1$ gives a $\T_1$-equivariant ball $\iota_0$ of capacity $c$ centred at $P$. As $\Emb^{S^1}(M,P,c)$ is path-connected by Theorem~\ref{thm:main-b}, we may take $\iota=\iota_0$, so that $\T_1$ is contained in $\Stab^{S^1}(\iota)$. Its lift to the blow-up is the toric action $\T_1(P)$, which extends $\widehat\alpha$ by Proposition~\ref{prop:capacity rule}. In the notation of Lemma~\ref{lemma:decomposable normalisation}, the five edges of its polygon carry the classes $D_0=E_2$, $E'$, $E_1=F'-E'$, $B'$ and $F'$, of self-intersections $-1$, $-1$, $-1$, $0$ and $0$. Hence $\T_1(P)$ is the toric action $\ttT_0$, and its integrable structure $\tilde J_0$ lies in the stratum $\tilde U_0$, the exceptional divisor being a $\tilde J_0$-holomorphic sphere in the class $D_0$.

By Lemma~\ref{lemma:decomposable normalisation} and the definition of the strata, the space of invariant compatible structures for which $[\exdivisor]$ has an embedded holomorphic representative is the open stratum
\[
\jj^{S^1}_{\tom}([\exdivisor])=\jj^{S^1}_{\tom}\cap\tilde U_0 ,
\]
which contains $\tilde J_0$. By Theorem~\ref{Thm:Homogenousstrata}, the orbit map identifies this stratum, up to homotopy, with the homogeneous space $\Symp^{S^1}_h(\bX_2,\tom)/\ttT_0$, where $\ttT_0$ is the stabilizer of $\tilde J_0$, that is, the group of equivariant K\"ahler isometries of $(\tom,\tilde J_0)$ acting trivially on homology. Part~(2) of Lemma~\ref{lemma:indecomposable exceptional class} then shows that the inclusion
\[
\ttT_0\into\Symp^{S^1}_h(\bX_2,\exdivisor,\tom)
\]
is a homotopy equivalence. Every element of $\Symp^{S^1}(\bX_2,\exdivisor,\tom)$ acts trivially on homology, since it fixes $E_2=[\exdivisor]$, permutes the exceptional classes $E_0$, $E_1$, $E_2$ of~\S\ref{subsec:blowup}, and cannot exchange $E_0$ and $E_1$, as $E_0\cdot E_2=1$ while $E_1\cdot E_2=0$. Hence $\Symp^{S^1}(\bX_2,\exdivisor,\tom)=\Symp^{S^1}_h(\bX_2,\exdivisor,\tom)$. By Proposition~\ref{prop:homotopy equivalence Stab(B) and Stab(E)}, $\Stab^{S^1}(\iota)$ is homotopy equivalent to $\Symp^{S^1}(\bX_2,\exdivisor,\tom)$. Moreover, since the blow-down identifies the tori $\ttT_0$ and $\T_1$, we obtain the homotopy equivalence $\T_1\into\Stab^{S^1}(\iota)$.

On the other hand, since the class $B_{+1}-2F$ has negative area $\lambda-2$, the space of compatible almost complex structures on $(\CCC,\oml)$ reduces to the single stratum $U_1$ of~\cite[Theorem 3.18]{ChPin-Memoirs}. As $a\neq0$, \cite[Theorem 6.1]{ChPin-Memoirs} then shows that the inclusion $\T_1\into\Symp^{S^1}_h(\CCC,\oml)$ is a homotopy equivalence, the space $\jj^{S^1}_{\oml}$ being contractible. Since every symplectomorphism of $\CCC$ acts trivially on homology~\cite[Lemma 3.21]{ChPin-Memoirs}, and since the centralizer $\Symp^{S^1}(\CCC,\oml)$, being connected, fixes the isolated fixed point $P$, the inclusion $\T_1\into\Symp^{S^1}(M,P,\oml)$ is a homotopy equivalence. In the commutative triangle formed by the inclusions
\[
\T_1\into\Stab^{S^1}(\iota)\into\Symp^{S^1}(M,P,\oml),
\]
two of the three maps are homotopy equivalences, hence so is the third. The evaluation fibration
\[
\Stab^{S^1}(\iota)\longrightarrow\Symp^{S^1}(M,P,\oml)\longrightarrow
\IEmb^{S^1}(M,P,c)
\]
then shows that $\IEmb^{S^1}(M,P,c)$, which is path-connected, is weakly contractible, hence contractible by Remark~\ref{rem:CW}.

Finally, the last assertion of Theorem~\ref{thm:embeddings by extensions} is proved as in the indecomposable case, the structure group being $\U(2)$ for $\alpha=S^1(1,1;1)$ and $\T$ otherwise.\qed

\subsection{Reading the homotopy type off the data}
\label{section:reading}

It remains to compute the pair $(N,\widehat N)$ from the data. By Theorem~\ref{thm:extoric} we may write $\alpha=S^1(a,b;r)$ in a presentation for which $p$ is one of the four vertices $P,Q,R,S$ of the polygon $\Delta_r$
of~\eqref{eq:polygon}. Note that, by Lemma~\ref{lemma:invariance}, nothing below depends on which presentation is chosen. The fixed point being isolated, $a\neq0$, and $b\notin\{0,r\}$ when $a=\pm1$. We write $n=|2b-r|$ and $\epsilon_n=\epsilon_r$.

\medskip
\noindent\textbf{Step 1: special cases.} 
If $M=\CCC$ and $c>\lambda-1$, or if the action is $S^1(1,1;0)$ on $(\SSS,\om_1)$, then $\IEmb^{S^1}(M,p,c)$ is contractible, by~\S\ref{subsubsec:decomposable case} and~\S\ref{subsubsec:exceptional action} respectively. The next steps concern the remaining cases.

\medskip
\noindent\textbf{Step 2: admissible capacities.} By Proposition~\ref{prop:capacity rule} and Corollary~\ref{cor:capacity is max}, an equivariant ball of capacity $c$ centred at $p$ exists if and only if
\[
c<c^{S^1}(M,p)=\max_{\T\in\Ext(\alpha)}c^{\T}(M,p),
\]
each $c^{\T}(M,p)$ being the smaller of the areas of the two invariant spheres through $p$, that is the smaller of the lengths of the two edges of the polygon of $\T$ at $p$. This is always $1$ at a vertex of the positive section, and $\lambda-\frac{m+\epsilon_m}{2}$ at a vertex of the negative section of $\T_m$ when that number is smaller than $1$.

\medskip
\noindent\textbf{Step 3: $N=|\Ext(\alpha)|$.} By Proposition~\ref{prop:two tori},
\[
N=\begin{cases}
2 & \text{if } a=\pm1 \text{ and } 2\lambda>n+\epsilon_n,\\
1 & \text{otherwise.}
\end{cases}
\]

\medskip
\noindent\textbf{Step 4: $\widehat N=|\Ext(\widehat\alpha_{p,c})|$.} By Theorem~\ref{Thm:Intersectionwithstrata_a=1} and Table~\ref{table:NumberOfStrata}, each of the tori of $\Ext(\alpha)$ survives the blow-up according to the position of $p$: without condition at a vertex of its positive section, and under $2\lambda>m+\epsilon_m+2c$ at a vertex of the negative section of $\T_m$. The table gives, for $a=\pm1$, which of the two possibilities occurs for each of the four vertices, according to the sign of $2b-r$.

\medskip
\noindent\textbf{Step 5: the two groups.} By \cite[Theorems 4.8 and 6.7]{ChPin-Memoirs} and by Theorems~\ref{Thm:Singlestratumhomotopytype} and~\ref{Thm:pushout centralizer SSS},
\[
\Symp^{S^1}_h(M,\oml)\simeq
\begin{cases}\T & \text{if } N=1,\\ \Omega S^3\times\mathbb{T}^3 & \text{if } N=2,\end{cases}
\qquad
\Symp^{S^1}_h(\bX_2,\omlc)\simeq
\begin{cases}\T & \text{if } \widehat N=1,\\ \Omega S^3\times\mathbb{T}^3 & \text{if } \widehat N=2,\end{cases}
\]
and by Theorem~\ref{thm:stabilizer} the second group is the stabilizer $\Stab^{S^1}(\iota)$ of the ball.

\medskip
\noindent\textbf{Step 6: the embedding spaces.}
Theorem~\ref{thm:embeddings by extensions} now reads the answer off the pair $(N,\widehat N)$. For the parametrised embeddings, the two weights of the circle at $p$, which decide the structure group, are the pairings of $(a,b)$ with the two edge directions at $p$ listed before Table~\ref{table:weights}. They are $\{a,b\}$ at $P$, $\{a,-b\}$ at $Q$, $\{ar-b,-a\}$ at $R$ and $\{-a,b-ar\}$ at $S$.

\begin{exmp}
Take $r=4$, $a=1$, $b=1$ on $(\SSS,\om_{2.1})$, so that $n=2$ and $\Lambda=\tfrac{41}{10}$, and let $p=Q$, a vertex of the negative section of $\T_4$, which has area $\tfrac{1}{10}$. By Lemma~\ref{lemma:dictionary} the vertex $Q$ corresponds to a vertex of the positive section of $\Delta_2$, so $c^{\T_4}(M,Q)=\tfrac{1}{10}$ while $c^{\T_2}(M,Q)=1$ and $c^{S^1}(M,Q)=1$. Take $c=1/2$. Step 3 gives $N=2$, since $2\lambda=\tfrac{42}{10}>2=n+\epsilon_n$. Step 4 gives $\widehat N=1$: the first torus is lost after blowing up, because $2\lambda=\tfrac{42}{10}<5=r+\epsilon_r+2c$. The second survives without condition. Hence the pair is $(2,1)$ and $\IEmb^{S^1}(M,Q,1/2)\simeq\Omega S^3\times S^1$. At the vertex $P$ of the same action, by contrast, both tori survive the blow-up and the space of embeddings is contractible.
\end{exmp}

\appendix
%%%%%%%%%%%%%%%%%%%%%%%%%%%%%%%%%%%%%%%%%%%%%%%%%%%%%%%%%%%%%%%%%%%%%%%%%%%%%%%%
\section{Equivariant Gauge Groups}\label{sec:equivariant gauge groups}
%%%%%%%%%%%%%%%%%%%%%%%%%%%%%%%%%%%%%%%%%%%%%%%%%%%%%%%%%%%%%%%%%%%%%%%%%%%%%%%%

Let $K$ be an abelian Lie group acting on the right on a principal bundle $\pi:P\to B$, and let $G$ be a Lie group acting on $B$, the action being lifted to a left action on $P$. The action of $G$ on $B$ need not be effective. Let $\Gg^G(P)$ denote the group of $G$-equivariant bundle automorphisms of $P$, that is, of the equivariant maps $u:P\to P$ with $\pi\circ u=\pi$, and let $\Inv(B,K)$ denote the space of $G$-invariant smooth maps $B\to K$. Both carry the $C^\infty$ topology.

\begin{lemma}\label{lemma:gauge group as invariant maps}
The assignment $u\mapsto\phi_u$, where $\phi_u:P\to K$ is determined by $x\cdot\phi_u(x)=u(x)$, is a homeomorphism $\Gg^G(P)\cong\Inv(B,K)$.
\end{lemma}

\begin{proof}
The map $\phi_u$ is well defined because $K$ acts freely and transitively on the fibres. It is $G$-invariant. Indeed, comparing
\[
u(h\cdot x)=h\cdot x\cdot\phi_u(h\cdot x)
\qquad\text{and}\qquad
u(h\cdot x)=h\cdot u(x)=h\cdot x\cdot\phi_u(x)
\]
and using freeness again gives $\phi_u(h\cdot x)=\phi_u(x)$. It is constant on the fibres, since $\phi_u(x\cdot k)=k^{-1}\phi_u(x)k=\phi_u(x)$ because $K$ is abelian, so it descends to a $G$-invariant map $B\to K$. Conversely, a $G$-invariant map $B\to K$ defines an equivariant automorphism of $P$ through the same formula, and the two assignments are continuous and inverse to one another.
\end{proof}

\begin{lemma}\label{lemma:invariant circle valued maps}
Assume in addition that $B$ is connected and simply connected, that $K=S^1$ and that $G$ is compact. Then $\Inv(B,S^1)$ deformation retracts onto its subgroup of constant maps, so that $\Inv(B,S^1)\simeq S^1$.
\end{lemma}

\begin{proof}
Since $B$ is simply connected, every smooth $f:B\to S^1$ can be written $f=e^{2\pi i\tilde f}$ with $\tilde f:B\to\R$ smooth, uniquely up to an integer constant. If $f$ is $G$-invariant then, for each $g\in G$, the function $x\mapsto\tilde f(gx)-\tilde f(x)$ is continuous, integer valued, hence equal to a constant $n(g)$. Since $n:G\to\Z$ is a homomorphism, it is trivial. Thus $\tilde f$ is $G$-invariant and $\Inv(B,S^1)\cong C^\infty(B)^G/\Z$,
the integers acting by addition of constants. Fix $x_0\in B$. The homotopy $\tilde f\mapsto(1-t)\tilde f+t\,\tilde f(x_0)$ is a deformation retraction of $C^\infty(B)^G$ onto the constants. Since it commutes with the addition of integers, it descends to a deformation retraction of $\Inv(B,S^1)$ onto $\R/\Z=S^1$.
\end{proof}

%%%%%%%%%%%%%%%%%%%%%%%%%%%%%%%%%%%%%%%%%%%%%%%%%%%%%%%%%%%%%%%%%%%%%%%%%%%%%%%%
\section{Graphs}
%%%%%%%%%%%%%%%%%%%%%%%%%%%%%%%%%%%%%%%%%%%%%%%%%%%%%%%%%%%%%%%%%%%%%%%%%%%%%%%%

We list half of the decorated graphs obtained by blowing up the $S^1(a,b,r)$ actions on $(\SSS, \om_\lambda)$ and $(\CCC, \om_\lambda)$. The other half are their upside-downs with appropriate label modification. 
Here $r=2k$ or $2k+1$, and $\eps_r \equiv r \pmod 2$. The blowup is performed along an embedded ball of capacity $c$. 
The values of the moment map $\mu$ are normalized by setting its minimum to zero. All fixed surfaces have genus $0$ (hence omitted) and area $A$. The vertex labels are only there for references, not part of the data for a decorated graph.    

\begin{figure}[H]
\begin{subfigure}{0.46\textwidth}
\begin{tikzpicture}
  [
    point/.style = {draw, circle, fill = black, inner sep = 1pt},
    every edge quotes/.style = {font=\tiny, color=black, auto=left, inner sep=1pt},
  ]

  % Nodes
  \node (Q') at (0.5,0.5) [point] {};
  \node[above right, xshift=-2pt, yshift=-3pt] at (Q') {\tiny $Q'$}; 
  \node[left] at (Q') {\tiny $\mu = b - bc$}; 

  \node (Q'') at (0,1.5) [point] {};
  \node[right, xshift=-2pt, yshift=0pt] at (Q'') {\tiny $Q''$}; 
  \node[left] at (Q'') {\tiny $\mu = b + ac$}; 

  \node (R') at (0,2.5) [point] {};
  \node[below right, xshift=-2pt, yshift=3pt] at (R') {\tiny $R'$}; 
  \node[left] at (R') {\tiny $\mu = b + a(\lambda - k - \epsilon_r)$};

  \node (S') at (2,3.5) [point] {};
  \node[below left, xshift=3pt, yshift=0pt] at (S') {\tiny $S'$}; 
  \node[above] at (S') {\tiny $\mu = a(\lambda+k)$};

  \node (P') at (2,0) [point] {};
  \node[above left, xshift=3pt, yshift=0pt] at (P') {\tiny $P'$}; 
  \node[below] at (P') {\tiny $\mu = 0$};

  % Edges
  \draw (Q')  edge ["$a + b$"]  (Q'');
  \draw (Q'') edge ["$a$"]      (R');
  \draw (R')  edge ["$ar - b$"] (S');
  \draw (S')   edge ["$a$"]      (P');
  \draw (P')   edge ["$b$"]      (Q');

\end{tikzpicture}	
\caption{\scriptsize}
\end{subfigure}
%%%%%%%%%%%%%%%%%%
\begin{subfigure}{0.46\textwidth}
\begin{tikzpicture}
	[
		point/.style = {draw, circle,  fill = black, inner sep = 1pt},
		every edge quotes/.style = {font=\tiny, color = black, auto=left, inner sep = 1pt},
	]
		\node (Q') at (0,1) [point]{};
	    \node[above right, xshift=-2pt, yshift=-3pt] at (Q') {\tiny $Q'$}; 
	    \node[left] at (Q') {\tiny $\mu = b$};

		\node (R') at (0,2) [point]{};
	    \node[right, xshift=-2pt, yshift=0pt] at (R') {\tiny $R'$} ; 
	    \node[left] at (R') {\tiny $\mu = b + a(\lambda-k-\epsilon_r-c)$};

		\node (R'') at (0.5,3) [point]{};
	    \node[below right, xshift=-2pt, yshift=3pt] at (R'') {\tiny $R''$}; 
	    \node[left] at (R'') {\shortstack{\tiny $\mu = b(1-c) + a(\lambda - k - \epsilon_r+r c)$}};

		\node (S') at (2,3.5) [point]{};
	    \node[below left, xshift=3pt, yshift=0pt] at (S') {\tiny $S'$}; 
	    \node[above] at (S') {\tiny $\mu = a(\lambda+k)$};

		\node (P') at (2,0) [point]{};
	    \node[above left, xshift=3pt, yshift=0pt] at (P') {\tiny $P'$}; 
	    \node[below] at (P') {\tiny $\mu = 0$};

		\draw (Q') edge ["$a$"] (R');
		\draw (R') edge ["$a(r+1)-b$"] (R'');
		\draw (R'') edge ["$ar-b$"] (S');
		\draw (S') edge ["$a$"] (P');
		\draw (P') edge ["$b$"] (Q'); 	
\end{tikzpicture}
\caption{\scriptsize}
\end{subfigure}
%%%%%%%%%%%%%%%%% 
\begin{subfigure}{0.46\textwidth}
\begin{tikzpicture}
	[
		point/.style = {draw, circle,  fill = black, inner sep = 1pt},
		every edge quotes/.style = {font=\tiny, color = black, auto=left, inner sep = 1pt},
	]
		\node (Q') at (0,1) [point]{};
	    \node[above right, xshift=-2pt, yshift=-3pt] at (Q') {\tiny $Q'$}; 
	    \node[left] at (Q') {\tiny $\mu = b$};

		\node (R') at (0,2) [point]{};
	    \node[right, xshift=-2pt, yshift=0pt] at (R') {\tiny $R'$} ; 
	    \node[left] at (R') {\tiny $\mu = b+a(\lambda-k-\epsilon_r)$};

		\node (S') at (1,3) [point]{};
	    \node[below right, xshift=-2pt, yshift=3pt] at (S') {\tiny $S'$}; 
	    \node[left] at (S') {\tiny $\mu = bc+ a(\lambda+k - rc)$};

		\node (S'') at (2,3.5) [point]{};
	    \node[below left, xshift=3pt, yshift=0pt] at (S'') {\tiny $S''$}; 
	    \node[above] at (S'') {\tiny $\mu = a(\lambda+k-c)$};

		\node (P') at (2,0) [point]{};
	    \node[above left, xshift=3pt, yshift=0pt] at (P') {\tiny $P'$}; 
	    \node[below] at (P') {\tiny $\mu = 0$};

		\draw (Q') edge ["$a$"] (R');
		\draw (R') edge ["$ar-b$"] (S');
		\draw (S') edge ["$a(r-1)-b$"] (S'');
		\draw (S'') edge ["$a$"] (P');
		\draw (P') edge ["$b$"] (Q');  	
\end{tikzpicture}		
\caption{\scriptsize $ar-b>a$} 
\end{subfigure}
%%%%%%%%%%%%%%%%%
\begin{subfigure}{0.46\textwidth}
\begin{tikzpicture}
	[
		point/.style = {draw, circle,  fill = black, inner sep = 1pt},
		every edge quotes/.style = {font=\tiny, color = black, auto=left, inner sep = 1pt},
	]
		\node (Q') at (0,1) [point]{};
	    \node[above right, xshift=-2pt, yshift=-3pt] at (Q') {\tiny $Q'$}; 
	    \node[left] at (Q') {\tiny $\mu = b$};

		\node (R') at (0,2) [point]{};
	    \node[right, xshift=-2pt, yshift=0pt] at (R') {\tiny $R'$} ; 
	    \node[left] at (R') {\tiny $\mu = b+a(\lambda-k-\epsilon_r)$};

		\node (S') at (1,3.5) [point]{};
	    \node[below, xshift=0pt, yshift=0pt] at (S') {\tiny $S'$}; 
	    \node[left] at (S') {\tiny $\mu = bc+ a(\lambda+k - rc)$};

		\node (S'') at (2,3) [point]{};
	    \node[below left, xshift=3pt, yshift=0pt] at (S'') {\tiny $S''$}; 
	    \node[right] at (S'') {\tiny $\mu = a(\lambda+k-c)$};

		\node (P') at (2,0) [point]{};
	    \node[above left, xshift=3pt, yshift=0pt] at (P') {\tiny $P'$}; 
	    \node[right] at (P') {\tiny $\mu = 0$};

		\draw (Q') edge ["$a$"] (R');
		\draw (R') edge ["$ar-b$"] (S');
		\draw (S') edge ["$a(1-r)+b$"] (S'');
		\draw (S'') edge ["$a$"] (P');
		\draw (P') edge ["$b$"] (Q');  	
\end{tikzpicture}			
\caption{\scriptsize $a>ar-b$}		
\end{subfigure}	
%%%%%%%%%%%%%%%%%
\begin{subfigure}{0.46\textwidth} 
	\centering
\begin{tikzpicture}
	[
		point/.style = {draw, circle,  fill = black, inner sep = 1pt},
		every edge quotes/.style = {font=\tiny, color = black, auto=left, inner sep = 1pt},
	]
		\node (Q') at (-0.5,1) [point]{};
	    \node[above right, xshift=-2pt, yshift=-3pt] at (Q') {\tiny $Q'$}; 
	    \node[left] at (Q') {\tiny $\mu = r-1$};

		\node (R') at (-0.5,2) [point]{};
	    \node[right, xshift=-2pt, yshift=0pt] at (R') {\tiny $R'$} ; 
	    \node[left] at (R') {\tiny $\mu = 1+\lambda-k-\epsilon_r$};

		\fill (0,3) ellipse [x radius=20pt, y radius=4pt] node[above=5pt]{\tiny $A=c$, $\mu=\lambda+k-c$};

		\node (P') at (0,0) [point]{};
	    \node[above, xshift=3pt, yshift=0pt] at (P') {\tiny $P'$}; 
	    \node[right] at (P') {\tiny $\mu = 0$};

		\draw (P') edge ["$r-1$"] (Q');  	
\end{tikzpicture}
\caption{\scriptsize $a=1$, $b=r-1>0$}\label{fig:a=1|b=r-1} 			
\end{subfigure}	
%%%%%%%%%%%%%%%%%
\begin{subfigure}{0.46\textwidth} 
	\centering
\begin{tikzpicture}
	[
		point/.style = {draw, circle,  fill = black, inner sep = 1pt},
		every edge quotes/.style = {font=\tiny, color = black, auto=left, inner sep = 1pt},
	]
		\node (Q') at (0,1) [point]{};
	    \node[above right, xshift=-2pt, yshift=-3pt] at (Q') {\tiny $Q'$}; 
	    \node[left] at (Q') {\tiny $\mu = b-bc$};

		\node (R') at (0,3) [point]{};
	    \node[right, xshift=-2pt, yshift=0pt] at (R') {\tiny $R'$} ; 
	    \node[left] at (R') {\tiny $\mu = b-bc+a(\lambda-k-\epsilon_r)$};

		\node (S') at (1,4) [point]{};
	    \node[below, xshift=-2pt, yshift=3pt] at (S') {\tiny $S'$}; 
	    \node[above] at (S') {\tiny $\mu = a(\lambda+k)-bc$};

		\node (P') at (2,1) [point]{};
	    \node[left, xshift=3pt, yshift=0pt] at (P') {\tiny $P'$}; 
	    \node[right] at (P') {\tiny $\mu = ac - bc$};

		\node (P'') at (1,0.5) [point]{};
	    \node[above, xshift=3pt, yshift=0pt] at (P'') {\tiny $P''$}; 
	    \node[below] at (P'') {\tiny $\mu = 0$};

		\draw (Q') edge ["$a$"] (R');
		\draw (R') edge ["$ar-b$"] (S');
		\draw (S') edge ["$a$"] (P');
		\draw (P') edge ["$a-b$"] (P'');
		\draw (P'') edge ["$b$"] (Q');  	
\end{tikzpicture}			
\caption{\scriptsize $a > b > 0$, $r>0$}\label{fig:a>b} 	
\end{subfigure}
%%%%%%%%%%%%%%%%%
\begin{subfigure}{0.46\textwidth} 
	\centering
\begin{tikzpicture}
	[
		point/.style = {draw, circle,  fill = black, inner sep = 1pt},
		every edge quotes/.style = {font=\tiny, color = black, auto=left, inner sep = 1pt},
	]
		\node (Q') at (0,1.5) [point]{};
	    \node[above right, xshift=-2pt, yshift=-3pt] at (Q') {\tiny $Q'$}; 
	    \node[left] at (Q') {\tiny $\mu = b-ac$};

		\node (R') at (0,2.5) [point]{};
	    \node[right, xshift=-2pt, yshift=0pt] at (R') {\tiny $R'$} ; 
	    \node[left] at (R') {\tiny $\mu = b+a(\lambda-k-\epsilon_r-c)$};

		\node (S') at (2,3) [point]{};
	    \node[below, xshift=-2pt, yshift=3pt] at (S') {\tiny $S'$}; 
	    \node[above] at (S') {\tiny $\mu = a(\lambda+k-c)$};

		\node (P') at (2,0) [point]{};
	    \node[above left, xshift=3pt, yshift=0pt] at (P') {\tiny $P'$}; 
	    \node[right] at (P') {\tiny $\mu = 0$};

		\node (P'') at (1,0.5) [point]{};
	    \node[above, xshift=3pt, yshift=0pt] at (P'') {\tiny $P''$}; 
	    \node[left] at (P'') {\tiny $\mu = bc-ac$};

		\draw (Q') edge ["$a$"] (R');
		\draw (R') edge ["$ar-b$"] (S');
		\draw (S') edge ["$a$"] (P');
		\draw (P') edge ["$b-a$"] (P'');
		\draw (P'') edge ["$b$"] (Q');  	
\end{tikzpicture}			
\caption{\scriptsize $ar> b>a>0$}\label{fig:b>a} 		
\end{subfigure}	
%%%%%%%%%%%%%%%%%
\begin{subfigure}{0.46\textwidth} 
	\centering
\begin{tikzpicture}
	[
		point/.style = {draw, circle,  fill = black, inner sep = 1pt},
		every edge quotes/.style = {font=\tiny, color = black, auto=left, inner sep = 1pt},
	]
		\node (Q') at (-0.5,1) [point]{};
	    \node[above right, xshift=-2pt, yshift=-3pt] at (Q') {\tiny $Q'$}; 
	    \node[left] at (Q') {\tiny $\mu = b-c$};

		\node (R') at (-0.5,2) [point]{};
	    \node[right, xshift=-2pt, yshift=0pt] at (R') {\tiny $R'$} ; 
	    \node[left] at (R') {\tiny $\mu = b+a(\lambda-k-\epsilon_r)-c$};

		\node (S') at (0,3) [point]{};
	    \node[below right, xshift=-2pt, yshift=3pt] at (S') {\tiny $S'$}; 
	    \node[left] at (S') {\tiny $\mu = a(\lambda+k)-c$};

		\fill (0,0) ellipse [x radius=20pt, y radius=4pt] node[below=5pt]{\tiny $A=c$, $\mu=0$};

		\draw (R') edge ["$r-1$"] (S');  	
\end{tikzpicture}
\caption{\scriptsize $a=b=1$, $r>1$}\label{fig:a=1|b=1} 			
\end{subfigure}	

\caption{When $a>0$, $b>0$, and $ar-b>0$}\label{fig:tikz-both}
\end{figure}

\begin{figure}[H] 

\begin{subfigure}{0.46\textwidth} 
	\centering
\begin{tikzpicture}
	[
		point/.style = {draw, circle,  fill = black, inner sep = 1pt},
		every edge quotes/.style = {font=\tiny, color = black, auto=left, inner sep = 1pt},
	]
		\node (Q') at (0,1) [point]{};
	    \node[above right, xshift=-2pt, yshift=-3pt] at (Q') {\tiny $Q'$}; 
	    \node[left] at (Q') {\tiny $\mu =b - bc$};

		\node (Q'') at (0,2) [point]{};
	    \node[right, xshift=-2pt, yshift=0pt] at (Q'') {\tiny $Q''$} ; 
	    \node[left] at (Q'') {\tiny $\mu=b+ac$};

		\node (R') at (1,4) [point]{};
	    \node[below right, xshift=-2pt, yshift=3pt] at (R') {\tiny $R'$}; 
	    \node[left] at (R') {\tiny $\mu = b+a(\lambda-k-\epsilon_r)$};

		\node (S') at (2,1.5) [point]{};
	    \node[below left, xshift=3pt, yshift=0pt] at (S') {\tiny $S'$}; 
	    \node[right] at (S') {\tiny $\mu = a(\lambda+k)$};

		\node (P') at (1,0) [point]{};
	    \node[above left, xshift=3pt, yshift=0pt] at (P') {\tiny $P'$}; 
	    \node[right] at (P') {\tiny $\mu = 0$};

		\draw (Q') edge ["$a+b$"] (Q'');
		\draw (Q'') edge ["$a$"] (R');
		\draw (R') edge ["$b-ar$"] (S');
		\draw (S') edge ["$a$"] (P');
		\draw (P') edge ["$b$"] (Q');  	
\end{tikzpicture}			
\caption{\scriptsize } 	
\end{subfigure}
%%%%%%%%%%%%%%%%%
\begin{subfigure}{0.46\textwidth} 
	\centering
\begin{tikzpicture}
	[
		point/.style = {draw, circle,  fill = black, inner sep = 1pt},
		every edge quotes/.style = {font=\tiny, color = black, auto=left, inner sep = 1pt},
	]
		\node (Q') at (0,2) [point]{};
	    \node[above right, xshift=-2pt, yshift=-3pt] at (Q') {\tiny $Q'$}; 
	    \node[left] at (Q') {\tiny $\mu = b$};

		\node (R') at (1,4) [point]{};
	    \node[below right, xshift=-2pt, yshift=3pt] at (R') {\tiny $R'$}; 
	    \node[left] at (R') {\tiny $\mu = b+a(\lambda-k-\epsilon_r)$};

		\node (S') at (2,3) [point]{};
	    \node[below left, xshift=3pt, yshift=0pt] at (S') {\tiny $S'$}; 
	    \node[right] at (S') {\tiny $\mu = bc+a(\lambda+k-rc)$};

		\node (S'') at (2,1) [point]{};
		\node[left, xshift=-2pt, yshift=0pt] at (S'') {\tiny $S''$} ; 
		\node[right] at (S'') {\tiny $\mu = a(\lambda+k-c)$};

		\node (P') at (1,0) [point]{};
	    \node[above left, xshift=3pt, yshift=0pt] at (P') {\tiny $P'$}; 
	    \node[right] at (P') {\tiny $\mu = 0$};

		\draw (Q') edge ["$a$"] (R');
		\draw (R') edge ["$b-ar$"] (S');
		\draw (S') edge ["$b-ar+a$"] (S'');		
		\draw (S'') edge ["$a$"] (P');
		\draw (P') edge ["$b$"] (Q');  	
\end{tikzpicture}			
\caption{\scriptsize }	
\end{subfigure}
%%%%%%%%%%%%%%%%
\begin{subfigure}{0.46\textwidth} 
	\centering
\begin{tikzpicture}
	[
		point/.style = {draw, circle,  fill = black, inner sep = 1pt},
		every edge quotes/.style = {font=\tiny, color = black, auto=left, inner sep = 1pt},
	]
		\node (Q') at (0,1) [point]{};
	    \node[above right, xshift=-2pt, yshift=-3pt] at (Q') {\tiny $Q'$}; 
	    \node[left] at (Q') {\tiny $\mu = b$};

		\node (R') at (0.5,3) [point]{};
	    \node[below right, xshift=-2pt, yshift=3pt] at (R') {\tiny $R'$}; 
	    \node[left] at (R') {\tiny $\mu = b+a(\lambda-k-\epsilon_r-c)$};		

		\node (R'') at (1.5,4) [point]{};
		\node[below, xshift=-2pt, yshift=0pt] at (R'') {\tiny $R''$} ; 
		\node[left] at (R'') {\tiny $\mu = b-bc+a(\lambda-k-\epsilon_r+rc)$};		

		\node (S') at (2,2) [point]{};
	    \node[below left, xshift=3pt, yshift=0pt] at (S') {\tiny $S'$}; 
	    \node[right] at (S') {\tiny $\mu = a(\lambda+k)$};

		\node (P') at (1,0) [point]{};
	    \node[above left, xshift=3pt, yshift=0pt] at (P') {\tiny $P'$}; 
	    \node[right] at (P') {\tiny $\mu = 0$};

		\draw (Q') edge ["$a$"] (R');
		\draw (R') edge ["$a-b+ar$"] (R'');				
		\draw (R'') edge ["$b-ar$"] (S');
		\draw (S') edge ["$a$"] (P');
		\draw (P') edge ["$b$"] (Q');  	
\end{tikzpicture}			
\caption{\scriptsize $a > b-ar > 0$, $b>0$}\label{fig:a>b-ar} 	
\end{subfigure}
%%%%%%%%%%%%%%%
\begin{subfigure}{0.46\textwidth} 
	\centering
\begin{tikzpicture}
	[
		point/.style = {draw, circle,  fill = black, inner sep = 1pt},
		every edge quotes/.style = {font=\tiny, color = black, auto=left, inner sep = 1pt},
	]
		\node (Q') at (0,1) [point]{};
	    \node[above right, xshift=-2pt, yshift=-3pt] at (Q') {\tiny $Q'$}; 
	    \node[left] at (Q') {\tiny $\mu = b$};

		\node (R') at (0.5,4) [point]{};
	    \node[below right, xshift=-2pt, yshift=3pt] at (R') {\tiny $R'$}; 
	    \node[above] at (R') {\tiny $\mu = b+a(\lambda-k-\epsilon_r-c)$};	

		\node (R'') at (1.5,3) [point]{};
		\node[left, xshift=-2pt, yshift=0pt] at (R'') {\tiny $R''$} ; 
		\node[right] at (R'') {\tiny $\mu = b-bc+a(\lambda-k-\epsilon_r+rc)$};		

		\node (S') at (2,2) [point]{};
	    \node[below left, xshift=3pt, yshift=0pt] at (S') {\tiny $S'$}; 
	    \node[right] at (S') {\tiny $\mu = a(\lambda+k)$};

		\node (P') at (1,0) [point]{};
	    \node[above left, xshift=3pt, yshift=0pt] at (P') {\tiny $P'$}; 
	    \node[right] at (P') {\tiny $\mu = 0$};

		\draw (Q') edge ["$a$"] (R');
		\draw (R') edge ["$b-ar-a$"] (R'');				
		\draw (R'') edge ["$b-ar$"] (S');
		\draw (S') edge ["$a$"] (P');
		\draw (P') edge ["$b$"] (Q');  	
\end{tikzpicture}			
\caption{\scriptsize $b-ar > a > 0$, $b>0$}\label{fig:b-ar>a}  	
\end{subfigure}
%%%%%%%%%%%%%%%%%
\begin{subfigure}{0.46\textwidth} 
	\centering
\begin{tikzpicture}
	[
		point/.style = {draw, circle,  fill = black, inner sep = 1pt},
		every edge quotes/.style = {font=\tiny, color = black, auto=left, inner sep = 1pt},
	]
		\node (Q') at (-0.5,1) [point]{};
	    \node[above right, xshift=-2pt, yshift=-3pt] at (Q') {\tiny $Q'$}; 
	    \node[left] at (Q') {\tiny $\mu = b$};

		\node (S') at (0.5,2) [point]{};
	    \node[left, xshift=-2pt, yshift=0pt] at (S') {\tiny $S'$} ; 
	    \node[right] at (S') {\tiny $\mu = a(\lambda+k)$};

		\fill (0,3) ellipse [x radius=20pt, y radius=4pt] node[above=5pt]{\tiny $A=c$, $\mu = b+a(\lambda-k-\epsilon_r)-c$};

		\node (P') at (0,0) [point]{};
	    \node[above, xshift=3pt, yshift=0pt] at (P') {\tiny $P'$}; 
	    \node[right] at (P') {\tiny $\mu = 0$};

		\draw (P') edge ["$r+1$"] (Q');  	
\end{tikzpicture}
\caption{\scriptsize $a=1$, $b=r+1$}\label{fig:a=1|b=r+1} 			
\end{subfigure}
%%%%%%%%%%%%%%%%%%%
\begin{subfigure}{0.46\textwidth} 
	\centering
\begin{tikzpicture}
	[
		point/.style = {draw, circle,  fill = black, inner sep = 1pt},
		every edge quotes/.style = {font=\tiny, color = black, auto=left, inner sep = 1pt},
	]
		\node (P'') at (0,1) [point]{};
	    \node[above right, xshift=-2pt, yshift=-3pt] at (P'') {\tiny $P''$}; 
	    \node[left] at (P'') {\tiny $\mu = bc-ac$};

		\node (Q') at (0,2) [point]{};
	    \node[right, xshift=-2pt, yshift=0pt] at (Q') {\tiny $Q'$} ; 
	    \node[left] at (Q') {\tiny $\mu = b-ac$};

		\node (R') at (1,4) [point]{};
	    \node[below right, xshift=-2pt, yshift=3pt] at (R') {\tiny $R'$}; 
	    \node[left] at (R') {\tiny $\mu = b+a(\lambda-k-\epsilon_r-c)$};

		\node (S') at (2,1.5) [point]{};
	    \node[below left, xshift=3pt, yshift=0pt] at (S') {\tiny $S'$}; 
	    \node[right] at (S') {\tiny $\mu = a(\lambda+k-c)$};

		\node (P') at (1,0) [point]{};
	    \node[above left, xshift=3pt, yshift=0pt] at (P') {\tiny $P'$}; 
	    \node[right] at (P') {\tiny $\mu = 0$};

		\draw (P'') edge ["$b$"] (Q');  			
		\draw (Q'') edge ["$a$"] (R');
		\draw (R') edge ["$b-ar$"] (S');
		\draw (S') edge ["$a$"] (P');
		\draw (P') edge ["$b-a$"] (P'');		
\end{tikzpicture}			
\caption{\scriptsize $b> a > 0$, $b > ar$}\label{fig:b>a>0} 	
\end{subfigure}
%%%%%%%%%%%%%%
\begin{subfigure}{0.46\textwidth} 
	\centering
\begin{tikzpicture}
	[
		point/.style = {draw, circle,  fill = black, inner sep = 1pt},
		every edge quotes/.style = {font=\tiny, color = black, auto=left, inner sep = 1pt},
	]
		\node (Q') at (0,2) [point]{};
	    \node[above right, xshift=-2pt, yshift=-3pt] at (Q') {\tiny $Q'$}; 
	    \node[left] at (Q') {\tiny $\mu = b-bc$};

		\node (R') at (1,4) [point]{};
	    \node[below right, xshift=-2pt, yshift=3pt] at (R') {\tiny $R'$}; 
	    \node[left] at (R') {\tiny $\mu = b-bc+a(\lambda-k-\epsilon_r)$};

		\node (S') at (2,3) [point]{};
	    \node[below left, xshift=3pt, yshift=0pt] at (S') {\tiny $S'$}; 
	    \node[right] at (S') {\tiny $\mu = a(\lambda+k)-bc$};

		\node (P') at (2,1) [point]{};
		\node[left, xshift=-2pt, yshift=0pt] at (P') {\tiny $P'$} ; 
		\node[right] at (P') {\tiny $\mu = ac-bc$};

		\node (P'') at (1,0) [point]{};
	    \node[above left, xshift=3pt, yshift=0pt] at (P'') {\tiny $P''$}; 
	    \node[right] at (P'') {\tiny $\mu = 0$};

		\draw (Q') edge ["$a$"] (R');
		\draw (R') edge ["$b$"] (S');
		\draw (S') edge ["$a$"] (P');		
		\draw (P') edge ["$a-b$"] (P'');
		\draw (P'') edge ["$b$"] (Q');  	
\end{tikzpicture}			
\caption{\scriptsize $a > b> 0$, $r=0$}\label{fig:a>b-B} 	
\end{subfigure}
%%%%%%%%%%%%%%%%%%%
\begin{subfigure}{0.46\textwidth} 
	\centering
\begin{tikzpicture}
	[
		point/.style = {draw, circle,  fill = black, inner sep = 1pt},
		every edge quotes/.style = {font=\tiny, color = black, auto=left, inner sep = 1pt},
	]
		\node (Q') at (-0.5,1) [point]{};
	    \node[above right, xshift=-2pt, yshift=-3pt] at (Q') {\tiny $Q'$}; 
	    \node[left] at (Q') {\tiny $\mu = b-c$};

		\node (R') at (0,3) [point]{};
	    \node[right, xshift=-2pt, yshift=0pt] at (R') {\tiny $R'$} ; 
	    \node[left] at (R') {\tiny $\mu = b+a(\lambda-k-\epsilon_r)-c$};

		\node (S') at (0.5,2) [point]{};
	    \node[below right, xshift=-2pt, yshift=3pt] at (S') {\tiny $S'$}; 
	    \node[left] at (S') {\tiny $\mu = a(\lambda+k)-c$};

		\fill (0,0) ellipse [x radius=20pt, y radius=4pt] node[below=5pt]{\tiny $A=c$, $\mu=0$};

\end{tikzpicture}
\caption{\scriptsize $a=b=1$, $r=0$}\label{fig:a=1|b=1-B} 			
\end{subfigure}	
\caption{When $a>0$, $b>0$, and $ar-b<0$}
\end{figure}

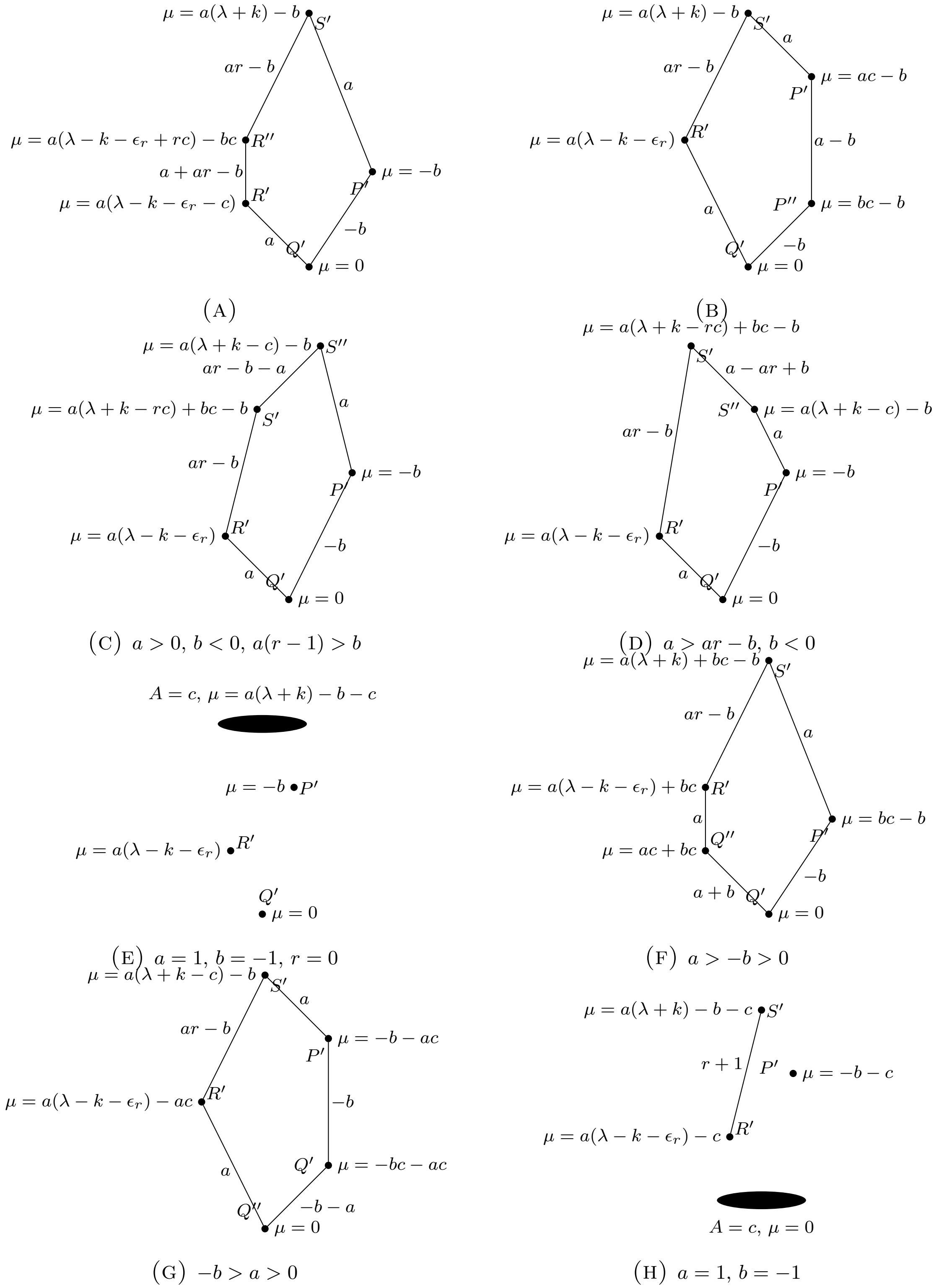
\begin{figure}[H]	

\begin{subfigure}{0.46\textwidth} 
	\centering
\begin{tikzpicture}
	[
		point/.style = {draw, circle,  fill = black, inner sep = 1pt},
		every edge quotes/.style = {font=\tiny, color = black, auto=left, inner sep = 1pt},
	]
		\node (R') at (0,1) [point]{};
	    \node[above right, xshift=-2pt, yshift=-3pt] at (R') {\tiny $R'$}; 
	    \node[left] at (R') {\tiny $\mu = a(\lambda-k-\epsilon_r-c)$};

		\node (R'') at (0,2) [point]{};
	    \node[right, xshift=-2pt, yshift=0pt] at (R'') {\tiny $R''$} ; 
	    \node[left] at (R'') {\tiny $\mu = a(\lambda-k-\epsilon_r+rc)-bc$};		

		\node (S') at (1,4) [point]{};
	    \node[below right, xshift=-2pt, yshift=3pt] at (S') {\tiny $S'$}; 
	    \node[left] at (S') {\tiny $\mu = a(\lambda+k)-b$};

		\node (P') at (2,1.5) [point]{};
	    \node[below left, xshift=3pt, yshift=0pt] at (P') {\tiny $P'$}; 
	    \node[right] at (P') {\tiny $\mu = -b$};

		\node (Q') at (1,0) [point]{};
	    \node[above left, xshift=3pt, yshift=0pt] at (Q') {\tiny $Q'$}; 
	    \node[right] at (Q') {\tiny $\mu = 0$};

		\draw (R') edge ["$a+ar-b$"] (R'');
		\draw (R'') edge ["$ar-b$"] (S');
		\draw (S') edge ["$a$"] (P');
		\draw (P') edge ["$-b$"] (Q');
		\draw (Q') edge ["$a$"] (R');  	
\end{tikzpicture}			
\caption{\scriptsize }
\end{subfigure}
%%%%%%%%%%%%%%%%
\begin{subfigure}{0.46\textwidth} 
	\centering
\begin{tikzpicture}
	[
		point/.style = {draw, circle,  fill = black, inner sep = 1pt},
		every edge quotes/.style = {font=\tiny, color = black, auto=left, inner sep = 1pt},
	]
		\node (R') at (0,2) [point]{};
	    \node[above right, xshift=-2pt, yshift=-3pt] at (R') {\tiny $R'$}; 
	    \node[left] at (R') {\tiny $\mu = a(\lambda-k-\epsilon_r)$};

		\node (S') at (1,4) [point]{};
	    \node[below right, xshift=-2pt, yshift=3pt] at (S') {\tiny $S'$}; 
	    \node[left] at (S') {\tiny $\mu = a(\lambda+k)-b$};

		\node (P') at (2,3) [point]{};
	    \node[below left, xshift=3pt, yshift=0pt] at (P') {\tiny $P'$}; 
	    \node[right] at (P') {\tiny $\mu = ac-b$};

		\node (P'') at (2,1) [point]{};
		\node[left, xshift=-2pt, yshift=0pt] at (P'') {\tiny $P''$} ; 
		\node[right] at (P'') {\tiny $\mu = bc-b$};

		\node (Q') at (1,0) [point]{};
	    \node[above left, xshift=3pt, yshift=0pt] at (Q') {\tiny $Q'$}; 
	    \node[right] at (Q') {\tiny $\mu = 0$};

		\draw (R'') edge ["$ar-b$"] (S');
		\draw (S') edge ["$a$"] (P');
		\draw (P') edge ["$a-b$"] (P'');		
		\draw (P'') edge ["$-b$"] (Q');
		\draw (Q') edge ["$a$"] (R');  	
\end{tikzpicture}			
\caption{\scriptsize }
\end{subfigure}
%%%%%%%%%%%%%%%%
\begin{subfigure}{0.46\textwidth} 
	\centering
\begin{tikzpicture}
	[
		point/.style = {draw, circle,  fill = black, inner sep = 1pt},
		every edge quotes/.style = {font=\tiny, color = black, auto=left, inner sep = 1pt},
	]
		\node (R') at (0,1) [point]{};
	    \node[above right, xshift=-2pt, yshift=-3pt] at (R') {\tiny $R'$}; 
	    \node[left] at (R') {\tiny $\mu = a(\lambda-k-\epsilon_r)$};

		\node (S') at (0.5,3) [point]{};
	    \node[below right, xshift=-2pt, yshift=3pt] at (S') {\tiny $S'$}; 
	    \node[left] at (S') {\tiny $\mu = a(\lambda+k-rc)+bc-b$};

		\node (S'') at (1.5,4) [point]{};
		\node[right, xshift=-2pt, yshift=0pt] at (S'') {\tiny $S''$} ; 
		\node[left] at (S'') {\tiny $\mu = a(\lambda+k-c)-b$};

		\node (P') at (2,2) [point]{};
	    \node[below left, xshift=3pt, yshift=0pt] at (P') {\tiny $P'$}; 
	    \node[right] at (P') {\tiny $\mu = -b$};

		\node (Q') at (1,0) [point]{};
	    \node[above left, xshift=3pt, yshift=0pt] at (Q') {\tiny $Q'$}; 
	    \node[right] at (Q') {\tiny $\mu = 0$};

		\draw (R') edge ["$ar-b$"] (S');
		\draw (S') edge ["$ar-b-a$"] (S'');			
		\draw (S'') edge ["$a$"] (P');
		\draw (P') edge ["$-b$"] (Q');
		\draw (Q') edge ["$a$"] (R');  	
\end{tikzpicture}			
\caption{\scriptsize $a > 0$, $b<0$, $a(r-1) > b$}\label{fig:a<ar-b} 	
\end{subfigure}
%%%%%%%%%%%%%%%%
\begin{subfigure}{0.46\textwidth} 
	\centering
\begin{tikzpicture}
	[
		point/.style = {draw, circle,  fill = black, inner sep = 1pt},
		every edge quotes/.style = {font=\tiny, color = black, auto=left, inner sep = 1pt},
	]
		\node (R') at (0,1) [point]{};
	    \node[above right, xshift=-2pt, yshift=-3pt] at (R') {\tiny $R'$}; 
	    \node[left] at (R') {\tiny $\mu = a(\lambda-k-\epsilon_r)$};

		\node (S') at (0.5,4) [point]{};
	    \node[below right, xshift=-2pt, yshift=3pt] at (S') {\tiny $S'$}; 
	    \node[above] at (S') {\tiny $\mu = a(\lambda+k-rc)+bc-b$};

		\node (S'') at (1.5,3) [point]{};
		\node[left, xshift=-2pt, yshift=0pt] at (S'') {\tiny $S''$} ; 
		\node[right] at (S'') {\tiny $\mu = a(\lambda+k-c)-b$};	

		\node (P') at (2,2) [point]{};
	    \node[below left, xshift=3pt, yshift=0pt] at (P') {\tiny $P'$}; 
	    \node[right] at (P') {\tiny $\mu = -b$};

		\node (Q') at (1,0) [point]{};
	    \node[above left, xshift=3pt, yshift=0pt] at (Q') {\tiny $Q'$}; 
	    \node[right] at (Q') {\tiny $\mu = 0$};

		\draw (R') edge ["$ar-b$"] (S');
		\draw (S') edge ["$a-ar+b$"] (S'');			
		\draw (S'') edge ["$a$"] (P');
		\draw (P') edge ["$-b$"] (Q');
		\draw (Q') edge ["$a$"] (R');  	
\end{tikzpicture}			
\caption{\scriptsize $a > ar-b$, $b<0$}\label{fig:a>ar-b}  	
\end{subfigure}
%%%%%%%%%%%%%%%%
\begin{subfigure}{0.46\textwidth} 
	\centering
\begin{tikzpicture}
	[
		point/.style = {draw, circle,  fill = black, inner sep = 1pt},
		every edge quotes/.style = {font=\tiny, color = black, auto=left, inner sep = 1pt},
	]
		\node (R') at (-0.5,1) [point]{};
	    \node[above right, xshift=-2pt, yshift=-3pt] at (R') {\tiny $R'$}; 
	    \node[left] at (R') {\tiny $\mu = a(\lambda-k-\epsilon_r)$};

		\node (P') at (0.5,2) [point]{};
	    \node[right, xshift=-2pt, yshift=0pt] at (P') {\tiny $P'$} ; 
	    \node[left] at (P') {\tiny $\mu = -b$};

		\fill (0,3) ellipse [x radius=20pt, y radius=4pt] node[above=5pt]{\tiny $A=c$, $\mu = a(\lambda+k)-b-c$};

		\node (Q') at (0,0) [point]{};
	    \node[above, xshift=3pt, yshift=0pt] at (Q') {\tiny $Q'$}; 
	    \node[right] at (Q') {\tiny $\mu = 0$};

\end{tikzpicture}
\caption{\scriptsize $a=1$, $b=-1$, $r=0$}\label{fig:a=1|b=-1} 			
\end{subfigure}
%%%%%%%%%%%%%%%%%
\begin{subfigure}{0.46\textwidth} 
	\centering
\begin{tikzpicture}
	[
		point/.style = {draw, circle,  fill = black, inner sep = 1pt},
		every edge quotes/.style = {font=\tiny, color = black, auto=left, inner sep = 1pt},
	]
		\node (Q'') at (0,1) [point]{};
	    \node[above right, xshift=-2pt, yshift=-3pt] at (Q'') {\tiny $Q''$}; 
	    \node[left] at (Q'') {\tiny $\mu = ac+bc$};

		\node (R') at (0,2) [point]{};
	    \node[right, xshift=-2pt, yshift=0pt] at (R') {\tiny $R'$} ; 
	    \node[left] at (R') {\tiny $\mu = a(\lambda-k-\epsilon_r)+bc$};

		\node (S') at (1,4) [point]{};
	    \node[below right, xshift=-2pt, yshift=3pt] at (S') {\tiny $S'$}; 
	    \node[left] at (S') {\tiny $\mu = a(\lambda+k)+bc-b$};

		\node (P') at (2,1.5) [point]{};
	    \node[below left, xshift=3pt, yshift=0pt] at (P') {\tiny $P'$}; 
	    \node[right] at (P') {\tiny $\mu = bc-b$};

		\node (Q') at (1,0) [point]{};
	    \node[above left, xshift=3pt, yshift=0pt] at (Q') {\tiny $Q'$}; 
	    \node[right] at (Q') {\tiny $\mu = 0$};

		\draw (R') edge ["$ar-b$"] (S');	
		\draw (S') edge ["$a$"] (P');
		\draw (P') edge ["$-b$"] (Q');
		\draw (Q') edge ["$a+b$"] (Q'');				
		\draw (Q'') edge ["$a$"] (R');  		
\end{tikzpicture}			
\caption{\scriptsize $a > -b > 0$}\label{fig:a>-b} 	
\end{subfigure}
%%%%%%%%%%%%%%%%%%%%%%
\begin{subfigure}{0.46\textwidth} 
	\centering
\begin{tikzpicture}
	[
		point/.style = {draw, circle,  fill = black, inner sep = 1pt},
		every edge quotes/.style = {font=\tiny, color = black, auto=left, inner sep = 1pt},
	]
		\node (R') at (0,2) [point]{};
	    \node[above right, xshift=-2pt, yshift=-3pt] at (R') {\tiny $R'$}; 
	    \node[left] at (R') {\tiny $\mu = a(\lambda-k-\epsilon_r)-ac$};

		\node (S') at (1,4) [point]{};
	    \node[below right, xshift=-2pt, yshift=3pt] at (S') {\tiny $S'$}; 
	    \node[left] at (S') {\tiny $\mu = a(\lambda+k-c)-b$};

		\node (P') at (2,3) [point]{};
	    \node[below left, xshift=3pt, yshift=0pt] at (P') {\tiny $P'$}; 
	    \node[right] at (P') {\tiny $\mu = -b-ac$};

		\node (Q') at (2,1) [point]{};
		\node[left, xshift=-2pt, yshift=0pt] at (Q') {\tiny $Q'$} ; 
		\node[right] at (Q') {\tiny $\mu = -bc-ac$};

		\node (Q'') at (1,0) [point]{};
	    \node[above left, xshift=3pt, yshift=0pt] at (Q'') {\tiny $Q''$}; 
	    \node[right] at (Q'') {\tiny $\mu = 0$};

		\draw (R') edge ["$ar-b$"] (S');	
		\draw (S') edge ["$a$"] (P');
		\draw (P') edge ["$-b$"] (Q');
		\draw (Q') edge ["$-b-a$"] (Q'');				
		\draw (Q'') edge ["$a$"] (R');  		
\end{tikzpicture}			
\caption{\scriptsize $-b> a > 0$}\label{fig:a<-bb} 	
\end{subfigure}
%%%%%%%%%%%%%%%%%%%%%
\begin{subfigure}{0.46\textwidth} 
	\centering
\begin{tikzpicture}
	[
		point/.style = {draw, circle,  fill = black, inner sep = 1pt},
		every edge quotes/.style = {font=\tiny, color = black, auto=left, inner sep = 1pt},
	]
		\node (R') at (-0.5,1) [point]{};
	    \node[above right, xshift=-2pt, yshift=-3pt] at (R') {\tiny $R'$}; 
	    \node[left] at (R') {\tiny $\mu = a(\lambda-k-\epsilon_r)-c$};

		\node (S') at (0,3) [point]{};
	    \node[right, xshift=-2pt, yshift=0pt] at (S') {\tiny $S'$} ; 
	    \node[left] at (S') {\tiny $\mu = a(\lambda+k)-b-c$};

		\node (P') at (0.5,2) [point]{};
	    \node[left, xshift=-2pt, yshift=3pt] at (P') {\tiny $P'$}; 
	    \node[right] at (P') {\tiny $\mu = -b-c$};

		\fill (0,0) ellipse [x radius=20pt, y radius=4pt] node[below=5pt]{\tiny $A=c$, $\mu=0$};

		\draw (R') edge ["$r+1$"] (S');  			

\end{tikzpicture}
\caption{\scriptsize $a=1$, $b=-1$}\label{fig:a=1|b=-1-B} 			
\end{subfigure}	
\caption{When $a>0$, $b<0$}
\end{figure}

\begin{figure}[H]

\begin{tikzpicture}
	[
		point/.style = {draw, circle,  fill = black, inner sep = 1pt},
		every edge quotes/.style = {font=\tiny, color = black, auto=left, inner sep = 1pt},
	]
		\node (R') at (-0.5,1) [point]{};
	    \node[above right, xshift=-2pt, yshift=-3pt] at (R') {\tiny $R'$}; 
	    \node[left] at (R') {\tiny $\mu = \lambda-k-\epsilon_r-c$};		

		\node (R'') at (0,2) [point]{};
	    \node[below right, xshift=-2pt, yshift=3pt] at (R'') {\tiny $R''$}; 
	    \node[left] at (R'') {\tiny $\mu = \lambda-k-\epsilon_r+cr$};		

		\node (S') at (0.5,2.5) [point]{};
	    \node[right, xshift=-2pt, yshift=0pt] at (S') {\tiny $S'$} ; 
	    \node[left] at (S') {\tiny $\mu = \lambda+k$};

		\fill (0,0) ellipse [x radius=20pt, y radius=4pt] node[below=5pt]{\tiny $A=1$, $\mu=0$};

		\draw (R'') edge ["$r$"] (S');  
		\draw (R') edge ["$r+1$"] (R'');  					

\end{tikzpicture}
%%%%%%%%%%%%%%%%%%%%
\begin{tikzpicture}
	[
		point/.style = {draw, circle,  fill = black, inner sep = 1pt},
		every edge quotes/.style = {font=\tiny, color = black, auto=left, inner sep = 1pt},
	]
		\node (R') at (-0.5,1) [point]{};
	    \node[above right, xshift=-2pt, yshift=-3pt] at (R') {\tiny $R'$}; 
	    \node[left] at (R') {\tiny $\mu = \lambda-k-\epsilon_r$}; 

		\node (S') at (0,2) [point]{};
	    \node[below right, xshift=-2pt, yshift=3pt] at (S') {\tiny $S'$}; 
	    \node[left] at (S') {\tiny $\mu = \lambda+k-rc$};		

		\node (S'') at (0.5,2.5) [point]{};
	    \node[right, xshift=-2pt, yshift=0pt] at (S'') {\tiny $S''$} ; 
	    \node[left] at (S'') {\tiny $\mu = \lambda+k-c$};

		\fill (0,0) ellipse [x radius=20pt, y radius=4pt] node[below=5pt]{\tiny $A=1$, $\mu=0$};

		\draw (R') edge ["$r$"] (S');  
		\draw (S') edge ["$r-1$"] (S'');  					

\end{tikzpicture}
%%%%%%%%%%%%%%%%%%%%%
\begin{tikzpicture}
	[
		point/.style = {draw, circle,  fill = black, inner sep = 1pt},
		every edge quotes/.style = {font=\tiny, color = black, auto=left, inner sep = 1pt},
	]
		\node (R') at (-0.5,1) [point]{};
	    \node[above right, xshift=-2pt, yshift=-3pt] at (R') {\tiny $R'$}; 
	    \node[left] at (R') {\tiny $\mu = \lambda_k-\epsilon_r$};

		\node (S') at (0,2) [point]{};
	    \node[below right, xshift=-2pt, yshift=3pt] at (S') {\tiny $S'$}; 
	    \node[left] at (S') {\tiny $\mu = \lambda+k$};		

		\node (P') at (0.5,0.5) [point]{};
	    \node[right] at (P') {\tiny $\mu = c$};

		\fill (0,0) ellipse [x radius=20pt, y radius=4pt] node[below=5pt]{\tiny $A=1-c$, $\mu=0$};

		\draw (R') edge ["$r$"] (S');  

\end{tikzpicture}
\caption{$a=1$, $b=0$}
\end{figure} 

\begin{figure}[H]

\begin{tikzpicture}
	[
		point/.style = {draw, circle,  fill = black, inner sep = 1pt},
		every edge quotes/.style = {font=\tiny, color = black, auto=left, inner sep = 1pt},
	]
		\node (P') at (0,0.5) [point]{};
	    \node[left] at (P') {\tiny $\mu = c$};		

		\fill (0,2.5) ellipse [x radius=20pt, y radius=4pt] node[above=5pt]{\tiny $A=\lambda-k-\epsilon_r$, $\mu=1$};

		\fill (0,0) ellipse [x radius=20pt, y radius=4pt] node[below=5pt]{\tiny $A=\lambda+k-c$, $\mu=0$};
\end{tikzpicture}
%%%%%%%%%%%%%%%%%%%%%%%%%%%%%%%%
\begin{tikzpicture}
	[
		point/.style = {draw, circle,  fill = black, inner sep = 1pt},
		every edge quotes/.style = {font=\tiny, color = black, auto=left, inner sep = 1pt},
	]
		\node (Q') at (0,2) [point]{};
	    \node[right] at (Q') {\tiny $\mu = 1-c$};		

		\fill (0,2.5) ellipse [x radius=20pt, y radius=4pt] node[above=5pt]{\tiny $A=\lambda-k-\epsilon_r-c$, $\mu=1$};

		\fill (0,0) ellipse [x radius=20pt, y radius=4pt] node[below=5pt]{\tiny $A=\lambda+k$, $\mu=0$};
\end{tikzpicture}
\caption{$a=0$, $b=1$}
\end{figure} 

\begin{figure}[H]

\begin{tikzpicture}
	[
		point/.style = {draw, circle,  fill = black, inner sep = 1pt},
		every edge quotes/.style = {font=\tiny, color = black, auto=left, inner sep = 1pt},
	]
		\node (Q'') at (-0.5,1.5) [point]{};
	    \node[above right, xshift=-2pt, yshift=-3pt] at (Q'') {\tiny $Q''$}; 
	    \node[left] at (Q'') {\tiny $\mu = r+c$}; 		

		\node (Q') at (0,0.5) [point]{};
	    \node[right, xshift=-2pt, yshift=3pt] at (Q') {\tiny $Q'$}; 
	    \node[left] at (Q') {\tiny $\mu = r(1-c)$};		

		\node (P') at (0.5,0) [point]{};
	    \node[right, xshift=-2pt, yshift=0pt] at (P') {\tiny $P'$} ; 
	    \node[left] at (P') {\tiny $\mu = 0$};

		\fill (0,2.5) ellipse [x radius=20pt, y radius=4pt] node[above=5pt]{\tiny $A=1$, $\mu=\lambda+k$};

		\draw (P') edge ["$r$"] (Q');  
		\draw (Q') edge ["$r+1$"] (Q'');  					

\end{tikzpicture}
%%%%%%%%%%%%%%%%%%
\begin{tikzpicture}
	[
		point/.style = {draw, circle,  fill = black, inner sep = 1pt},
		every edge quotes/.style = {font=\tiny, color = black, auto=left, inner sep = 1pt},
	]
		\node (Q') at (-0.5,1.5) [point]{};
	    \node[above right, xshift=-2pt, yshift=-3pt] at (Q') {\tiny $Q'$}; 
	    \node[left] at (Q') {\tiny $\mu = \lambda-k-\epsilon_r-c$}; 

		\node (P'') at (0,0.5) [point]{};
	    \node[right, xshift=-2pt, yshift=3pt] at (P'') {\tiny $P''$}; 
	    \node[left] at (P'') {\tiny $\mu = (r-1)c$};		

		\node (P') at (0.5,0) [point]{};
	    \node[right, xshift=-2pt, yshift=0pt] at (P') {\tiny $P'$} ; 
	    \node[left] at (P') {\tiny $\mu = 0$};

		\fill (0,2.5) ellipse [x radius=20pt, y radius=4pt] node[above=5pt]{\tiny $A=1$, $\mu=\lambda+k-c$};

		\draw (P'') edge ["$r$"] (Q');  
		\draw (P') edge ["$r-1$"] (P'');  					

\end{tikzpicture}
%%%%%%%%%%%%%%%%%
\begin{tikzpicture}
	[
		point/.style = {draw, circle,  fill = black, inner sep = 1pt},
		every edge quotes/.style = {font=\tiny, color = black, auto=left, inner sep = 1pt},
	]
		\node (Q') at (-0.5,1.5) [point]{};
	    \node[above right, xshift=-2pt, yshift=-3pt] at (Q') {\tiny $Q'$}; 
	    \node[left] at (Q') {\tiny $\mu = r$};

		\node (P') at (0,0.5) [point]{};
	    \node[below right, xshift=-2pt, yshift=3pt] at (P') {\tiny $P'$}; 
	    \node[left] at (P') {\tiny $\mu = 0$};		

		\node (R') at (0.5,2) [point]{};
	    \node[right] at (R') {\tiny $\mu = \lambda+k-c$};

		\fill (0,2.5) ellipse [x radius=20pt, y radius=4pt] node[above=5pt]{\tiny $A=1-c$, $\mu=\lambda+k$};

		\draw (P') edge ["$r$"] (Q');  

\end{tikzpicture}
\caption{$a=1$, $b=r$} 
\end{figure}

%%%%%%%%%%%%%%%%%%%%%%%%%%%%%%%%%%%%%%%%%%%%%%%%%%%%%%%%%%%%%%%%%%%%%%%%%%%%%%%%
\section*{Acknowledgements}
%%%%%%%%%%%%%%%%%%%%%%%%%%%%%%%%%%%%%%%%%%%%%%%%%%%%%%%%%%%%%%%%%%%%%%%%%%%%%%%%

MP is supported by NSERC Discovery Grants RGPIN-2020-06428 and RGPIN-2026-07645.  PC was funded by the EoS grant ``Geometry and Beyond" (EOS project 40007524) and partially funded by the EPSRC Grant EP/W015749/1. We would like to thank Yael Karshon for helpful discussions. 

%%%%%%%%%%%%%%%%%%%%%%%%%%%%%%%%%%%%%%%%%%%%%%%%%%%%%%%%%%%%%%%%%%%%%%%%%%%%%%%%
\section*{Use of artificial intelligence}
%%%%%%%%%%%%%%%%%%%%%%%%%%%%%%%%%%%%%%%%%%%%%%%%%%%%%%%%%%%%%%%%%%%%%%%%%%%%%%%%

In accordance with the EMS \emph{Code of Practice for Mathematical Publication}
(approved 19 March 2025), which asks authors to acknowledge the use of AI or other automated tools employed in the creation of content by a declaration describing the tools and the nature of their use, and with the corresponding AMS policy adapted from the COPE position statement on authorship and AI tools, we disclose the following.

The mathematical ideas and proofs were generated and mostly written by humans (with the first named author giving several talks about the results since 2023). An AI assistant (Anthropic's Claude) was used during the preparation of this paper in the following ways: as an aid in \LaTeX{} editing and typesetting; for cross-referencing and version control; for checking the references, and in polishing the presentation. Mathematically, it was used to correct a few gaps in the original manuscript, to uniformize notation, and to streamline some of the arguments. The authors take full responsibility for the mathematical content -- the choice of problems, the arguments, and their exposition --  and for the final form of this paper.

%%%%%%%%%%%%%%%%%%%%%%%%%%%%%%%%%%%%%%%%%%%%%%%%%%%%%%%%%%%%%%%%%%%%%%%%%%%%%%%%
\bibliographystyle{amsalpha}
\bibliography{bibliography}
%%%%%%%%%%%%%%%%%%%%%%%%%%%%%%%%%%%%%%%%%%%%%%%%%%%%%%%%%%%%%%%%%%%%%%%%%%%%%%%%

\end{document}